\documentclass{article}
\usepackage[utf8]{inputenc} % allow utf-8 input
\usepackage[T1]{fontenc}    % use 8-bit T1 fonts
\usepackage[colorlinks,urlcolor=blue,citecolor=blue,linkcolor=blue]{hyperref}       % hyperlinks
\usepackage{url}            % simple URL typesetting
\usepackage{booktabs}       % professional-quality tables
\usepackage{amsfonts}       % blackboard math symbols
\usepackage{microtype}      % microtypography

\usepackage{nicefrac} 
\usepackage{bm}
\usepackage{amsmath}
\usepackage{amsthm}
\usepackage{amssymb}
\usepackage{thmtools,thm-restate}
\usepackage{thm-restate}
\usepackage{paralist}
\usepackage{comment}
\usepackage{bbm}
\usepackage{float}
\usepackage{balance}
\usepackage{stmaryrd}
\usepackage{multicol}
\usepackage[export]{adjustbox}
\usepackage{multirow}
\usepackage{tabularx}
\usepackage{mathtools}
\usepackage{mdframed}

\usepackage{enumitem}
\setitemize{noitemsep,topsep=0pt,parsep=0pt,partopsep=0pt}

\usepackage{chemmacros}

\newcommand{\diag}{\operatorname{diag}}

\newcommand{\innp}[1]{\left\langle #1 \right\rangle}

\newcommand{\vx}{\mathbf{x}}

\newcommand{\cx}{\mathcal{X}}

\newcommand{\vy}{\mathbf{y}}

\newcommand{\defeq}{\stackrel{\mathrm{\scriptscriptstyle def}}{=}}

\newcommand{\rr}{\mathbb{R}}
\newcommand{\norm}[1]{\left\| #1 \right\|}

\newcommand{\abs}[1]{\lvert #1 \rvert}

\newcommand*{\vsepfbox}[1]{%
  \begingroup
    \sbox0{\fbox{#1}}%
    \setlength{\fboxrule}{0pt}%
    \mbox{\kern-\fboxsep\fbox{\unhbox0}\kern-\fboxsep}%
  \endgroup
}
\newcommand*{\vertbar}{\rule[-1ex]{0.5pt}{2.5ex}}
\newcommand*{\horzbar}{\rule[.5ex]{2.5ex}{0.5pt}}

\newenvironment{nscenter}
 {\parskip=0pt\par\nopagebreak\centering}
 {\par\noindent\ignorespacesafterend}

\theoremstyle{plain} \numberwithin{equation}{section}
\newtheorem{theorem}{Theorem}[section]
\numberwithin{theorem}{section}

\newtheorem{proposition}[theorem]{Proposition}

\theoremstyle{definition}
\newtheorem{definition}[theorem]{Definition}

\newtheorem{remark}[theorem]{Remark}

\theoremstyle{plain}

\DeclareMathOperator*{\trace}{trace}
\DeclareMathOperator*{\argmin}{argmin}

\DeclareMathOperator{\vertex}{\mathrm{vert}}
\DeclareMathOperator{\co}{\mathrm{conv}}
\DeclareMathOperator{\vect}{vec}
\graphicspath{{imgs/}}
\makeatletter
\def\mathcolor#1#{\@mathcolor{#1}}
\def\@mathcolor#1#2#3{%
  \protect\leavevmode
  \begingroup
    \color#1{#2}#3%
  \endgroup
}
\makeatother

\newcommand{\hrulealg}[0]{\vspace{1mm} \hrule \vspace{1mm}}

\usepackage{subfig}
\usepackage{graphicx,wrapfig,lipsum}

\usepackage[ruled, linesnumbered]{algorithm2e} % http://www.ctan.org/pkg/algorithm2e
\usepackage{xcolor} % http://www.ctan.org/tex-archive/macros/latex/contrib/xcolor
\usepackage{tikz}
\usepackage{todonotes}
\usetikzlibrary{fit,calc}
\newcommand*{\tikzmk}[1]{\tikz[remember picture,overlay,] \node (#1) {};\ignorespaces}
\newcommand{\boxit}[1]{\tikz[remember picture,overlay]{\node[yshift=3pt,fill=#1,opacity=.25,fit={(A)($(B)+(.90\linewidth,.8\baselineskip)$)}] {};}\ignorespaces}

\newcommand{\boxitshort}[1]{\tikz[remember picture,overlay]{\node[yshift=3pt,fill=#1,opacity=.25,fit={(C)($(D)+(.87\linewidth,.8\baselineskip)$)}] {};}\ignorespaces}

\colorlet{pink}{red!40}
\colorlet{lightblue}{blue!30}
\colorlet{lightgreen}{green!30}

\DontPrintSemicolon

\usepackage{jmlr2e}
\usepackage[margin=1in]{geometry}

\usepackage{natbib}
\setcitestyle{authoryear,round,citesep={;},aysep={,},yysep={;}}
\renewcommand{\cite}[1]{\citep{#1}}

\usepackage[vvarbb]{newtxmath}

\begin{document}

\title{CINDy: Conditional gradient-based Identification of Non-linear Dynamics -- Noise-robust recovery}

\author{\name Alejandro Carderera \email  \href{mailto:alejandro.carderera@gatech.edu}{alejandro.carderera@gatech.edu}\\
       \addr Georgia Institute of Technology, USA\\
       Zuse Institute Berlin, Germany 
       \AND
       \name Sebastian Pokutta \email \href{mailto:pokutta@math.tu-berlin.de}{pokutta@math.tu-berlin.de} \\
       \addr Technische Universit\"at Berlin, Germany\\
       Zuse Institute Berlin, Germany
       \AND
       \name Christof Sch\"utte \email \href{mailto:schuette@mi.fu-berlin.de}{schuette@mi.fu-berlin.de} \\
       \addr Freie Universit\"at Berlin, Germany\\
       Zuse Institute Berlin, Germany
       \AND
       \name Martin Weiser \email \href{mailto:weiser@zib.de}{weiser@zib.de} \\
       \addr Zuse Institute Berlin, Germany
       }
       
\maketitle

\begin{abstract}
Governing equations are essential to the study of nonlinear dynamics, often enabling the prediction of previously unseen behaviors as well as the inclusion into control strategies. 
The discovery of governing equations from data thus has the potential to transform data-rich fields where well-established dynamical models remain
unknown. This work contributes to the recent trend in data-driven sparse identification of nonlinear dynamics
of finding the best sparse fit to observational data in a large library of potential nonlinear models. We propose an efficient first-order Conditional Gradient algorithm for solving the underlying optimization problem. In comparison to the most prominent alternative framework, the new framework shows significantly improved performance on several essential issues like sparsity-induction, structure-preservation, noise robustness, and sample efficiency. We demonstrate these advantages on several
dynamics from the field of synchronization, particle dynamics, and enzyme chemistry.
%We demonstrate these advantages on several
%examples with quite different dynamical behaviors from synchronization, particle dynamics, and enzyme chemistry.
\end{abstract}

\section{Introduction}
Many of the developments of physics have stemmed from our ability to
describe natural phenomena in terms of differential equations. These
equations have helped build our understanding of natural phenomena in
fields as wide-ranging as classical mechanics, electromagnetism, fluid
dynamics, neuroscience and quantum mechanics. They have also enabled key
technological advances such as the combustion engine, the laser, or
the transistor.

The modern age of Machine Learning and Big Data has heralded an age of
\emph{data-driven} models, in which the phenomena we explain are
described in terms of statistical relationships and static data. Given
sufficient data, we are able to train neural networks to classify, or
to predict, with high accuracy, without the underlying model having
any apparent knowledge of how the data was generated, or its
structure. This makes the task of classifying, or predicting, on
out-of-sample data a particularly challenging task. On the other hand,
there has been a recent surge in interest in recovering the
differential equations with which the data, often coming from a
physical system, have been generated. This enables us to better
understand how the data is generated, and to better predict on
out-of-sample data, as opposed to using other learning approaches. Moreover, 
learning governing equations also permits understanding the mechanisms underlying the observed dynamical behavior; this is key to further scientific progress.  

The seminal work of \cite{schmidt2009distilling} used symbolic
regression to search the space of mathematical expressions, in order
to find one that adequately fits the data. This entails randomly
combining mathematical operations, analytical functions, state
variables and constants and selecting those that show promise. These
are later randomly expanded and combined in search of an expression
that represents the data sufficiently well. Related to this approach
is the \emph{Approximate Vanishing Ideal} Algorithm
\cite{heldt2009approximate}, based on the combination of Gr\"{o}bner
and Border bases with total least-squares regression, where a set of
polynomials over (arbitrary) basis functions is successively expanded
to capture all relations approximately satisfied by the data. A more
recent algorithm, known as the \emph{Sparse Identification of
  Nonlinear Dynamics} (SINDy) algorithm assumes that we have access to
a library of predefined basis functions, and the problem becomes that
of finding a linear combination of basis functions that best predicts
the data at hand. This is done using sequentially-thresholded
least-squares, in order to recover a sparse linear combination of
basis functions (and potentially the coordinate system) that is able to represent the underlying phenomenon
well \cite{brunton2016discovering,Champion22445}. This algorithm works extremely
well when using noise-free data, but often produces dense solutions
when the data is contaminated with noise. There have been several
suggestions to deal with this, from more noise-robust non-convex
problem formulations \cite{schaeffer2017sparse}, to problem
formulations that involve both learning the dynamic, and the noise
contaminating the underlying data \cite{rudy2019deep,
  kaheman2020automatic}. Neither of these approaches is computationally efficient
for high-dimensional problems. The former having the additional
drawback that the problem formulation is non-convex. The latter, on the other hand, requires solving an optimization problem whose dimension increases linearly with the number of samples in the training data (as it involves learning the noise vector associated with each data point), as opposed to simply increasing linearly with the dimension of the phenomena and the size of the library of basis functions.

%This stands in contrast to dynamic relationships which can be described by differential equations. The difference between these approaches can best be viewed through an analogy - consider the important breakthroughs carried out by Johannes Kepler and Isaac Newton. The former developed a data-driven model of elliptical orbits to describe the trajectories of celestial bodies, using what were at the time state-of-the-art measurements of the position of the planets. However, Kepler did not explain the dynamical relationships that gave rise to these elliptical trajectories, or how celestial bodies reacted when perturbed. On the other hand, Newton developed Newtonian mechanics, a set of dynamical differential equations that relate momentum and energy, that allow us to discover the underlying processes that generate these elliptic orbits. These equations allow us to make predictions on regimes for which we have no data, and have provided us with a robust set of tools with which we have even gone so far as to land a spacecraft on the moon.

\subsection{Contributions}
In this paper we present the \emph{Conditional gradient-based
  Identification of Non-linear Dynamics} framework, dubbed CINDy, in
homage to the influential SINDy framework presented in
\cite{brunton2016discovering}, which uses a sparsity-inducing optimization
algorithm to solve convex formulations of the sparse
recovery problem. CINDy uses a first-order convex optimization algorithm
based on the \emph{Conditional Gradient} (CG) algorithm
\cite{polyak66cg} (also known as the Frank-Wolfe algorithm
\cite{fw56}), and brings together many of the advantages of existing
sparse recovery techniques into a single framework. As documented in detail below, we compared CINDy to the most prominent alternative frameworks for solving the respective learning problem (SINDy, FISTA, IPM, SR3) with the following results:
\begin{enumerate}
    \item \textbf{Sparsity-inducing.} The CG-based optimization algorithm has an
      implicit bias for sparse solutions through the way it builds its
      iterates. Other existing approaches are forced to ensure
      sparsity through thresholding, or through problem formulations
      that encourage sparsity. This has a major impact on the
      \emph{structural generalization behavior}, where CINDy
      significantly outperforms other methods leading to much more accurate
      trajectory predictions in the presence of noise. 
    \item \textbf{Structure-preserving dynamic.} The CINDy framework
      can easily incorporate underlying symmetries and conservation
      laws into the learning problem, resulting in learned dynamics
      consistent with the true physics, with minimal impact on the
      running time of the algorithm but significantly reducing sample
      complexity (due to reduced degrees of freedom) and
      improved generalization performance. 
\item \textbf{Noise robustness.} When it comes to recovery of dynamics in the
  presence of noise, we demonstrate a significant
  advantage over SINDy of about one to two orders of
  magnitude in recovery error \emph{with respect to the true dynamic},
  rather than just out-of-sample errors. This is largely due to the sparsity induced by
  the underlying CG optimization algorithm.
\item \textbf{Sample efficiency and large-scale learning.} We 
demonstrate that, given a certain noise level, CINDy will
  require significantly fewer samples to recover the dynamic with a
  given accuracy. Moreover, being a first-order method our approach
  naturally allows for the learning of large-scale dynamics, allowing
  even the use of stochastic first-order information in the extremely
  large-scale regime.
\item \textbf{Black-box implementation.} We provide an implementation
  of CINDy that can be used as a black-box not requiring any
  specialized knowledge in CG methods. The source code is made
  available under \url{https://github.com/ZIB-IOL}. We hope that this
  stimulates research in the use of CG-based algorithm for sparse
  recovery.
\end{enumerate}

\subsection{Preliminaries}

We denote vectors using bold lower-case letters, and matrices using
upper-case letters. We will use $x_i$ to refer to the $i$-th element
of the vector $\vx$, and $X_{i,j}$ to refer to the element on the
$i$-th row and $j$-th column of the matrix $X$. Let $\norm{\vx}$ and
$\norm{\vx}_1$ denote the $\ell_2$ and $\ell_1$ norm of $\vx$
respectively, furthermore, let $\norm{\vx}_0$ denote the $\ell_0$
norm\footnote{Technically, the \(\ell_0\) norm is not a norm.}, which
is the number of non-zero elements in $\vx$. Moreover, given a matrix
$X \in \rr^{n \times m}$ for $p,q \geq 1$ let
$\norm{X}_{p,q} = \left(\sum_{j=1}^m \left(\sum_{i=1}^n\abs{X_{i,j}}^p
  \right)^{q/p} \right)^{1/q}$ denote the $\ell_{p,q}$ norm of $X$. We
will use $\norm{X}_{F} = \norm{X}_{2,2}$ to refer to the familiar
\emph{Frobenius norm} of a matrix, and $\norm{X}_0$ to refer to the
number of non-zero elements in $X$. Given a matrix
$X\in\rr^{m \times n}$ let $\vect \left( X \right) \in \rr^{mn}$
denote the vectorization of the matrix $X$, that is the stacking
$\vect \left( X \right) =
[X_{1,1},\cdots,X_{m,1},\cdots,X_{1,2},\cdots,X_{m,2},\cdots,X_{n,1},
\cdots,X_{m,n}]^T$.  Given a non-empty set $\mathcal{S}\subset \rr^n$
we refer to its \emph{convex hull} as $\co\left( \mathcal{S}
\right)$. The trace of the square matrix $X\in \rr^{n\times n}$ will
be denoted by $\trace \left( X\right)$. We use $\dot{\vx}(t)$ to
denote the derivative of $\vx(t)$ with respect to time, denoted by
$t$, that is, $\dot{\vx}(t) = \frac{d \vx(t)}{dt}$. Given two integers
$i\in \mathbb{Z}$ and $j \in \mathbb{Z}$ with $i \leq j$ we use
$\llbracket i, j \rrbracket$ to denote the set
$\left\{ k \in\mathbb{Z} \mid i \leq k \leq j \right\}$. The vector
with all entries equal to one is denoted by
$\mathbf{1}_d\in\rr^d$. Lastly, we use $\Delta_d$ to denote the unit
probability simplex of dimension $d$, that is, the set
$\Delta_d = \left\{ \vx \in \rr^d \mid \mathbf{1}_d^T\vx = 1, \vx \geq
  0 \right\}$.
  
Throughout the text we will distinguish between \emph{problem formulation}, that is, the specific form of the mathematical optimization problem we are trying to solve, and the \emph{optimization algorithm} used to solve that problem formulation. Throughout the text we refer to CINDy and SINDy as \emph{frameworks}, which are the result of applying a specific \emph{optimization algorithm} to a particular \emph{problem formulation}. In Section~\ref{Section:LearningSparseDynamics} we will largely focus on the problem formulation, while in Section~\ref{sec:frank-wolfe-algor} we will focus mainly on the optimization algorithm. We try to make this distinction to highlight the fact that:

\begin{center}
\emph{The success of any learning framework is the product of coupling  \\ an appropriate problem formulation to a suitable optimization algorithm.}
\end{center}

\section{Learning sparse dynamics} \label{Section:LearningSparseDynamics}

Many physical systems can be described in terms of ordinary
differential equations of the form
$\dot{\vx}(t) = F\left(\vx(t)\right)$, where $\vx(t) \in \rr^d$ denotes
the state of the system at time $t$ and $F: \rr^d \rightarrow \rr^d$
can usually be expressed as a linear combination of simpler \emph{ansatz
functions} $\psi_i: \rr^d \rightarrow \rr$ belonging to a dictionary
$\mathcal{D} = \left\{\psi_i \mid i \in \llbracket 1, n \rrbracket
\right\}$. This allows us to express the dynamic followed by the
system as
$\dot{\vx}(t) = F\left(\vx(t)\right) = \Xi^T \bm{\psi}(\vx(t))$ where
$\Xi \in \rr^{n \times d}$ is a -- typically sparse -- matrix
$\Xi = \left[\xi_1, \cdots, \xi_d \right]$ formed by column vectors
$\xi_i \in \rr^n$ for $i \in \llbracket 1, d \rrbracket$ and
$\bm{\psi}(\vx(t)) = \left[ \psi_1(\vx(t)), \cdots, \psi_n(\vx(t))
\right]^T \in \rr^{n}$. We can therefore write:

\begin{align}
\dot{\vx}(t) = \begin{bmatrix}
\horzbar & \xi_1 & \horzbar\\
 & \vdots & \\
\horzbar & \xi_d & \horzbar
\end{bmatrix} 
\begin{bmatrix}
\psi_1(\vx(t)) \\
\vdots \\
\psi_n(\vx(t))
\end{bmatrix}.
\end{align}
Alternatively, one could also consider that for any $t\geq t_1$ we can write $\vx(t) = \vx(t_1) + \int_{t_1}^t \dot{\vx}(\tau) d\tau = \vx(t_1) + \int_{t_1}^t \Xi^T \bm{\psi}(\vx(\tau)) d\tau = \vx(t_1) + \Xi^T \int_{t_1}^t \bm{\psi}(\vx(\tau)) d\tau$. In matrix form this results in:
\begin{align}
\vx(t) - \vx(t_1) =  \begin{bmatrix}
\horzbar & \xi_1 & \horzbar\\
 & \vdots & \\
\horzbar & \xi_d & \horzbar
\end{bmatrix} 
 \begin{bmatrix}
\int_{t_1}^{t} \psi_1(\vx(\tau)) d\tau  \\
\vdots \\
\int_{t_1}^{t} \psi_n(\vx(\tau)) d\tau
\end{bmatrix},
\end{align}
 In the absence of noise, if we are given a series of data points from the physical system $\left\{ \vx(t_i), \dot{\vx}(t_i) \right\}_{i=1}^m$, then we know that:
 \begin{align*}
\begin{bmatrix}
\vertbar & & \vertbar\\
\dot{\vx}(t_1) & \cdots & \dot{\vx}(t_m)\\
\vertbar & & \vertbar
\end{bmatrix} = 
\begin{bmatrix}
\horzbar & \xi_1 & \horzbar\\
 & \vdots & \\
\horzbar & \xi_d & \horzbar
\end{bmatrix} 
\begin{bmatrix}
\vertbar & & \vertbar\\
\bm{\psi}\left(\vx(t_1)\right) & \cdots & \bm{\psi}\left(\vx(t_m)\right)\\
\vertbar & & \vertbar
\end{bmatrix}.
\end{align*}
 Or alternatively, viewing the dynamic from an integral perspective, we have that:
 \begin{align*}
\begin{bmatrix}
\vertbar & & \vertbar\\
\vx(t_2) -\vx(t_1)  & \cdots & \vx(t_m) - \vx(t_1)\\
\vertbar & & \vertbar
\end{bmatrix} = 
\begin{bmatrix}
\horzbar & \xi_1 & \horzbar\\
 & \vdots & \\
\horzbar & \xi_d & \horzbar
\end{bmatrix} 
 \begin{bmatrix}
\int_{t_1}^{t_2} \psi_1(\vx(\tau)) d\tau & \cdots & \int_{t_1}^{t_m} \psi_1(\vx(\tau)) \\
\vdots & \ddots & \vdots\\
\int_{t_1}^{t_2} \psi_n(\vx(\tau)) d\tau &\cdots & \int_{t_1}^{t_m} \psi_n(\vx(\tau))
\end{bmatrix}.
\end{align*}
If we collect the data in matrices $\delta X = \left[ \vx(t_2) - \vx(t_1),\cdots, \vx(t_m) - \vx(t_1)\right] \in\rr^{d\times m-1}$, $\dot{X} = \left[ \dot{\vx}(t_1),\cdots, \dot{\vx}(t_m)\right] \in\rr^{d\times m}$, $\Psi\left(X\right) = \left[ \bm{\psi}(\vx(t_1)),\cdots, \bm{\psi}(\vx(t_m))\right]\in\rr^{n\times m}$, and $\Gamma(X)\in \rr^{n \times m -1}$ with $\Gamma(X)_{i,j} = \int_{t_1}^{t_{j+1}} \psi_i(\vx(\tau)) d\tau$, we can view the dynamic from two perspectives:
\begin{mdframed}[linewidth=0.5mm]
\begin{minipage}{0.5\textwidth}
\begin{nscenter}
\textbf{Differential approach}
\end{nscenter}
\begin{align*}
    \dot{X} = \Xi^T \Psi(X)
\end{align*}
\end{minipage}
\vline
\begin{minipage}{0.5\textwidth}
\begin{nscenter}
\textbf{Integral approach}
\end{nscenter}
\begin{align*}
    \delta X = \Xi^T \Gamma(X)
\end{align*}
\end{minipage}
\end{mdframed}
Consequently, when we try to recover the sparsest dynamic that fits this dynamic, we can attempt to solve one of two problems, which we present in tandem:
\begin{mdframed}[linewidth=0.5mm]
\begin{minipage}{0.5\textwidth}
\begin{nscenter}
\textbf{Differential approach}
\end{nscenter}
\begin{gather*}
\argmin\limits_{\substack{  \dot{X} = \Omega^T \Psi(X)  \\ \Omega \in \rr^{n \times d}}} \norm{\Omega}_0. \label{eq:l0_minimization_differential}
\end{gather*}
\end{minipage}
\vline
\begin{minipage}{0.5\textwidth}
\begin{nscenter}
\textbf{Integral approach}
\end{nscenter}
\begin{gather}
 \argmin\limits_{\substack{  \delta X = \Omega^T \Gamma(X)  \\ \Omega \in \rr^{n \times d}}} \norm{\Omega}_0. \label{eq:l0_minimization}
\end{gather}
\end{minipage}
\end{mdframed}
Note that in the previous problem formulation we are implicitly
assuming that we can compute $\Gamma(X)_{i,j}$, which is usually never
the case. In practice we have to resort to approximating the integrals
using quadrature over the given data, that is, for example
$\Gamma(X)_{i,j} \approx \frac{1}{2}\sum_{k=1}^j
(\psi_i(\vx(t_k))+\psi_i(\vx(t_{k+1})))$. If we have access to
$\dot{X}$ and $X$, it will make sense to attack the problem from a
differential perspective, but if we only have access to $X$, and we
have to estimate $\dot{X}$ from data, there are occasions where we can
benefit from the integral approach, as we can potentially estimate
$\Gamma(X)$ more accurately than $\dot{X}$; this can be true in
particular in the presence of noise. Henceforth, we use
$\dot{X}$ and $\Gamma(X)$ to denote the approximate matrices computed
using numerical rules, as opposed to the exact differential and
integral matrices. Unfortunately, the problems shown in
Equations~\eqref{eq:l0_minimization} are notoriously difficult NP-hard
combinatorial problems, due to the presence of the $\ell_0$ norm in
the objective function of the minimization problem of both
optimization problems \cite{juditsky2020statistical}. Moreover, if the
data points are contaminated by noise, leading to noisy matrices
$\dot{Y}$, $\delta Y$, $\Psi(Y)$ and $\Gamma(Y)$, depending on the
expressive power of the basis functions $\psi_i$ for
$i \in \llbracket 1, n\rrbracket$, it may not even be possible (or desirable) to
satisfy $\dot{Y} = \Omega^T \Psi(Y)$ or
$\delta Y = \Omega^T \Gamma(Y)$ for any $\Omega \in \rr^{n\times
  d}$. Thus one can attempt to solve, for a suitably chosen
$\varepsilon > 0$:
\begin{mdframed}[linewidth=0.5mm]
\begin{minipage}{0.5\textwidth}
\begin{nscenter}
\textbf{Differential approach}
\end{nscenter}
\begin{gather*}
\argmin\limits_{\substack{  \norm{\dot{Y} - \Omega^T \Psi(Y)}_F \leq \varepsilon \\ \Omega \in \rr^{n \times d}}} \norm{\Omega}_0.
\end{gather*}
\end{minipage}
\vline
\begin{minipage}{0.5\textwidth}
\begin{nscenter}
\textbf{Integral approach}
\end{nscenter}
\begin{gather}
\argmin\limits_{\substack{  \norm{ \delta Y - \Omega^T \Gamma(Y)}_F \leq \varepsilon \\ \Omega \in \rr^{n \times d}}} \norm{\Omega}_0. \label{eq:l0_minimization_noisy}
\end{gather}
\end{minipage}
\end{mdframed}

The most popular sparse recovery framework, dubbed SINDy
\cite{brunton2016discovering}, solves a component-wise relaxation of a
problem very closely related to the differential problem formulation shown in
Equation~\eqref{eq:l0_minimization_noisy}
\cite{zhang2019convergence}. Each step of the SINDy algorithm consists
of a least-squares step and a thresholding step. The coefficients that
have been thresholded are discarded in future iterations, making the
least-squares problem progressively smaller. More specifically this
process, when applied to one of the components of the problem,
converges to (one of) the local minimizers of:
\begin{gather}
\argmin\limits_{\xi_j \in \rr^{d}} \sum_{i=1}^m \norm{\dot{x}_j(t_i) - \xi_j^T \bm{\psi}(\vx(t_i))}_2^2 + \alpha \norm{\xi_j}_0, \label{eq:l0_minimization_differential_SINDy}
\end{gather}
for a suitably chosen $\alpha \geq 0$ \cite{zhang2019convergence} and for $j \in \llbracket 1,d \rrbracket$. That is, the CINDy framework is the result of applying the sequentially-thresholded
least-squares optimization algorithm to the non-convex problem formulation in Equation~\eqref{eq:l0_minimization_differential_SINDy}. This methodology was later extended to partial differential equations by appropriately modifying the problem formulation and using an optimization algorithm that alternated between ridge-regression steps (as opposed to least-squares steps) and thresholding steps in \cite{rudy2017data}.

In another seminal paper \citet{schaeffer2017sparse} framed the sparse recovery problem from an integral perspective for the first time, using the Douglas-Rachford algorithm \cite{combettes2011proximal} to solve the non-convex integral problem formulation in Equation~\eqref{eq:l0_minimization_noisy}. They showed experimentally that when the data is contaminated with noise and information about the derivatives has to be computed numerically, it can be advantageous to use the integral approach, as opposed to the differential approach, as the numerical integration is more robust to noise than numerical differentiation.

However, both problem formulations in
Equation~\eqref{eq:l0_minimization_noisy} remain non-convex, and so as
is often done in optimization, we can attempt to \emph{convexify} the
problematic term in the problem formulation, namely substituting the
$\ell_0$ norm for the $\ell_1$ norm. Note that the smallest value of $p\geq 0$ that results in the norm $\norm{\cdot}_{p,p}$ being convex is $p = 1$.  This leads us to a problem, known as \emph{basis pursuit denoising} (BPD) \cite{chen1998atomic}, which can be written as:
\begin{mdframed}[linewidth=0.5mm]
\begin{minipage}{0.5\textwidth}
\begin{nscenter}
\textbf{BPD Differential approach}
\end{nscenter}
\begin{gather*}
\argmin\limits_{\substack{  \norm{\dot{Y} - \Omega^T \Psi(Y)}^2_F \leq \epsilon  \\ \Omega \in \rr^{n \times d}}}  \norm{\Omega}_{1,1}
\end{gather*}
\end{minipage}
\vline
\begin{minipage}{0.5\textwidth}
\begin{nscenter}
\textbf{BPD Integral approach}
\end{nscenter}
\begin{gather}
\argmin\limits_{\substack{  \norm{\delta Y - \Omega^T \Gamma(Y)}^2_F \leq \epsilon  \\ \Omega \in \rr^{n \times d}}}  \norm{\Omega}_{1,1} \label{eq:l1_minimization_noisy2}
\end{gather}
\end{minipage}
\end{mdframed}
for appropriately chosen $\epsilon > 0$. The formulation shown in Equation~\eqref{eq:l1_minimization_noisy2} initially developed by the signal processing community, is intimately tied to the \emph{Least Absolute Shrinkage and Selection Operator} (LASSO) regression formulation \cite{tibshirani1996regression}, developed in the statistics community, which takes the form:
\begin{mdframed}[linewidth=0.5mm]
\begin{minipage}{0.5\textwidth}
\begin{nscenter}
\textbf{LASSO Differential approach}
\end{nscenter}
\begin{gather*}
\argmin\limits_{\substack{ \norm{\Omega}_{1,1} \leq \alpha  \\ \Omega \in \rr^{n \times d}}}  \norm{\dot{Y} - \Omega^T \Psi(Y)}^2_F
\end{gather*}
\end{minipage}
\vline
\begin{minipage}{0.5\textwidth}
\begin{nscenter}
\textbf{LASSO Integral approach}
\end{nscenter}
\begin{gather}
\argmin\limits_{\substack{ \norm{\Omega}_{1,1} \leq \alpha  \\ \Omega \in \rr^{n \times d}}}  \norm{\delta Y - \Omega^T \Gamma(Y)}^2_F \label{eq:l1_minimization_noisy}
\end{gather}
\end{minipage}
\end{mdframed}
In fact, the differential approach to the LASSO problem formulation shown in Equation~\eqref{eq:l1_minimization_noisy} was used in \cite{schaeffer2017learning} in conjunction with the Douglas-Rachford algorithm \cite{combettes2011proximal} to solve the sparse recovery problem. A variation of the LASSO problem formulation was also used to recover the governing equations in chemical reaction systems \cite{hoffmann2019reactive} using a sequential quadratic optimization algorithm. The following proposition formalizes the relationship between the BPD and the LASSO problems.
\begin{proposition}\label{prop:relation_BPD_LASSO}~\citet{foucart2017mathematical}[Proposition 3.2]
\begin{enumerate}
    \item If $\Xi$ is the unique minimizer of the BPD problem shown in Equation~\eqref{eq:l1_minimization_noisy2} with $\epsilon > 0$, then there exists an $\alpha \geq 0$ such that $\Xi$ is the unique minimizer of the LASSO problem shown in Equation~\eqref{eq:l1_minimization_noisy}.
    \item If $\Xi$ is a minimizer of the LASSO problem shown in Equation~\eqref{eq:l1_minimization_noisy} with $\alpha > 0$, then there exists an $\epsilon \geq 0$ such that $\Xi$ is a minimizer of the BPD problem shown in Equation~\eqref{eq:l1_minimization_noisy2}. 
\end{enumerate}
\end{proposition}
Both problems shown in Equation~\eqref{eq:l1_minimization_noisy2} and
\eqref{eq:l1_minimization_noisy} have a convex objective function and
a convex feasible region, which allows us to use the powerful tools
and guarantees of convex optimization. These will be the problem formulations on which we will focus to build our framework. Note that these two formulations can also
be recast as an unconstrained optimization problem (via Lagrange dualization) in which the
$\ell_1$ norm has been added to the objective function (see
\citet{foucart2017mathematical} and \citet{borwein2010convex} for more
details). Moreover, there is a significant body of theoretical
literature, both from the statistics and the signal processing
community, on the conditions for which we can successfully recover the
support of $\Xi$ (see e.g., \citet{wainwright2009sharp}), the uniqueness of the
LASSO solutions (see e.g., \citet{tibshirani2013lasso}), or the robust
reconstruction of phenomena from incomplete data
(see e.g., \citet{candes2006robust}), to name but a few results.

\begin{remark}[From learning ODE's to learning PDE's]
    Section~\ref{Section:LearningSparseDynamics} so far has only dealt with the case where the dynamic is expressed as a \emph{ordinary differential equation} (ODE). This framework can also be extended to deal with the case of a dynamic expressed as a \emph{partial differential equation} (PDE), by simply adding the necessary partial derivatives as ansatz functions to the regression problem \cite{schaeffer2017learning, rudy2017data}.
\end{remark}

\subsection{Incorporating
  structure} \label{section:IncorporatingStructure} Conservation laws
are a fundamental pillar of our understanding of physical
systems. These conservation laws stem from differentiable symmetries
that are present in nature \cite{noether1918invariante}. Imposing
these symmetry constraints in our sparse regression problem can
potentially lead to better generalization performance under noise,
reduced sample complexity, and to learned dynamics that are consistent
with the symmetries present in the real world. Our approach allows for
arbitrary polyhedral constraints to be added, i.e., linear inequality
and equality constraints; boundedness will be ensured automatically
due to the \(\ell_1\) norm constraint. In particular, there are two
large classes of structural constraints that can be easily encoded
into our learning problem.

\subsubsection{Conservation properties} \label{Section:Convervation}

From a differential perspective, we often observe in dynamical systems that certain relations hold between the elements of $\dot{\vx}(t)$. Such is the case in chemical reaction dynamics, where if we denote the rate of change of the $i$-th species by $\dot{x}_i(t)$, we might observe relations of the form $a_j\dot{x}_j(t) + a_k\dot{x}_k(t) = 0$ due to mass conservation, which relate the $j$-th and $k$-th species being studied.
%We can also observe these phenomena in chemical reaction systems, where the rate of change of the different compounds involved can often be related. 
In the case where these relations are linear, we know that for some $\mathcal{J}\subseteq \llbracket 1,d \rrbracket$ and all $t\geq 0$ we can write:
\begin{align}
\sum\limits_{j\in \mathcal{J}} a_j \dot{x}_j(t) = c. \label{eq:conservation_derivatives}
\end{align}
We can encode Equation~\eqref{eq:conservation_derivatives} into our learning problem by using the fact that $\dot{x}_j(t) = \xi_j^T \bm{\psi}(\vx(t))$ and imposing that for all data points $i \in \llbracket 1,m \rrbracket$
\begin{align*}
\sum\limits_{j\in \mathcal{J}} a_j \xi_j^T \bm{\psi}(\vx(t_i)) = c,
\end{align*}
which can be expressed more succinctly as
\begin{align*}
\sum\limits_{j\in \mathcal{J}} a_j \xi_j^T \Psi(X) = c\mathbf{1}_m.
\end{align*}
This involves the addition of $m$ linear constraints into our learning
problem, which in the absence of noise does not pose any
problems. However, when the data
$\left\{ \dot{\vx}(t_i), \vx(t_i) \right\}_{i=1}^m$ is contaminated by
noise, and we only have access to
$\left\{ \dot{\vy}(t_i), \vy(t_i) \right\}_{i=1}^m$, it is futile to
assume that $\sum\limits_{j\in \mathcal{J}} a_j \dot{y}_j(t_i) = c$
for all $i\in \llbracket 1, m \rrbracket$ or that
$\sum\limits_{j\in \mathcal{J}} a_j \xi_j^T \Psi(Y) =
c\mathbf{1}_m$. In this case, it is more reasonable to assume that the
derivatives are approximately preserved, and instead impose for some
$\varepsilon > 0$ and all $i \in \llbracket 1,m \rrbracket$ that:
\begin{align*}
\left\lvert \sum\limits_{j\in \mathcal{J}} a_j \xi_j^T \bm{\psi}(\vy(t_i)) - c \right\rvert \leq \varepsilon.
\end{align*}
The addition of this constraint to the problem in
Equation~\eqref{eq:l1_minimization_noisy} preserves the convexity of
the original problem. Moreover, the feasible region of the
optimization problem remains polyhedral.

\subsubsection{Symmetry between variables} \label{section:symmetry}
\begin{figure}[th!]
\centering
  \includegraphics[width=0.5\linewidth]{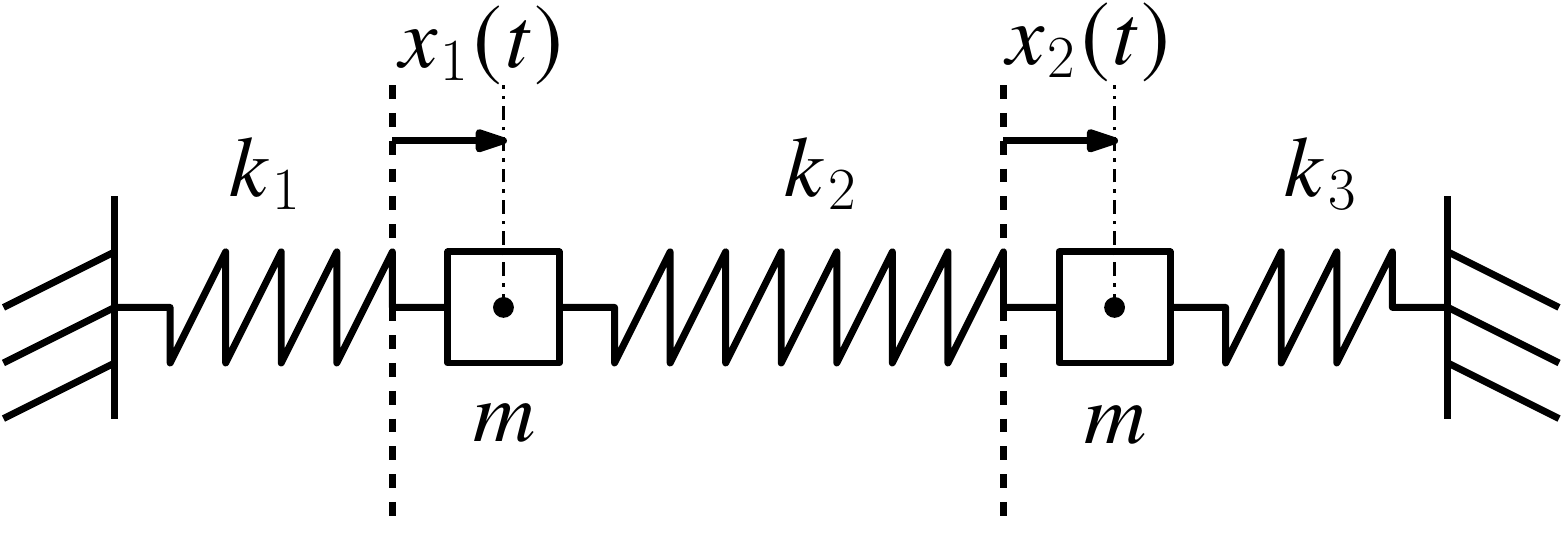}
 \caption{Two spring-coupled masses.}
  \label{fig:two_spring_mass}
\end{figure}

One of the key assumptions used in many-particle quantum systems is the fact the particles being studied are indistinguishable. And so it makes sense to assume that the effect that the $i$-th particle exerts on the $j$-th particle is the same as the effect that the $j$-th particle exerts on the $i$-th particle. The same can be said in classical mechanics for a collection of identical masses, where each mass is connected to all the other masses through identical springs. As an example, consider the system formed by two spring-coupled masses depicted in Figure~\ref{fig:two_spring_mass}.

Here we denote the displacement of the center of mass of the $i$-th body from its equilibrium position at rest by $x_i(t)$. This allows us to express the dynamical evolution of the system by $m \ddot{x}_1 (t) = -k_1 x_1(t) + k_2(x_2(t) - x_1(t))$ and $m \ddot{x}_2 (t) = -k_2(x_2(t) - x_1(t)) - k_3 x_2(t)$. Suppose we  are given access to a series of noisy data points $\left\{ \ddot{\vy}(t_i), \vy(t_i) \right\}_{i=1}^m$ and we want to learn the dynamic $\ddot{\vy}(t) = \left[ \ddot{y}_1(t), \ddot{y}_2(t)\right]$ with a dictionary $\mathcal{D} = \{ \psi_1 (\vy) = 1, \psi_2 (\vy) = y_1, \psi_3 (\vy) = y_2 \}$ of basis functions of polynomials of degree up to one. The problem is analogous to that of learning $\dot{\vx}(t)$ and can be framed similarly to that of Equation~\eqref{eq:l1_minimization_noisy}, substituting $\dot{Y}$ for $\ddot{Y}$. If we use $\xi_j\left( \psi(\vx)\right)$ to refer to the coefficient in $\xi_j$, where $\xi_j$ is the $j$-th column of $\Omega$, associated with the basis function $\psi(\vx)$, we can write:

\begin{align}
\hat{\Xi} = \argmin\limits_{\substack{ \norm{\Omega}_{1,1} \leq \tau  \\ \Omega \in \rr^{n \times d}}}  \norm{\begin{bmatrix}
\ddot{y}_1(t_1) &\cdots & \ddot{y}_1(t_m)\\
\ddot{y}_2(t_1) & \cdots & \ddot{y}_2(t_m)
\end{bmatrix} - 
\begin{bmatrix}
\xi_1\left( 1 \right) & \xi_1\left( y_1 \right) & \xi_1\left( y_2 \right)  \\
\xi_2\left( 1 \right) & \xi_2\left( y_1 \right) & \xi_2\left( y_2 \right) 
\end{bmatrix}
 \begin{bmatrix}
1 & \cdots & 1 \\
y_1 (t_1)& \cdots & y_1(t_m)\\
y_2 (t_1) &\cdots & y_2(t_m)
\end{bmatrix}}^2_F \label{eq:l1_minimization_noisy_structure}.
\end{align}
Where we have that
$\xi_1 = \left[ \xi_1\left( 1 \right),\xi_1\left( y_1 \right),
  \xi_1\left( y_2 \right) \right]$ and
$\xi_2 = \left[ \xi_2\left( 1 \right), \xi_2\left( y_1 \right),
  \xi_2\left( y_2 \right) \right]$. In light of the structure of the
system and its symmetry, it makes sense to add to the learning problem
the constraint $\xi_1\left( y_2 \right) = \xi_2\left( y_1 \right)$,
that is, the effect of $y_1(t)$ on $\ddot{y}_2(t)$ is the same as the
effect of $y_2(t)$ on $\ddot{y}_1(t)$. These constraints can also be
readily applied in the integral formulation of the LASSO recovery problem.

\subsubsection{Existing structured sparse recovery frameworks} \label{Section:Existing algorithms}

In the context of the sparse recovery of dynamics, several optimization algorithms have been proposed to enforce \emph{linear equality constraints} in the problem formulation, as opposed to the more general linear inequality constraints. In the \emph{Constrained Sparse Garlerkin Regression} framework \citet{loiseau2018constrained} propose an optimization algorithm that alternates between solving a quadratic problem subject to linear equality constraints using an \emph{Interior-Point Method} (IPM) \cite{nesterov1994interior}, and thresholding the coefficients in a similar way as is done in the optimization of the SINDy framework. The satisfaction of the linear equality constraints is imposed by the IPM, while the sparsity of the dynamic is enforced by the thresholding step. Much like the SINDy framework, it does not produce sparse solutions in the presence of mild noise, although it successfully incorporates constraints. We remark that one could substitute the use of the IPM by a step that explicitly solves the Karush-Kuhn-Tucker (KKT) conditions of the quadratic problem subject to linear equality constraints.

Another popular approach is to use the \emph{Sparse
  Relaxed Regularized Regression} (SR3) framework \cite{zheng2018unified, champion2020unified}, which formulates a relaxation of the regularized problem. For example, focusing on the differential formulation subject to linear equality constraints, the SR3 framework would result in a problem formulation of the form:
  \begin{align}
\argmin\limits_{\substack{ \Omega \in \mathcal{P}  \\ \Omega, W \in \rr^{n \times d}}}  \norm{\dot{Y} - \Omega^T \Psi(Y)}^2_F  + \alpha \norm{W}_0 +  \frac{1}{ \nu}\norm{\Omega - W}^2_F \label{eq:SR3_problem},
\end{align}
  where we use $\mathcal{P}$ to refer to a polytope that enforces the appropriate linear constraints on $\Omega$, $\alpha > 0$ controls the regularization of the relaxed variable $W$ and $\nu >0$ controls the penalty between $\Omega$ and $W$. Note that in the limit of values of $\alpha$ approaching zero in Equation~\eqref{eq:SR3_problem} we would recover the non-convex problem formulation solved by SINDy, shown in Equation~\eqref{eq:l0_minimization_differential_SINDy}. In order to tackle the SR3 problem formulation \citet{champion2020unified} propose the use of a proximal gradient descent optimization algorithm, which consists of alternatively minimizing $\Omega$ given a fixed $W$, which can be done through solving the KKT conditions, and then updating $W$ by applying the proximal operator associated with the $\ell_0$ norm to $W$. This is repeated until some measure of convergence is reached. The drawback in this problem formulation is that it requires careful control of the trade-off between the conditioning and the fidelity to the original problem and the sparsity imposed on the problem. Moreover, it is designed to impose linear equality constraints, as opposed to more general inequality constraints.

\section{Conditional Gradient algorithms}
\label{sec:frank-wolfe-algor}
Now that we have selected the LASSO problem formulation shown in Equation~\eqref{eq:l1_minimization_noisy} for our framework, we focus on the optimization algorithm used to find an approximate solution to the problem formulation. For simplicity, let us assume that we are dealing with the differential formulation of the problem in Equation~\eqref{eq:l1_minimization_noisy}, thus we would like to solve:
\begin{align}
 \argmin\limits_{\substack{ \norm{\Omega}_{1,1} \leq \alpha  \\ \Omega \in \rr^{n \times d}}} f(\Omega) \label{eq:l1_minimization_noisy_2},
\end{align}
 where $f(\Omega) = \norm{\dot{Y} - \Omega^T \Psi(Y)}^2_F$. This can be done using first-order projection-based algorithms such as \emph{gradient descent} or \emph{accelerated gradient descent}. Using the former, the iterate at iteration $k+1$ can be expressed, for a suitably chosen step size $\gamma_k>0$ as:
\begin{align}
  \Omega_{k+1} &= \argmin\limits_{\substack{ \norm{\Omega}_{1,1} \leq \tau  \\ \Omega \in \rr^{n \times d}}}  \norm{\Omega - \left(\Omega_{k} - \gamma_k \nabla f\left(\Omega_k \right) \right)}^2_F \\
  & = \argmin\limits_{\substack{ \norm{\Omega}_{1,1} \leq \tau  \\ \Omega \in \rr^{n \times d}}}  \norm{\Omega -\Omega_{k} - 2 \gamma_k\Psi(Y) \left(\dot{Y} - \Omega_k^T\Psi(Y) \right)^T }^2_F. \label{eq:projected_gradient_descent}
\end{align}

Fortunately, the quadratic problem shown in Equation~\eqref{eq:projected_gradient_descent} can be solved exactly with complexity $\mathcal{O}(nd)$ \cite{condat2016fast} (as this is equivalent to projecting a flattened version of the matrix onto the $\ell_1$ polytope of dimension $nd$). If we were to add $L$ additional linear constraints to the problem in Equation~\eqref{eq:l1_minimization_noisy} to reflect the underlying structure of the dynamical system through symmetry and conservation, we would arrive at a polytope $\mathcal{P}$ of the form
\begin{align*}
\mathcal{P} = \left\{ \Omega \in \rr^{n \times d} \mid \norm{\Omega}_{1,1} \leq \alpha, \trace( A_l^T \Xi ) \leq b_l, l\in \llbracket 1,L\rrbracket \right\},
\end{align*}
with $A_l \in \rr^{n\times d}$ and $b_l \in \rr$ for all $l \in \llbracket 1,L\rrbracket$. 
So in this case, with additional structural constraints, the problem would transform into:
\begin{align}
\argmin\limits_{\substack{ \Omega \in \mathcal{P}  \\ \Omega \in \rr^{n \times d}}}  \norm{\dot{Y} - \Omega^T \Psi(Y)}^2_F \label{eq:l1_minimization_noisy_additional_constraints}.
\end{align}

Unfortunately, in general there is no closed-form solution to the projection operator onto $\mathcal{P}$, and so in order to use projection-based algorithms to solve the optimization problem, one has to compute these projections approximately. Note that computing a projection onto $\mathcal{P}$ is equivalent to solving a quadratic problem over $\mathcal{P}$, which can be as expensive as solving the original quadratic problem shown in Equation~\eqref{eq:l1_minimization_noisy_additional_constraints}. In light of this difficulty, one can opt to solve the optimization problem using projection-free algorithms like the \emph{Conditional Gradients} (CG) algorithm \cite{polyak66cg} (also known as the \emph{Frank-Wolfe} (FW) algorithm \cite{fw56}, shown in Algorithm~\ref{algo:CG} with exact line search).

%\begin{algorithm}[th!]
%\SetKwInOut{Input}{Input}\SetKwInOut{Output}{Output}
%\Input{Initial point $\Omega_1 \in \mathcal{P}$.}
%\Output{Point $\Omega_{K+1} \in \mathcal{P}$.}
%\hrulealg
%\For{$k = 1$ to $K$}{
%$V_k \leftarrow \argmin\limits_{\Omega \in  \mathcal{P}} \trace\left(\Omega^T\Psi(Y) \left(\dot{Y} - \Omega_k^T\Psi(Y) \right)^T\right)$ \label{alg:linear_problem}\;
%$D_k \leftarrow V_k - \Omega_k$\;
%$\gamma_k \leftarrow -\frac{1}{2}\trace \left((\dot{Y} - \Omega_k\Psi(Y))\Psi(Y)^T D_k \right)/ \norm{D_k^T \Psi(Y)}_F^2 $\;
%$\Omega_{k+1} \leftarrow \Omega_k + \gamma_k D_k$\;}
%\caption{Frank-Wolfe (FW) algorithm applied to Problem~\eqref{eq:l1_minimization_noisy_additional_constraints}} \label{algo:CG}
%\end{algorithm}

\begin{minipage}{0.50\textwidth}
\begin{algorithm}[H]
\SetKwInOut{Input}{Input}\SetKwInOut{Output}{Output}
\Input{Initial point $\Omega_1 \in \mathcal{P}$.}
\Output{Point $\Omega_{K+1} \in \mathcal{P}$.}
\hrulealg
\For{$k = 1$ to $K$}{
$\nabla f \left( \Omega_k \right) \leftarrow 2\Psi(Y) \left(\dot{Y} - \Omega_k^T\Psi(Y) \right)^T$ \;
$V_k \leftarrow \argmin\limits_{\Omega \in  \mathcal{P}} \trace\left(\Omega^T\nabla f \left( \Omega_k \right) \right)$ \label{alg:linear_problem}\;
$D_k \leftarrow V_k - \Omega_k$\;
$\gamma_k \leftarrow \min\left\{-\frac{1}{2} \frac{\trace \left( D_k^T \nabla f \left( \Omega_k \right) \right)}{\norm{D_k^T \Psi(Y)}_F^2},1\right\} \label{alg:line_search} $\;
$\Omega_{k+1} \leftarrow \Omega_k + \gamma_k D_k$\;}
\caption{CG algorithm applied to \eqref{eq:l1_minimization_noisy_additional_constraints}} \label{algo:CG}
\end{algorithm}
\end{minipage}
% \vline 
\begin{minipage}{0.50\textwidth}
\begin{figure}[H]
\centering
  \includegraphics[width=0.70\linewidth]{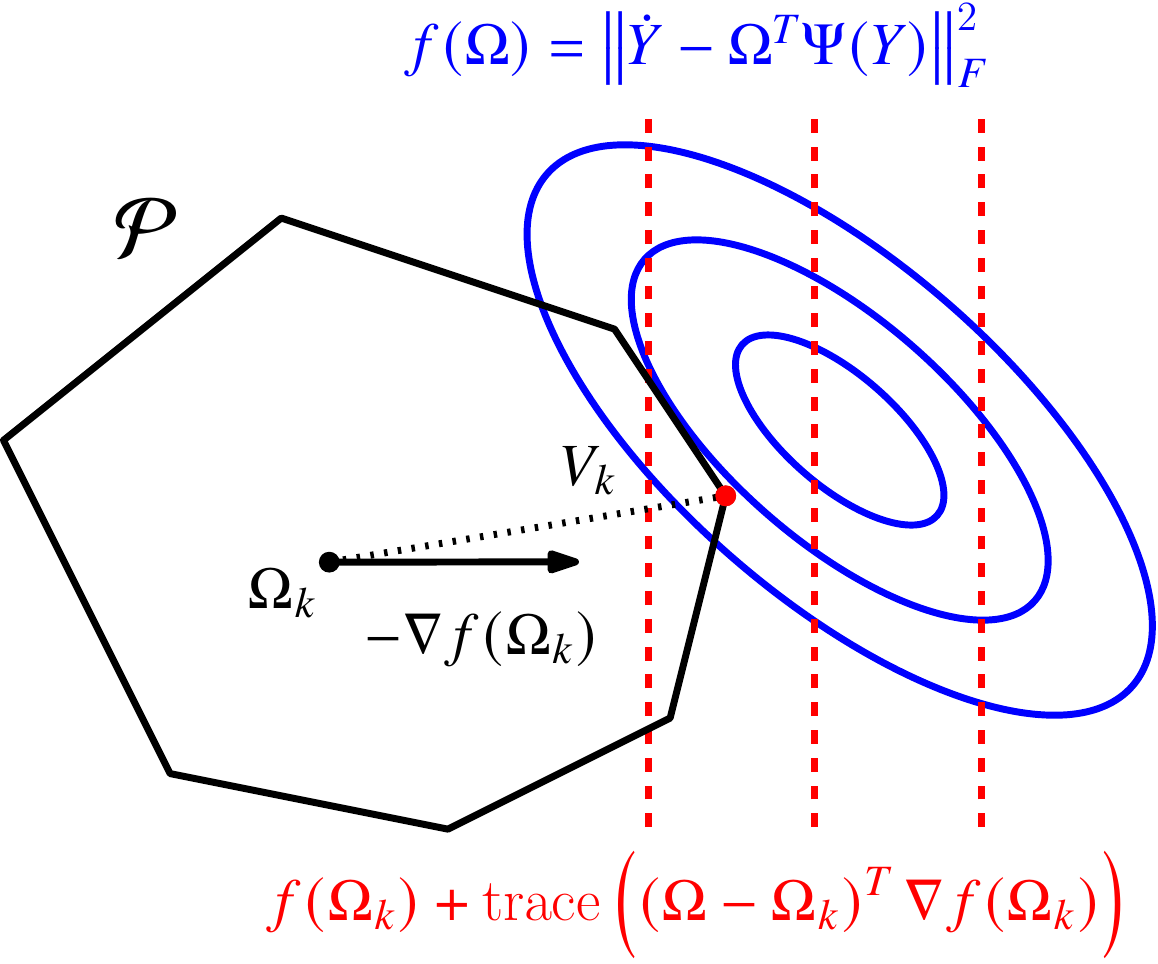}
 \caption{CG algorithm schematic.}
  \label{fig:FWschematic}
\end{figure}
\end{minipage} \medskip

This family of algorithms (including Algorithm~\ref{algo:CG}) requires
solving a linear optimization problem over a polytope
(Line~\ref{alg:linear_problem} of Algorithm~\ref{algo:CG}) at each
iteration, instead of a quadratic problem. As the iterates are
obtained as a convex combination of the current iterate $\Omega_k$ and
the solution of the linear optimization problem over $\mathcal{P}$,
denoted by $V_k$ (Line~\ref{alg:linear_problem} of
Algorithm~\ref{algo:CG}), thus always ensuring feasibility, these
methods are \emph{projection-free}. The direction $V_k - \Omega_k$ is
the direction that best approximates (in the inner product sense) the negative of the gradient of
the objective function at the current iterate $\Omega_k$, in this case
$ -\nabla f(\Omega_k) = -2\Psi(Y) \left(\dot{Y} - \Omega_k^T\Psi(Y)
\right)^T$, if we restrict ourselves to moving towards vertices of the
polytope $\mathcal{P}$. To be more precise:
\begin{align*}
   \trace\left( -\left( V_k - \Omega_k\right)^T \nabla f(\Omega_k) \right)&= \max\limits_{ \Omega \in \mathcal{P}} \trace\left(-\left(\Omega - \Omega_k\right)^T \nabla f(\Omega_k)\right) \\
    & = \max\limits_{ \Omega \in \mathcal{P}} \trace\left(-2\left(\Omega - \Omega_k\right)^T\Psi(Y) \left(\dot{Y} - \Omega_k^T\Psi(Y) \right)^T\right).
\end{align*}
This can be seen as equivalent to moving along the direction given by the vertex which minimizes a linear approximation of the objective function at the current iterate $\Omega_k$ over the polytope $\mathcal{P}$ (see Figure~\ref{fig:FWschematic}, where we have denoted the objective function as $f(\Omega) = \norm{\dot{Y} - \Omega^T \Psi(Y)}^2_F$), that is:
\begin{align*}
    V_k  & = \argmin\limits_{ \Omega \in \mathcal{P}} f(\Omega_k) + \trace\left(\left(\Omega - \Omega_k\right)^T \nabla f(\Omega_k)\right) \\
    & = \argmin\limits_{ \Omega \in \mathcal{P}} \norm{\dot{Y} - \Omega_k^T \Psi(Y)}^2_F + 2\trace\left(\left(\Omega - \Omega_k\right)^T\Psi(Y) \left(\dot{Y} - \Omega_k^T\Psi(Y) \right)^T\right).
\end{align*}

%\begin{figure}[th!]
%\centering
%  \includegraphics[width=0.4\linewidth]{Images/ApproximationSchematic_CG.pdf}
% \caption{Schematic of step taken by the Frank-Wolge algorithm.}
%  \label{fig:FWschematic}
%\end{figure}S

Once the vertex $V_k$ has been found, the exact line search solution
is computed in Line~\ref{alg:line_search} to find the step size
$\gamma_k$ that results in the greatest decrease in primal gap, that is,
$\gamma_k = \argmin_{\gamma \in [0,1]} f(\Omega_k + \gamma (V_k -
\Omega_k))$. Fortunately, as the function being minimized (see
Equation~\eqref{eq:l1_minimization_noisy_additional_constraints}) is a
quadratic, there is a closed form expression for the optimal step
size. Note that the step size in Line~\ref{alg:line_search} is always
non-negative and the clipping ensures that we build convex
combinations; the clipping is active (if at all) only in the very
first iteration by standard arguments (see, e.g., \cite{braun2021conditional}).

\begin{remark}
  If we assume that the starting point $\Omega_1$ is a vertex of the
  polytope, then we know that the iterate $\Omega_k$ can be expressed
  as a convex combination of at most $k$ vertices of
  $\mathcal{P}$. This is due to the fact that the algorithm can pick
  up no more than one vertex per iteration. Note that the CG algorithm
  applied to the problem shown in
  Equation~\eqref{eq:l1_minimization_noisy_2}, where the feasible
  region is the $\ell_1$ ball without any additional constraints,
  picks up at most one basis function in the $k$-th iteration, as
  $V_k^T \bm{\psi}(\vx(t)) = \pm \tau \psi_i(\vx(t)$ for some
  $i \in \llbracket 1,n \rrbracket$. This means that if we use the
  CG algorithm to solve a problem over the $\ell_1$ ball, we
  encourage sparsity not only through the regularization provided by
  the $\ell_1$ ball constraint in the problem formulation, but also \emph{through the specific nature of the
    CG algorithm} independently of the size of the feasible
  region. In practice, when using, e.g., early termination due to some
  stopping criterion, this results in the CG algorithm
  producing sparser solutions than projection-based algorithms (such
  as projected gradient descent, which typically use dense updates) when applied to
  Problem~\ref{eq:l1_minimization_noisy_additional_constraints},
  despite the fact that both algorithms converge to the same solution
  if the problem is strictly convex.

  Thus, in addition to the trade-off between reconstruction accuracy
  and sparsity offered by LASSO problem formulations parametrized by the size
  $\alpha$ of the \(\ell_1\) ball, we have the same trade-off in terms
  of the iteration count. Similar to iterative or semi-iterative
  regularization methods such as Landweber's method for ill-posed
  linear problems in Hilbert spaces~\cite{Hanke1991}, reconstruction
  accuracy improves in iterations, while the norm of the
  solution, in our case the $l_0$ norm, tends to grow. The trade-off
  can be decided by a termination criterion, e.g., a sufficiently small
  residual in Morozow's discrepancy principle~\cite{Morozov1966} or,
  in our case, a sufficiently small primal gap.
\end{remark}

One of the interesting properties of CG algorithms is the fact that, since $f$ is convex, at each iteration we can compute the \emph{Frank-Wolfe} gap, an upper bound on the primal gap, at no extra cost.

\begin{definition}[Frank-Wolfe gap]
The Frank-Wolfe gap of the function $f$ over the feasible region $\mathcal{P}$ evaluated at $\Omega_k$, denoted by $g_{\mathcal{P}}\left( \Omega\right)$, is given by:
\begin{align*}
    g_{\mathcal{P}}\left( \Omega_k\right) = \max_{\Omega \in \mathcal{P}} \trace \left( \left( \Omega_k - \Omega\right)^T \nabla f(\Omega_k)\right). %\label{eq:FW_gap}
\end{align*}
To see why this quantity provides an upper bound on the primal gap $f\left(\Omega_k\right) - \min_{\Omega \in \mathcal{P}} f\left(\Omega\right)$, when $f\left( \Omega\right)$ is convex, note that if we denote $\Omega^* = \argmin_{\Omega \in \mathcal{P}} f\left(\Omega\right)$, then
\begin{align}
  f\left(\Omega_k\right) -  f\left(\Omega^*\right)  &\leq \trace \left( \left( \Omega_k - \Omega^*\right)^T \nabla f(\Omega_k)\right) \label{eq:convexity_CG} \\
  &\leq \max\limits_{\Omega \in \mathcal{P}} \trace \left( \left( \Omega_k - \Omega\right)^T \nabla f(\Omega_k)\right) \\
  & = \trace \left( \left( \Omega_k - V_k\right)^T \nabla f(\Omega_k)\right),
\end{align}
holds with Equation~\eqref{eq:convexity_CG} following from convexity of $f$.
\end{definition}
To sum up the advantages of CG algorithms when applied to structured sparse LASSO recovery problem formulations:
\begin{enumerate}
    \item Sparsity is encouraged through a two-fold approach: through the $\ell_1$ regularization used in the problem formulation and through the use of the Conditional Gradient algorithms, which are sparse in nature.
    \item Linear equality and inequality constraints can be added
      easily and naturally to the constraint set of the problem to
      reflect symmetry or conservation assumptions. These additional
      constraints can be efficiently managed due to the fact that
      there are extremely efficient algorithms to solve linear
      programs over polytopes.
\end{enumerate}
These characteristics make the class of Conditional Gradient methods
extremely attractive versus projection-based algorithms to solve
LASSO recovery problem formulations.

\begin{remark}
  We remark that other optimization algorithms can be used to solve the constrained LASSO problem formulation, such as the \emph{Alternating Direction Method of Multipliers} (see \citet{james2012constrained, gaines2018algorithms} for an overview), however, they do not have an algorithmic bias towards sparse solution, as opposed to CG algorithms.
\end{remark}

\subsection{Fully-Corrective Conditional Gradients} \label{section:FCCG}

For the recovery of sparse dynamics from data, one of the most
interesting algorithms in terms of sparsity is the
\emph{Fully-Corrective Conditional Gradient} (FCCG) algorithm
(Algorithm~\ref{algo:FCFW}). This algorithm picks up a vertex $V_k$
from the polytope $\mathcal{P}$ at each iteration
(Line~\ref{alg:linear_problemFCFW} of Algorithm~\ref{algo:FCFW}) and
reoptimizes over the convex hull of $\mathcal{S}_{k} \bigcup V_k$
(Line~\ref{alg:quadratic_problemFCFW} of Algorithm~\ref{algo:FCFW}),
which is the union of the vertices picked up in previous iterations,
and the new vertex $V_k$. The reoptimization step can potentially
remove a large number of unnecessary vertices picked up in earlier
iterations.

\begin{algorithm}[th!]
\SetKwInOut{Input}{Input}\SetKwInOut{Output}{Output}
\Input{Initial point $\Omega_1 \in \mathcal{P}$.}
\Output{Point $\Omega_{K+1} \in \mathcal{P}$.}
\hrulealg
$\mathcal{S}_{1}  \leftarrow \emptyset$  \;
\For{$k = 1$ to $K$}{
$\nabla f \left( \Omega_k \right) \leftarrow 2 \Psi(Y) \left(\dot{Y} - \Omega_k^T\Psi(Y) \right)^T$ \;
$V_k \leftarrow \argmin\limits_{\Omega \in  \mathcal{P}} \trace\left(\Omega^T\nabla f \left( \Omega_k \right) \right)$ \label{alg:linear_problemFCFW}\;
$\mathcal{S}_{k+1}\leftarrow \mathcal{S}_{k} \bigcup V_k$\;
$\Omega_{k+1} \leftarrow \argmin\limits_{\Omega \in \co\left( \mathcal{S}_{k+1} \right) } \norm{\dot{Y} - \Omega^T \Psi(Y)}^2_F$\label{alg:quadratic_problemFCFW}\;}
\caption{Fully-Corrective Conditional Gradient (CG) algorithm applied to Problem~\eqref{eq:l1_minimization_noisy_additional_constraints}} \label{algo:FCFW}
\end{algorithm}

The reoptimization subproblem shown in
Line~\ref{alg:quadratic_problemFCFW} of Algorithm~\ref{algo:FCFW} is a
quadratic problem over a polytope $\mathcal{P}$, like the original
problem in
Equation~\eqref{eq:l1_minimization_noisy_additional_constraints}. However,
it can be rewritten as an optimization problem over the unit
probability simplex of dimension $k $, as the cardinality of the set
$\mathcal{S}_{k+1}$ satisfies $\abs{\mathcal{S}_{k+1}} = k$. To see
this, note that given a set 
$\mathcal{S}_{k+1} \subseteq \vertex \left( \mathcal{P} \right)$ we
can express any $\Omega \in \co\left( \mathcal{S}_{k+1} \right)$ as
$\Omega = \sum_{i = 1}^{k} \lambda_{i} V_i$ for some
$\bm{\lambda}= \left[\lambda_1, \cdots, \lambda_{k}\right]\in
\Delta_{k}$ and $V_i \in \mathcal{S}_{k+1}$ for all
$i\in \llbracket 1,k \rrbracket$. This leads to:
\begin{align*}
\min_{\Omega \in \co\left( \mathcal{S}_{k+1} \right) } \norm{\dot{Y} - \Omega^T \Psi(Y)}_F^2 & = \min_{\Omega \in \co\left( \mathcal{S}_{k+1} \right) } \left( \trace \left( \dot{Y}^T \dot{Y}\right) - 2 \trace \left( \Omega^T\Psi(Y) \dot{Y}^T \right) +  \trace \left( \Omega^T\Psi(Y) \Psi(Y)^T \Omega \right) \right) \\
& =  \min_{\bm{\lambda} \in \Delta^{k} } \left( \norm{\dot{Y}}_F^2 - 2\sum_{i=1}^{k} \lambda_i \trace \left( V_i^T\Psi(Y) \dot{Y}^T \right) + \sum_{i=1}^{k} \sum_{j=1}^{k} \lambda_i \lambda_j \trace \left( V_i^T\Psi(Y) \Psi(Y)^T V_j \right)\right).
\end{align*}
Which can be expressed more succinctly if we denote $\Lambda_{k+1} = \left[ \vect \left( V_1^T \Psi(Y) \right), \cdots,  \vect \left( V_{k}^T \Psi(Y) \right) \right] \in \rr^{dm\times k}$ and we write:
\begin{align}
\min_{\Omega \in \co\left( \mathcal{S}_{k+1} \right) } \norm{\dot{Y} - \Omega^T \Psi(Y)}_F^2 & = \min_{\bm{\lambda} \in \Delta^{k} } \left( \norm{\dot{Y}}_F^2 - 2\vect\left( \dot{Y}\right)^T \Lambda_{k+1} \bm{\lambda} + \bm{\lambda}^T \Lambda_{k+1}^T \Lambda_{k+1} \bm{\lambda} \right) \\
& = \min_{\bm{\lambda} \in \Delta^{k} } \norm{\Lambda_{k+1} \bm{\lambda} - \vect\left( \dot{Y}\right)}^2. \label{eq:unit_probability_simplex} 
\end{align}
So in order to solve the optimization problem in
Line~\ref{alg:quadratic_problemFCFW} of Algorithm~\ref{algo:FCFW} we
would need to solve the optimization problem shown in
Equation~\eqref{eq:unit_probability_simplex} (which is also convex, as
convexity is invariant under affine maps) and take $\Omega = \sum_{i =
  1}^{k} \lambda_{i} V_i$. While the original quadratic problem over
$\mathcal{P}$, shown in
Equation~\eqref{eq:l1_minimization_noisy_additional_constraints}, has
dimensionality $n\times d$, the quadratic problem shown in
Line~\ref{alg:quadratic_problemFCFW} of Algorithm~\ref{algo:FCFW} has
dimensionality $k$ when it is solved in $\bm{\lambda}$-space, which
leads to improved convergence due to reduced problem dimensionality. %In practice the FCFW algorithm usually converges in a number of iterations such that $k << nd$, which means that the optimization problem shown in Line~\ref{alg:quadratic_problemFCFW} of Algorithm~\ref{algo:FCFW} is much cheaper to solve than the original problem, or the projection problem shown in Equation~\eqref{eq:projected_gradient_descent}.

\begin{remark}
  If the polytope \(\mathcal{P}\) being considered is simply the
  $\ell_1$ ball, then computing $\Lambda_{k+1}$ requires at most $mk$
  multiplications, since in this case $\norm{V_i}_0=1$ for all
  $i\in \llbracket 1, k \rrbracket$, as $V_i$ is simply one of the
  vertices of the $\ell_1$ ball (which is a polytope), and so computing
  $\vect \left( V_i^T \Psi(Y) \right)$ requires at most $m$
  multiplications for each $i$. This means that $\Lambda_{k+1}$ is
  sparse, as $\norm{\Lambda_{k+1}}_0 \leq mk$, which allows us to
  efficiently compute
  $\Lambda_{k+1}^T \Lambda_{k+1} \in \rr^{k \times k}$ and
  $\vect\left( \dot{Y}\right)^T \Lambda_{k+1} \in \rr^{k}$.
\end{remark}

\begin{remark}
 If additional constraints are added to the the $\ell_1$ ball, in general, we cannot make any statements about the sparsity of $\Lambda_{k+1}$, other that in numerical experiments we observe that $\norm{\Lambda}_0 \ll dmk$.
\end{remark}

 Due to the fact that there are efficient algorithms to compute projections onto the probability simplex of dimension $k$ with complexity $\mathcal{O}\left(k\right)$ \cite{condat2016fast} we can use accelerated projected gradient descent to solve the subproblems in Line~\ref{alg:quadratic_problemFCFW} of Algorithm~\ref{algo:FCFW} \cite{nesterov1983method, nesterov2018lectures} (shown in Algorithm~\ref{algo:NAGD} and Algorithm~\ref{algo:NAGD_strcvx} in Appendix~\ref{NAGD}). Solving the problem shown in Line~\ref{alg:quadratic_problemFCFW} of Algorithm~\ref{algo:FCFW} to optimality at each iteration is computationally prohibitive, and so ideally we would like to solve the problem in Line~\ref{alg:quadratic_problemFCFW} to $\varepsilon_k$-optimality.

\subsection{Blended Conditional Gradients} \label{section:BCG} 

This leads to the question: How should we choose $\varepsilon_k$ at
each iteration $k$, if we want to find an $\varepsilon$-optimal
solution to the problem shown in
Equation~\eqref{eq:l1_minimization_noisy_additional_constraints}?
Computing a solution to the problem shown in
Line~\ref{alg:quadratic_problemFCFW} to accuracy
$\varepsilon_k = \varepsilon$ at each iteration might be way too
computationally expensive. Conceptually, we need relatively inaccurate
solutions for early iterations where
$\Omega^* \notin \co \left(\mathcal{S}_{k+1}\right)$, requiring only
accurate solutions when
$\Omega^* \in \co \left(\mathcal{S}_{k+1}\right)$. At the same time we
do not know whether we have found \(\mathcal{S}_{k+1}\) so that
$\Omega^* \in \co \left(\mathcal{S}_{k+1}\right)$.

The rationale behind the \emph{Blended Conditional Gradient} (BCG)
algorithm \cite{braun2019blended} (the variant used for our specific
problem is shown in Algorithm~\ref{algo:BCG}) is to provide an
explicit value of the accuracy $\varepsilon_k$ needed at each
iteration starting with rather large \(\varepsilon_k\) in early iterations
and progressively getting more accurate when approaching the optimal
solution; the process is controlled by an optimality gap
measure. In some sense one might think of BCG as a practical version
of FCCG with stronger convergence guarantees and much faster
real-world performance.

\begin{algorithm}[th!]
\SetKwInOut{Input}{Input}\SetKwInOut{Output}{Output}
\SetKwComment{Comment}{$\triangleright$\ }{}
%\SetKwComment{Comment}{//}{}
\Input{Initial point $\Omega_0 \in \mathcal{P}$.}
\Output{Point $\Omega_{K+1} \in \mathcal{P}$.}
\hrulealg
$\Omega_1 \leftarrow \argmin\limits_{\Omega \in  \mathcal{P}} \trace\left(\Omega^T\nabla f \left( \Omega_0 \right) \right)$\;
$\Phi  \leftarrow \trace \left( \left( \Omega_0 - \Omega_1\right)^T \nabla f(\Omega_0)\right)/2 $  \;
$\mathcal{S}_{1}  \leftarrow \left\{ \Omega_1 \right\}$  \;
\For{$k = 1$ to $K$}{
\tikzmk{A}
Find $\Omega_{k+1}$ such that $g_{\co \left( \mathcal{S}_{k}\right)} (\Omega_{k+1}) \leq \Phi$ \label{alg:approx_minimization} \Comment*[r]{Solve problem approximately}
\tikzmk{B} \boxit{pink}
\tikzmk{A}
$\nabla f \left( \Omega_{k+1} \right) \leftarrow 2\Psi(Y) \left(\dot{Y} - \Omega_{k+1}^T\Psi(Y) \right)^T$ \;
$V_{k+1} \leftarrow \argmin\limits_{\Omega \in  \mathcal{P}} \trace\left(\Omega^T\nabla f \left( \Omega_{k+1} \right) \right)$ \label{alg:linear_problemBCG}\;
$g_{\mathcal{P}} \left( \Omega_{k+1}\right) \leftarrow\trace \left( \left( \Omega_{k+1} - V_{k+1}\right)^T \nabla f(\Omega_{k+1})\right)$ \label{alg:FW_gapBCG}\;
\uIf{$g_{\mathcal{P}} \left( \Omega_{k+1}\right) \leq \Phi$}{
\tikzmk{C}
$\Phi \leftarrow g_{\mathcal{P}} \left( \Omega_{k+1}\right)/2$ \label{alg:update_accuracy_BCG}\Comment*[r]{Update accuracy} 
\tikzmk{D} \boxitshort{lightblue}
\tikzmk{C}
$\mathcal{S}_{k+1} \leftarrow \mathcal{S}_k$ \;
$\Omega_{k+1} \leftarrow \Omega_k$ \;
}
\Else{
\tikzmk{C}
$\mathcal{S}_{k+1} \leftarrow \mathcal{S}_k \bigcup V_{k+1}$ \label{alg:expand_active_set_BCG}\Comment*[r]{Expand active set}
\tikzmk{D} \boxitshort{lightgreen}
\tikzmk{C}
$D_k \leftarrow V_{k + 1} - \Omega_k$\label{alg:descent_direction_BCG} \;
$\gamma_k \leftarrow \min\left\{-\frac{1}{2}\trace \left( D_k^T \nabla f \left(
    \Omega_k \right) \right)/ \norm{D_k^T
  \Psi(Y)}_F^2,1\right\}$ \label{alg:stepSize} \;
$\Omega_{k+1} \leftarrow \Omega_k + \gamma_k D_k$ \label{alg:update_BCG}
}
}
\caption{\texttt{CINDy}: Blended Conditional Gradient (BCG) algorithm variant applied to Problem~\eqref{eq:l1_minimization_noisy_additional_constraints}} \label{algo:BCG}
\end{algorithm}

The algorithm \emph{approximately} minimizes $f\left( \Omega \right)$
over $\co(\mathcal{S}_k)$ in Line~\ref{alg:approx_minimization} of
Algorithm~\ref{algo:BCG}. This problem is analogous to the one shown
in Equation~\eqref{eq:unit_probability_simplex} and can be solved in
the space of $\lambda$ barycentric coordinates. The approximate
minimization is carried out until the Frank-Wolfe gap satisfies
$g_{\co \left( \mathcal{S}_{k}\right)} (\Omega_{k+1}) \leq \Phi$. The
algorithm then computes the Frank-Wolfe gap over $\mathcal{P}$ in
Lines~\ref{alg:linear_problemBCG}-\ref{alg:FW_gapBCG}, that is
$g_{\mathcal{P}}\left(\Omega_{k+1}\right)$. If this is smaller than
the accuracy $\Phi$ to which we are computing the solutions in
Line~\ref{alg:approx_minimization}, we increase the accuracy to which
we compute the solutions in Line~\ref{alg:update_accuracy_BCG} by
taking $\Phi = g_{\mathcal{P}}\left(\Omega_{k+1}\right)/2$. This means
that as we get closer to the solution of the optimization problem, and
the gap $g_{\mathcal{P}}\left(\Omega_{k+1}\right)$ decreases, we
increase the accuracy to which we solve the problems over
$\co(\mathcal{S}_{k})$. If on the other hand
$g_{\mathcal{P}}\left(\Omega_{k+1}\right)$ is larger than $\Phi$,
expanding the active set promised more progress than continuing
optimizing over
$\co(\mathcal{S}_k)$. Thus, we potentially expand the active set in
Line~\ref{alg:expand_active_set_BCG}, and we perform a standard CG
step with exact line search in
Lines~\ref{alg:descent_direction_BCG}-\ref{alg:update_BCG}. Regarding
the step size in Line~\ref{alg:line_search} the same comments apply as
in Section~\ref{sec:frank-wolfe-algor}: it is always non-negative,
ensures convex combinations, and the clipping, if active, is active
only in the very first iteration.

The BCG algorithm enjoys robust theoretical convergence guarantees,
and exhibits very fast convergence in practice. Moreover, it generally
produces solutions with a high level of sparsity in the experiments
essentially identical to those produced by FCCG. This makes the BCG
algorithm a powerful alternative to the sequentially-thresholded
least-squares approach followed in the SINDy algorithm
\cite{brunton2016discovering} (or the sequentially-thresholded ridge
regression in \citet{rudy2017data}).

\section{Numerical experiments} \label{Section:Numerical_experiments}

We benchmark the CINDy framework (Algorithm~\ref{algo:BCG}) using the LASSO problem formulations presented in
Equation~\eqref{eq:l1_minimization_noisy} with the following
algorithms. Our main benchmark here is the SINDy framework, however we
included three other popular optimization methods for further
comparison. 

\paragraph{\textbf{SINDy}:} We use a SINDy framework implementation based on the Python \href{https://github.com/snagcliffs/PDE-FIND}{\texttt{PDE-FIND}} Github repository from \cite{rudy2017data} (which originally used ridge-regression, as opposed to the least-squares regression used in \cite{brunton2016discovering}).

\paragraph{\textbf{SR3}:} We use a Sparse Relaxed Regularized Regression (SR3) framework implementation based on the Python \href{https://github.com/kpchamp/SINDySR3}{\texttt{SINDySR3}} Github repository from \cite{champion2020unified}, with some modifications. As this framework admits regularization through the $\ell_0$ and the $\ell_1$ norm we test against both. We only show the results for the constrained version of the problem (where we impose additional structure), as we achieved the best performance with those. Note that for this framework we have to tune both the strength of the regularization, and the relaxation parameter, namely $\alpha$ and $\nu$ in Equation~\eqref{eq:SR3_problem}.

\paragraph{\textbf{FISTA}:} The \emph{Fast Iterative Shrinkage-Thresholding Algorithm} \cite{beck2009fast}, commonly known as FISTA, is a first-order accelerated method commonly used to solve LASSO problems which are equivalent to the ones shown in Equation~\eqref{eq:l1_minimization_noisy} (in which the $\ell_1$ norm appears as a regularization term in the objective function, as opposed to a constraint).
\paragraph{\textbf{IPM}:} \emph{Interior-Point Methods} (IPM) \cite{nesterov1994interior} are an extremely powerful class of convex optimization algorithms, able to reach a highly-accurate solution in a small number of iterations. The algorithms rely on the resolution of a linear system of equations at each iteration, which can be done efficiently if the underlying system is sparse. Unfortunately this is not the case for our LASSO formulations, which makes this algorithm impractical for large problems. We will use the path-following primal form interior-point method for quadratic problems described in \cite{andersen2011interior}, and implemented in Python's \href{https://cvxopt.org/}{\texttt{CVXOPT}}, to solve the LASSO problem in Equation~\eqref{eq:l1_minimization_noisy} (with and without additional constraints, as described in Section~\ref{section:IncorporatingStructure}).

We use CINDy (c) and CINDy to refer to the results achieved by the CINDy framework with and without the additional constraints described in Section~\ref{section:IncorporatingStructure}. Likewise, we use IPM (c), IPM, SR3 (c-$\ell_0$) and SR3 (c-$\ell_1$) to refer to the results achieved by IPM with and without additional constraints, and SR3 with constraints using the $\ell_0$ and $\ell_1$ regularization, respectively. We have not added structural constraints to the formulation in FISTA, as we would need to compute non-trivial proximal/projection operators, making the algorithm computationally expensive.

\begin{remark}[Hyperparameter selection for the CINDy framework]
  In the experiments we have not tuned the $\ell_1$ paramater in the
  LASSO formulation for the CINDy algorithm (neither in the integral
  nor the differential formulation), simply relying on
  $\alpha = 2\norm{\dot{Y}\Psi(Y)^{\dagger}}_{1,1}$ and
  $\alpha = 2\norm{\delta Y\Gamma(Y)^{\dagger}}_{1,1}$ in the
  differential and integral formulation for all the experiments. With
  this choice, purposefully, all computed solutions are located in the
  interior of the feasible region. It is important to note that due to
  this choice, \emph{sparsity in the recovered dynamics is therefore due to the implicit
  regularization by the optimization algorithm used in the CINDy framework, namely BCG} and not due to
  binding constraints of the LASSO problem formulation.
\end{remark}

\begin{remark}[Hyperparameter selection for SINDy, SR3, FISTA and IPM]
  We have selected the threshold coefficient for SINDy, the $\ell_1$ regularization and relazation parameters of SR3, the $\ell_1$ regularization parameters of FISTA, and the IPM algorithm based on performance on validation data. In the differential and integral formulations we have selected the hyperparameters that gave the smallest value of $\norm{\dot{Y}_{\text{validation}} - \Omega^T \Psi(Y_{\text{validation}})}^2_F$ and $\norm{\delta Y_{\text{validation}} - \Omega^T \Gamma(Y_{\text{validation}})}^2_F$ respectively. More concretely, we have used either \href{http://hyperopt.github.io/hyperopt/}{\texttt{Hyperopt}} \citep{bergstra2013making} or \href{https://github.com/scipy/scipy}{\texttt{Scipy's}} \citep{2020SciPy-NMeth} minimization functions to find the best set of hyperparameters for each algorithms, where we find the parameters that minimize the validation loss.
  % \todo[inline]{MW: Just as a comment -- this might be seen as unfair parametrization if we later claim that CINDy produces sparser results. Sparsity has simply not been a criterion for steering the competition, whereas this might be the case when selecting the CINDy stopping criterion. On the other hand, CINDy shows slightly superior $\mathcal{E}_R$ accuracy, which suggests that this is not the main point.}
  Other criteria could be chosen to increase the level of sparsity of
  the solution, at the expense of accuracy in inferring
  derivatives/trajectories. This would involve deciding how to weight the accuracy and the sparsity of the returned solutions when scoring these solutions, and it is unclear how to do so. Note that in the BCG algorithm, the
  sparsity-accuracy compromise is instead parametrized in terms of the
  stopping criterion instead of thresholds or $\ell_1$ bounds. We
  hasten to stress however, that the sparsity of CINDy (and more
  precisely of the BCG algorithm) is due to \emph{extremely sparse updates} in
  each iteration, whereas some of the other algortihms' updates are naturally
  \emph{dense} and sparsity is only realized by means of
  postprocessing these updates by using sparse regularization techniques.
\end{remark}

\begin{remark}[Stopping criterion]
  We use $200$ rounds of thresholding and least-squares for each run of the SINDy algorithm (or until no more coefficients are thresholded in a given iteration). For the SR3 algorithm, we use the existing stopping criterion in the original implementation from \cite{champion2020unified}. The threshold value that yields the best accuracy in terms of testing data is outputted afterwards. During the selection of the FISTA hyperparameter, the algorithm is run for a sufficiently large number of iterations until the primal progress made is below a tolerance of $10^{-8}$. The same is done when we run the FISTA algorithm with the final hyperparameter selected. The IPM algorithm is run with the default stopping criterion parameters. Lastly, the CINDy algorithm is run until the Frank-Wolfe gap, an upper bound on the primal gap, is below a tolerance of $10^{-6}$.
\end{remark}

For each physical model in this section we generate a set of $T$ points, simulating the physical model $c$ different times, to generate $c$ experiments. First, a random starting  point for the $j$-th experiment is generated. This random starting point is used to generate a set of $T/c$ points equally-spaced in time $\vx^j(t_i)$ with $i \in \llbracket 1, T/c \rrbracket$ using a high-order Runge-Kutta scheme and ensuring that the discretization error $\|\vx^j(t_i) - \vx(t_i)\|$ is below a tolerance of $10^{-13}$. The samples are then contaminated with i.i.d.\ Gaussian noise. If we denote the noisy data point at time $t_i$ for the $j$-th experiment by $\vy^j(t_i)$, we have 
\begin{align*}
    \vy^j(t_i) = \vx^j(t_i) + \eta \mathcal{N}\left(0, \Sigma \right),
\end{align*}
where $\mathcal{N}\left(0, \Sigma \right)$ denotes the $d$-dimensional multivariate Gaussian distribution centered at zero with covariance matrix $\Sigma$, where
\begin{gather}
    \Sigma = \diag\left(\frac{1}{T} \sum_{i = 1}^{T/c}\sum_{j=1}^c\left(x^j_1(t_i) - \mu_1 \right)^2, \cdots,\frac{1}{T} \sum_{i = 1}^{T/c}\sum_{j=1}^c\left(x^j_d(t_i) - \mu_d \right)^2\right) \\
    \mu_k = \frac{1}{T} \sum_{i = 1}^{T/c}\sum_{j=1}^c x^j_k(t_i).
\end{gather}
We denote the noise level that we vary in our experiments by $\eta$. Note that the $k$-th element on the diagonal of $\Sigma$ is simply the sample variance of the $k$-th component of points generated in the experiments.

 Both $\Gamma (Y)$ and $\Psi (Y)$ are normalized so that their rows
 have unit variance, to make the learning process easier for all the
 algorithms. For the training of the algorithm, $70\%$ of the data
 points are used, while $20\%$ are used for validation and selecting
 the combination that gives the best hyperparameters (the $\ell_1$
 radius in the FISTA and IPM algorithms, or the threshold in SINDy's
 sequential thresholded least-squares algorithm). Lastly, the $10\%$
 remaining data points are used for evaluating the output of each
 algorithm, and are referred to as the testing set.

 We would like to stress that, while it might seem that we are in
 the overdetermined regime in terms of the number of samples used in
 training, this is not accurate as due to the evolutionary nature,
 i.e., evolving the dynamic in time, for each of the experiments the
 samples obtained \emph{within} one experiment are highly correlated.

The approximate derivatives are computed from noisy data using local polynomial interpolation \cite{knowles2014methods}. The matrix $\Gamma (Y)$ used in the integral formulation of the sparse recovery problem was computed through the integration of local polynomial interpolations. The same matrices and training-validation-testing split is used in each experiment for all the algorithms, to make a fair comparison. Section~\ref{appx:derivs_and_ints} in the Appendix shows the difference in accuracy that can be achieved with the different methods when estimating the derivatives and integrals for one physical model. Given the disparity in accuracy that can be achieved between the estimation of derivatives and integrals, we have decided to present the results for the integral and differential formulation separately.

\subsection{Benchmark metrics}

We benchmark the algorithms in terms of the following metrics, in a
similar spirit as was done in \cite{kaheman2020automatic}. Given a
dictionary of basis functions $\mathcal{D}$ of cardinality $n$, an
associated exact dynamic $\Xi\in \rr^{n \times d}$ such that
$\dot{\vx}(t) = \Xi^T \bm{\psi}\left( \vx(t)\right)$, and a dynamic
$\Omega\in \rr^{n \times d}$ outputted by the algorithms, we define:
\begin{definition}[Recovery error] The \emph{recovery error} $\mathcal{E}_{R}$ of the algorithm is given by:
\begin{align*}
    \mathcal{E}_{R} \defeq \norm{\Omega - \Xi}_F.
\end{align*}
\end{definition}
Given noisy data points $Y_\text{testing}, \dot{Y}_\text{testing}$ from the testing set we define:
\begin{definition}[Derivative inference error] The \emph{derivative inference error} $\mathcal{E}_{D}$ of the algorithm is given by:
\begin{align*}
    \mathcal{E}_{D} \defeq \norm{(\Omega - \Xi)^T \Psi\left(Y_\text{testing}\right)}_F.
\end{align*}
This measure aims at quantifying how well the learned dynamics will infer the true derivatives at $Y$.
\end{definition}

\begin{definition}[Trajectory inference error] The \emph{trajectory inference error} $\mathcal{E}_{T}$ of the algorithm is given by:
\begin{align*}
    \mathcal{E}_{T} \defeq \norm{ (\Omega- \Xi)^T \Gamma\left(Y_\text{testing}\right)}_F.
\end{align*}
This measure aims at quantifying how well the learned dynamics will infer the trajectory using the approximate matrix $\Gamma\left(Y\right)$, compared to the true dynamic $\Xi$.
\end{definition}

In order to gauge how well a given algorithm is able to recover the true support of a dynamic, we define:

\begin{definition}[Extraneous terms] The \emph{extraneous terms} of a given dynamic $\Omega\in\rr^{d\times n}$ with respect to its true counterpart $\Xi\in\rr^{d\times n}$ is defined as:
\begin{align*}
\mathcal{S}_E  \defeq \big| \left\{ \Omega_{i,j} \mid \Omega_{i,j} \neq 0, \Xi_{i,j} = 0, i \in \llbracket 1,d \rrbracket,  j \in \llbracket 1,n \rrbracket\right\} \big|.
\end{align*}
This metric simply counts the terms picked up in $\Omega$ that are not present in the true physical model, represented by $\Xi$, i.e.\ it counts the false positives.
\end{definition}

\begin{definition}[Missing terms] The \emph{missing terms} of a given dynamic $\Omega\in\rr^{d\times n}$ with respect to its true counterpart $\Xi\in\rr^{d\times n}$ is defined as:
\begin{align*}
\mathcal{S}_M  \defeq \big| \left\{ \Omega_{i,j} \mid \Omega_{i,j} = 0, \Xi_{i,j} \neq 0, i \in \llbracket 1,d \rrbracket,  j \in \llbracket 1,n \rrbracket\right\} \big|.
\end{align*}
This metric simply counts the terms that have not been picked up in $\Omega$ that actually participate in governing the true physical model, represented by $\Xi$, i.e.\ it counts the false negatives.
\end{definition}

\subsection{Kuramoto model} \label{section:kura}
\begin{figure*}[]
    \centering
    \hspace{\fill}
    \subfloat{{\includegraphics[width=15.6cm]{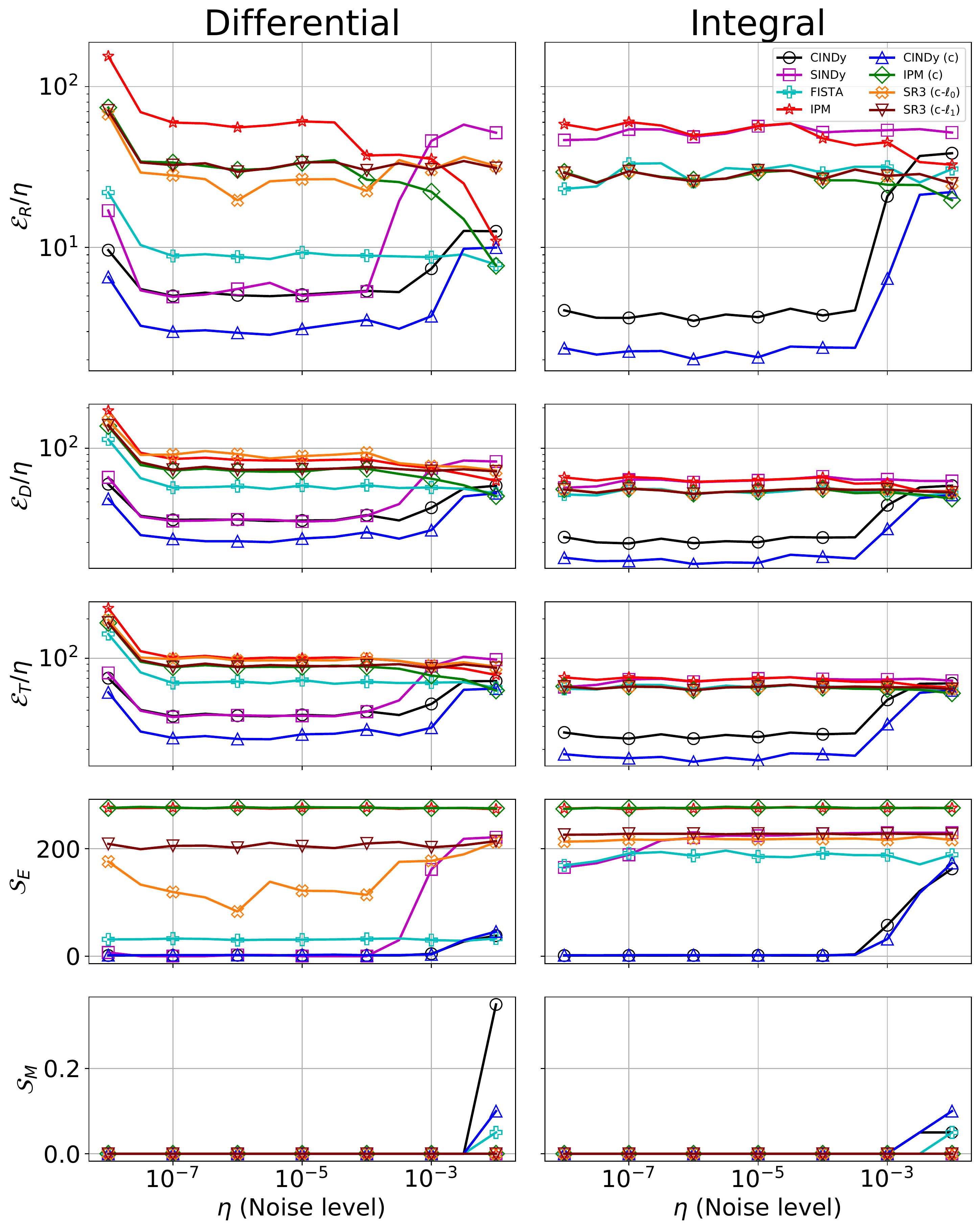} }}%
    \hspace{\fill}
    \caption{\textbf{Sparse recovery of the Kuramoto model: } Algorithm comparison for a Kuramoto model of dimension $d=5$, with a differential formulation shown on the left column, and with an integral formulation on the right column. The first, second, third, fourth and fifth rows of images indicate a comparison of $\mathcal{E}_R/\eta$, $\mathcal{E}_D/\eta$, $\mathcal{E}_T/\eta$, $\mathcal{S}_E$, and $\mathcal{S}_M$, as we vary the noise level $\eta$, respectively.}%
    \label{fig:kuramoto5_v5}%
\end{figure*}

The Kuramoto model describes a large collection of $d$ weakly coupled identical oscillators, that differ in their natural frequency $\omega_i$ \cite{kuramoto1975self} (see Figure~\ref{fig:kuramoto5_v5}). This dynamic is often used to describe synchronization phenomena in physics, and has been previously used in the numerical experiments of a tensor-based algorithm for the recovery of large dynamics \cite{gelss2019multidimensional}. If we denote by $x_i$ the angular displacement of the $i$-th oscillator, then the governing equation with external forcing (see \cite{acebron2005kuramoto}) can be written as:
\begin{align*}
  \dot{x}_i & =  \omega_i + \frac{K}{d}\sum_{j=1}^d \sin \left( x_j - x_i\right) + h\sin \left( x_i\right) \\
  & =  \omega_i + \frac{K}{d}\sum_{j=1}^d \left[\sin \left( x_j \right) \cos \left( x_i \right) - \cos \left( x_j \right) \sin \left( x_i \right) \right]+ h\sin \left( x_i\right),
\end{align*}
for $i\in \llbracket 1,d \rrbracket$, where $d$ is the number of oscillators (the dimensionality of the problem), $K$ is the coupling strength between the oscillators and $h$ is the external forcing parameter. The exact dynamic $\Xi$ can be expressed using a dictionary of basis functions formed by sine and cosine functions of $x_i$ for $i\in \llbracket 1,d\rrbracket$, and pairwise combinations of these functions, plus a constant term. To be more precise, the dictionary used is
\begin{align*}
   \mathcal{D} = \left\{\prod_{i=1}^d \sin(x_i)^{a_i} \prod_{i=1}^d \cos(x_i)^{b_i}  \mid  a_i, b_i \in \llbracket 0, 1 \rrbracket,  i\in \llbracket 1, d \rrbracket, 0 \leq \sum_{i=1}^d (a_i + b_i) \leq 2\right\}.
\end{align*}
Which has a cardinality of $1+d+2d^2$. Note however, that the data is contaminated with noise, and so we observe $\vy$ as opposed to $\vx$. For the system we choose the natural frequency $\omega_i \sim \mathcal{U}[0, 1]$ for $i\in \llbracket 1, d \rrbracket$. The random starting point for each instance of the experimental data used is chosen as $\vx^j(t_0) \sim \mathcal{U}[0, 2\pi]^d$. This starting point is used to generate a trajectory according to the exact dynamic using a high-order Runge-Kutta scheme for a maximum time $t_{T/c}$ of $10$ seconds. This trajectory is then contaminated with noise, in accordance with the description in the previous section.

\begin{remark}
  \label{remark:robustness_full_rank}
Note that if we were to include $\cos(y_i)^2$ and $\sin(y_i)^2$ for $i\in \llbracket 1,d\rrbracket$ in $D$, the matrix $\Psi(Y)$ built with this library would not have full rank, as $\cos(y_i)^2 + \sin(y_i)^2 = 1$, and the library $D$ already includes a built-in constant. In our experiments we have observed that when using a dictionary that does not have full rank, the thresholded least-squares algorithm SINDy tends to produce solutions that include constant, $\cos(y_i)^2$ and $\sin(y_i)^2$ terms, whereas the dynamics returned by the CINDy algorithm tends to only include constant terms, thereby providing a more parsimonious representation of the dynamic. If we add an $\ell_2$ regularization term to the optimization algorithm, resulting in a thresholded ridge regression algorithm, the resulting dynamic tends toward higher parsimony, but at the expense of accuracy in predicting derivatives and trajectories.
\end{remark}

We also test the performance of the CINDy algorithm and the IPM algorithm with the addition of symmetry constraints. We use $\xi_j\left( \psi(\vx)\right)$ to refer to the coefficient in $\xi_j$, where $\xi_j$ is the $j$-th column of $\Xi$, associated with the basis function $\psi(\vx)$. The underlying rationale behind the constraints is that as the particles are identical, except for their intrinsic frequency, the effect of $x_i$ on $x_j$ should be the same as the effect of $x_j$ on $x_i$. In both the integral and the differential formulation we impose that for all $i, j \in \llbracket 1, d\rrbracket$:
\begin{gather*}
    \xi_j\left(\sin \left( x_i \right)\right) = \xi_i\left(\sin \left( x_j \right)\right) \\
    \xi_i\left(\cos \left( x_i \right)\right) = \xi_i\left(\cos \left( x_j \right)\right) \\
    \xi_j\left(\cos \left( x_i \right)\right) = \xi_i\left(\cos \left( x_j \right)\right) \\
     \xi_j\left(\sin \left( x_i \right)\cos \left( x_j \right)\right) = \xi_i\left(\sin \left( x_j \right) \cos \left( x_i \right)\right) \\
      \xi_j\left(\cos \left( x_i \right)\sin \left( x_j \right)\right) = \xi_i\left(\cos \left( x_j \right) \sin \left( x_i \right)\right) \\
       \xi_j\left(\sin \left( x_i \right)\sin \left( x_j \right)\right) = \xi_i\left(\sin \left( x_j \right) \sin \left( x_i \right)\right) \\
      \xi_j\left(\cos \left( x_i \right)\cos \left( x_j \right)\right) = \xi_i\left(\cos \left( x_j \right) \cos \left( x_i \right)\right),
\end{gather*}
which are simple linear constraints that can easily be added to the linear optimization oracle used in Line~\ref{alg:linear_problemBCG} of the CINDy algorithm (Algorithm~\ref{algo:BCG}). 

\begin{figure*}[]
    \centering
    \hspace{\fill}
    \subfloat{{\includegraphics[width=15.6cm]{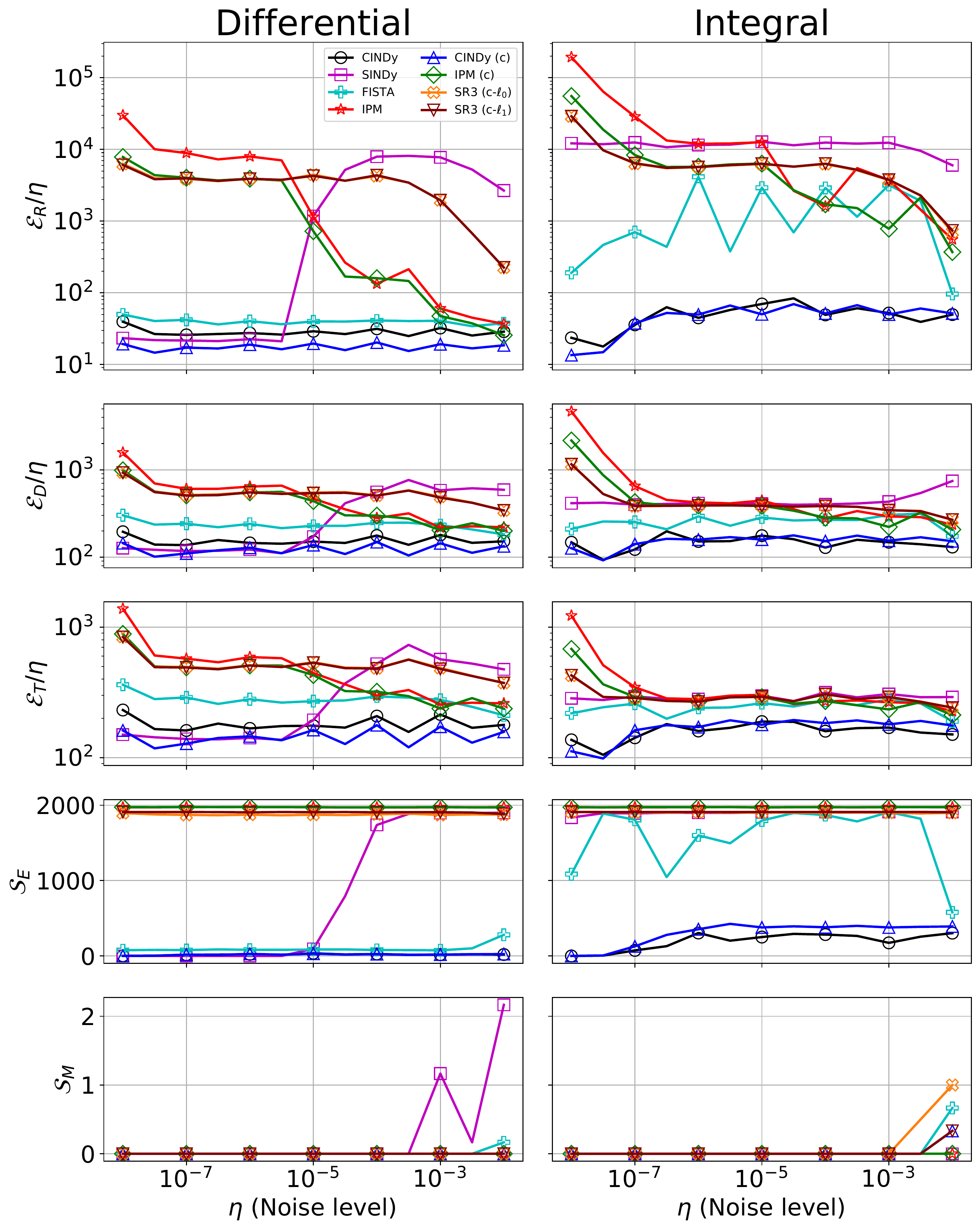} }}%
    \hspace{\fill}
    \caption{\textbf{Sparse recovery of the Kuramoto model: } Algorithm comparison for a Kuramoto model of dimension $d=10$, with a differential formulation shown on the left column, and with an integral formulation on the right column. The first, second, third, fourth and firth rows of images indicate a comparison of $\mathcal{E}_R/\eta$, $\mathcal{E}_D/\eta$, $\mathcal{E}_T/\eta$, $\mathcal{S}_E$, and $\mathcal{S}_T$, as we vary the noise level $\eta$, respectively.}%
    \label{fig:kuramoto10_v5}%
\end{figure*}

The images in Figure~\ref{fig:kuramoto5_v5} and \ref{fig:kuramoto10_v5} show the recovery results for $K = 2$ and $h = 0.2$ and two different values for the dimension, $d = 5$ and $d = 10$, respectively. A total of $6000$ points were used to infer the dynamic, spread over $40$ experiments for a maximum time of $10$ seconds, for both cases. The values of the noise level $\eta$ ranged from $10^{-8}$ to $10^{-2}$. The derivatives used where computed using differentiation of local polynomial approximations, and the integrals using integration of local polynomial approximations. Each test was performed $20$ times. The graphs indicate with lines the average value obtained for $\mathcal{E}_{R}/\eta$, $\mathcal{E}_{D}/\eta$, $\mathcal{E}_{T}/\eta$, $\mathcal{S}_E$ and $\mathcal{S}_E$ for a given noise level and algorithm. When plotting we divide the recovery error $\mathcal{E}_{R}$, the derivative inference error $\mathcal{E}_{D}$, and the trajectory inference error $\mathcal{E}_{T}$ by the noise level $\eta$ in order to better visualize the different performance of the algorithms.

For the case with $d = 5$ (see Figure~\ref{fig:kuramoto5_v5}) we can observe that in the differential formulation the CINDy, CINDy (c) and SINDy frameworks achieve the smallest recovery error $\mathcal{E}_{R}$ for noise levels below $10^{-4}$, and the performance of the SINDy framework degrades after that, relative to that of the best-performing algorithms. Whereas the CINDy and CINDy (c) frameworks are among the best performing in terms of $\mathcal{E}_{R}$ for all noise levels. Regarding the $\mathcal{E}_{D}$ and $\mathcal{E}_{T}$ errors, we can see very similar results as to those found in the $\mathcal{E}_{R}$ images. Regarding the sparsity achieved through the algorithms, we can see that CINDy  and CINDy (c) consistently tend to produce the sparsest solutions, as measured with $\mathcal{S}_E$, whereas IPM and IPM (c) tend to pick up all the terms in the dictionary. For interior point methods, dense solutions are of course to be expected if no solution rounding is performed. Note also that the performance of SINDy degrades above a noise level of $10^{-4}$, where it starts to produce dense solutions. Lastly, on average none of the algorithms tends to miss more than one of the basis functions that are present in the exact dynamic, as measured by $\mathcal{S}_M$, this is especially remarkable for the CINDy and CINDy (c) frameworks, which consistently have the lowest number of extra terms. Overall, in the differential formulation experiments we observe that the CINDy and CINDy (c) frameworks produce the most accurate solutions, as measured by $\mathcal{E}_{R}$, consistently producing among the best dynamics for all the noise levels, while producing dynamics that are much sparser than the ones produced by the other algorithms being considered.

In terms of integral formulation, the CINDy and CINDy (c) frameworks
are on average more accurate in terms of $\mathcal{E}_{R}$ than any of
the algorithms tested for noise levels below $10^{-3}$, by a larger margin than in the differential
formulation. For noise levels above $10^{-3}$ the results start to look somewhat similar for all the algorithms. Regarding the $\mathcal{E}_{D}$ and $\mathcal{E}_{T}$
error, we can see that CINDy and CINDy (c) produce the most accurate solutions in terms of inferring derivatives or trajectories for noise levels below $10^{-3}$. Note that in this case there is also a large difference in sparsity between the algorithms, as all the algorithms except the CINDy and CINDy (c) frameworks tend to pick up a large number of extraneous basis
functions.

For the case with $d = 10$ (see Figure~\ref{fig:kuramoto10_v5}) in the
differential case, we can see very good performance from SINDy,
CINDy, CINDy (c), and FISTA in $\mathcal{E}_R$ for noise
levels below $10^{-5}$. However, the performance for SINDy degrades
above $10^{-5}$, in terms of recovery accuracy (there is a difference
of more than two orders of magnitude between SINDy and CINDy), and in
terms of $\mathcal{S}_E$, as the algorithm starts picking up most of
the available basis functions from the dictionary. Regarding the integral
formulation, there is a large difference in the performance of CINDy
and CINDy (c), and the rest of the frameworks. For most noise levels
these two aforementioned algorithms are more than two orders of
magnitude more accurate in terms of $\mathcal{E}_R$, while maintaining
extremely sparse solutions. This performance difference is key in
accurately simulating out-of-sample trajectories. We also remark on the fact that the IPM, IPM (c), SR3 (c-$\ell_0$) and SR3 (c-$\ell_1$) improve in terms of the ratio $\mathcal{E}_R/\eta$ as we increase $\eta$, remaining above the CINDy and CINDy (c) frameworks however. We also note that in terms of $\mathcal{E}_R$ there is a benefit to the use of additional structural constraints in the formulation, which can be seen when comparing the results of the IPM and IPM (c) runs, and the CINDy and CINDy (c) runs.

We highlight that in the experiments in this section the SINDy framework can provide good performance for low noise levels in some instances (as can be seen in the results for both $d=5$ and $d=10$ in the differential formulation), however as we increase the noise level its performance significantly degrades. Lastly, note that with the SR3 variants we were only able to obtain significantly better performance over SINDy for medium-low noise levels for the experiments with $d = 10$ in terms of $\mathcal{E}_R$ (however obtaining solutions that were not sparse in terms of $\mathcal{S}_E$), whereas the advantage in using SR3 over SINDy is not apparent in the experiments for $d = 5$.

\subsubsection{Sample Efficiency} \label{section:kuramoto:sample_efficiency}

One of the most crucial aspects of many modern learning problems is
the quantity of data needed to train a given model to achieve a
certain target accuracy. When training data is expensive to gather, it
is usually advantageous to use models or frameworks that require the
least amount of training data to reach a given target accuracy on
validation data. As we can see in Figure~\ref{fig:kuramoto5_v5}, in
both the differential and integral formulation all the algorithms
perform similarly when tested on noisy validation data, as is shown in
the second and third rows of images, however, there are disparities in
how they perform when measuring the performance against the exact
dynamic, as seen in the first row of images. For example, in the
differential formulation the CINDy and CINDy (c) algorithms perform
noticeably better than the SINDy algorithm for higher noise levels,
from $10^{-4}$ to $10^{-2}$, and in the integral formulation the CINDy
and CINDy (c) algorithms perform noticeably better than the SINDy
algorithm for noise levels below $10^{-3}$. From a sample efficiency
perspective, this suggests that in both these regimes where CINDy has
an advantage, it will require fewer samples than SINDy to reach a
target accuracy.

\begin{figure*}[]
    \centering
    \vspace{-10pt}
    \hspace{\fill}
    \subfloat{{\includegraphics[width=8cm]{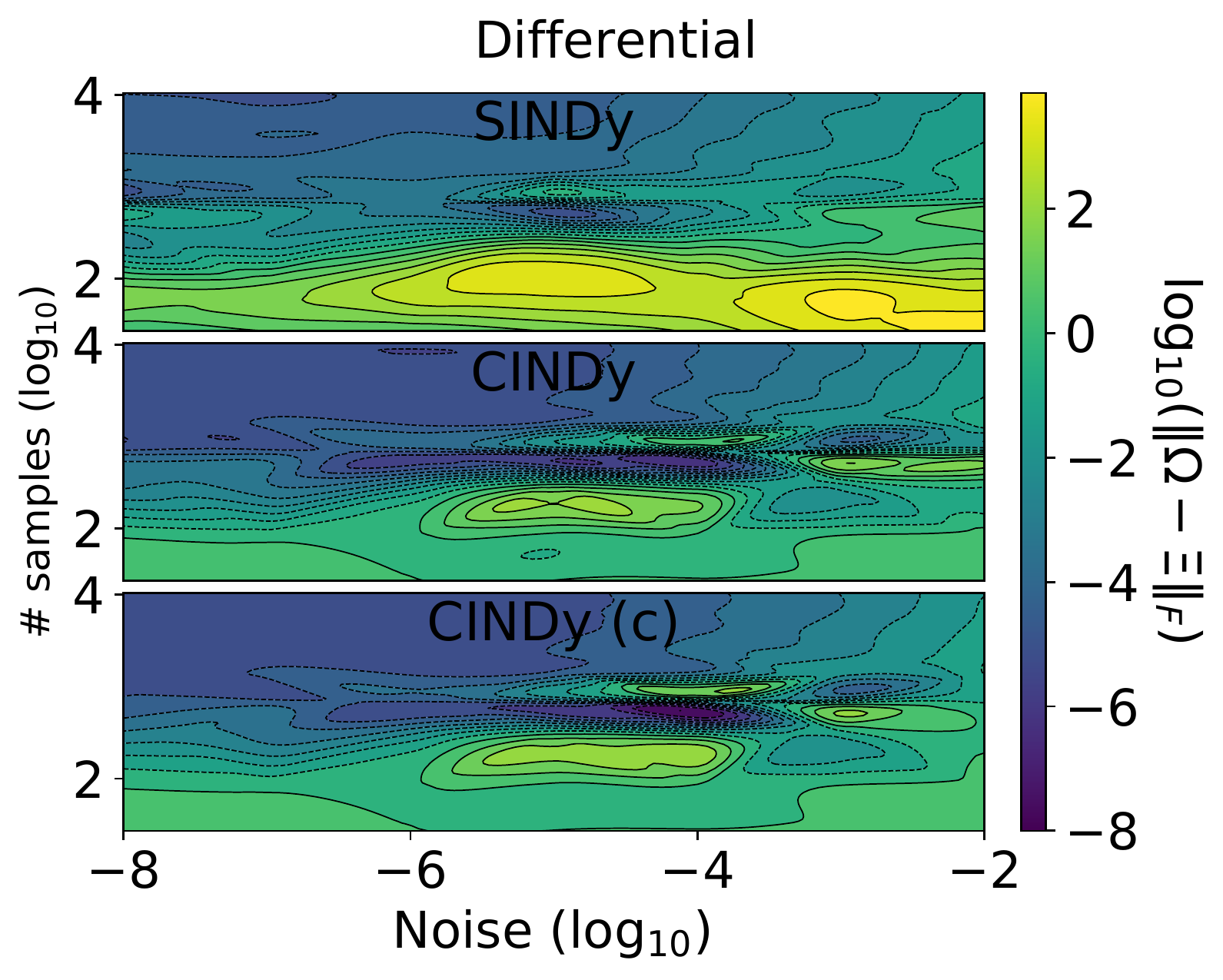} }\label{fig:kuramoto:diff:dim5_sampleeff_acc}}%
    %\qquad
    \hspace{\fill}
    \subfloat{{\includegraphics[width=8cm]{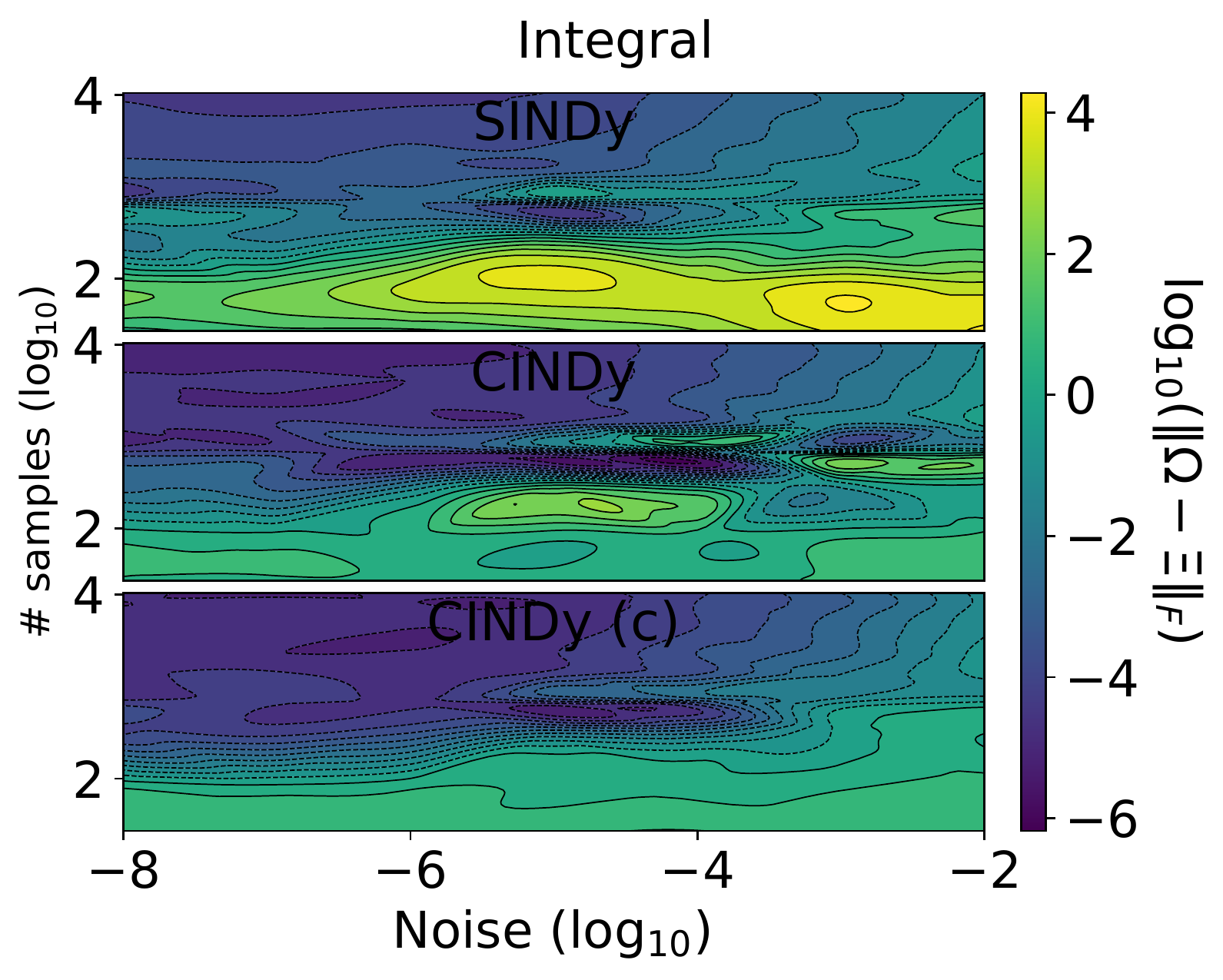} }\label{fig:kuramoto:int:dim5_sampleeff_acc}}%
    \hspace{\fill}
    \caption{\textbf{Sample efficiency of the sparse recovery of the
        Kuramoto model: } Algorithm comparison for a Kuramoto model of
      dimension $d=5$ for both formulations.}%
    \label{fig:kuramoto_sampleeff}%
  \end{figure*}

This is confirmed in Figure~\ref{fig:kuramoto_sampleeff}, which shows
a heat map of $\log(\mathcal{E}_R)$ for different noise levels
(x-axis) and different numbers of training samples (y-axis), when
using the Kuramoto model of dimension $d = 5$ for benchmarking. If one
takes a look at the differential formulation, one can see that in the
low training sample regime both CINDy and CINDy (c) perform better
than SINDy at higher noise levels. It should also be noted that CINDy
(c) performs better than CINDy, which is expected, as the
  introduction of extra constraints \emph{lowers the dimensionality} of our
  learning problem, for which we now have to learn fewer
  parameters. Thus the addition of constraints brings two advantages:
(1) it outputs dynamics that are consistent with the
underlying physics of the phenomenon and (2) it can potentially require
fewer data samples to train. Similar conclusions can be drawn when
inspecting the results for the integral formulation, for example
focusing again on the low training sample regime. We provide an extended
analysis over a broader range of sample sizes in
Appendix~\ref{appx:section:sample_efficiency} for completeness.

\subsubsection{Simulation of learned trajectories} \label{section:kuramoto:trajectories}

By observing the results for $\mathcal{E}_{T}$ in Figure~\ref{fig:kuramoto10_v5} it would seem at first sight that all the
frameworks (except the IPM and IPM (c) for low noise
levels) will perform similarly when inferring trajectories from an
initial position. However, when we simulate the Kuramoto system from a
given initial position, the algorithms have very different
performances. This is due to the fact that while the single point
evaluations might have rather similar errors (on average), which means
nothing else but that they generalize similarly on the specific evaluations
(as expected as this was the considered objective function) they do differ very
much in their \emph{structural generalization behavior}: all
frameworks but CINDy and CINDy (c) pick up wrong terms to explain the
dynamic, which then in the trajectory evolution, due to compounding,
lead to significant mismatches.

\begin{figure}[h!]
	\centering
	\includegraphics[width=14cm]{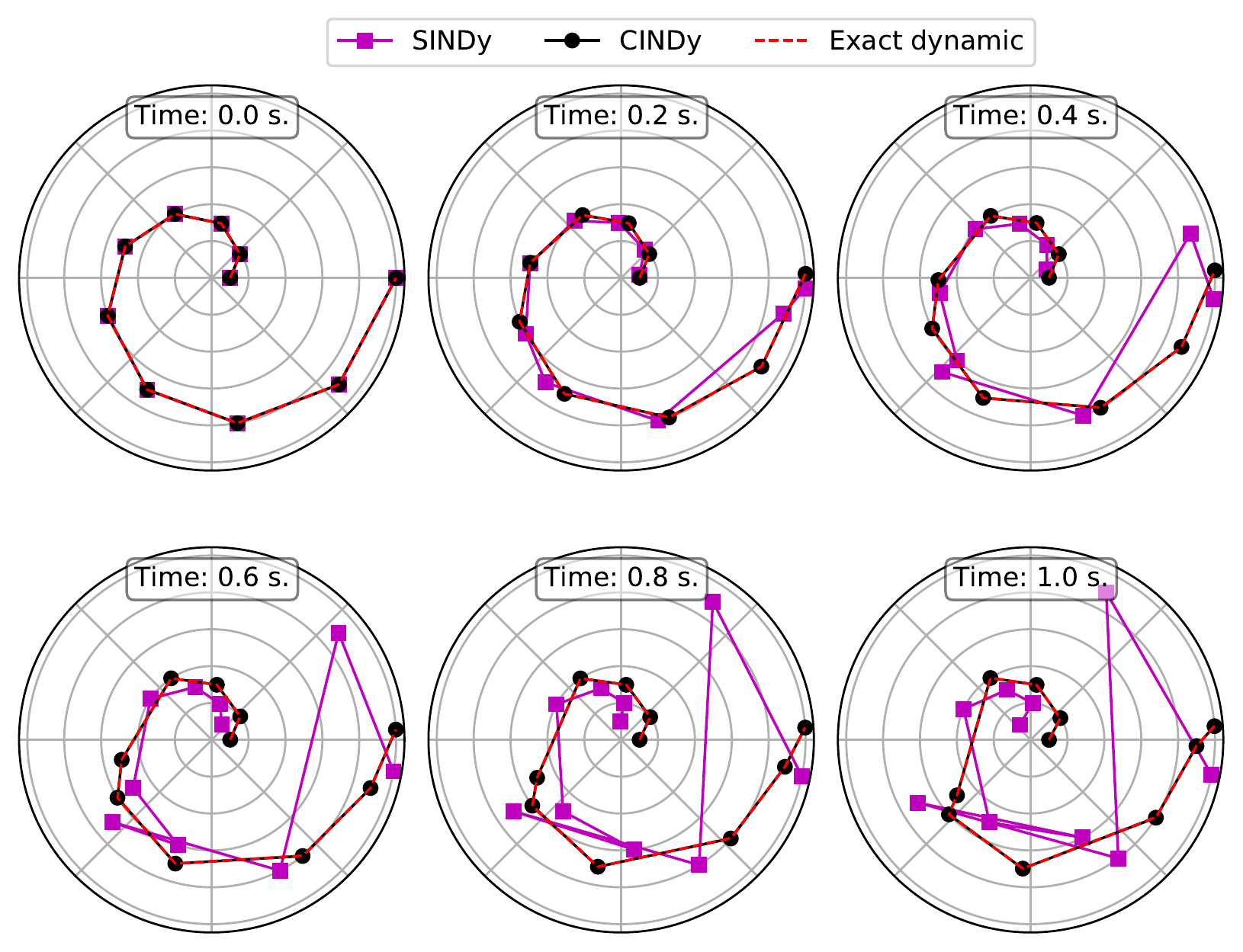} 
	\caption{\textbf{Trajectory comparison:} Kuramoto model of dimension $d = 10$.} \label{kuramoto:trajectories_example}
\end{figure}

In Figure~\ref{kuramoto:trajectories_example} we show the results after simulating the dynamics learned by the CINDy and SINDy algorithm from the integral formulation for a Kuramoto model with $d = 10$ and a noise level of $10^{-3}$. In order to see more easily the differences between the algorithms and the position of the oscillators, we have placed the $i$-th oscillator at a radius of $i$, for $i\in \llbracket 1,d \rrbracket$. This is contrary to how this dynamic is usually visualized, with all the particles oscillating with the same radii.

\subsection{Fermi-Pasta-Ulam-Tsingou model} \label{section:FPUT}

The Fermi-Pasta-Ulam-Tsingou model describes a one-dimensional system of $d$ identical particles, where neighboring particles are connected with springs, subject to a nonlinear forcing term \cite{fermi1955studies}. This computational model was used at Los Alamos to study the behaviour of complex physical systems over long time periods. The prevailing hypothesis behind the experiments was the idea that these systems would eventually exhibit ergodic behaviour, as opposed to the approximately periodic behaviour that some complex physical systems seemed to exhibit. This is indeed the case, as this model transitions to an ergodic behaviour, after seemingly periodic behaviour over the short time scale. This dynamic has already been used in \citet{gelss2019multidimensional}. The equations of motion that govern the particles, when subjected to cubic forcing terms is given by
\begin{align*}
 \ddot{x}_i = \left(x_{i+1} - 2 x_i + x_{i-1} \right) + \beta \left[ \left( x_{i+1} - x_i \right)^3 - \left( x_{i} - x_{i-1} \right)^3 \right],  
\end{align*}
where $i\in \llbracket 1, d\rrbracket$ and $x_{i}$ refers to the displacement of the $i$-th particle with respect to its equilibrium position. The exact dynamic $\Xi$ can be expressed using a dictionary of monomials of degree up to three. To be more precise, the dictionary used is:
\begin{align*}
   \mathcal{D} = \left\{\prod_{i=1}^d x_i^{a_i} \mid  a_i \in \llbracket 0, 3 \rrbracket,  i\in \llbracket 1, n \rrbracket, 0 \leq \sum_{i=1}^d a_i \leq 3\right\},
\end{align*}
which has cardinality $\binom{d + 3}{3}$. Using this dictionary we know that the exact dynamic satisfies $\norm{\Xi}_0 = 10d - 8$. As in the previous example, we can impose a series of linear constraints between the dynamics of neighboring particles (as the particles are identical). In both the integral and the differential formulation we impose that for all $i \in \llbracket 1, d\rrbracket$ and all $j \in \{ i+1, i-1\}$ with $0\leq j \leq d$,
\begin{gather*}
    \xi_j\left( x_i^{a} x_j^b \right) = \xi_i\left( x_i^b x_j^a\right)
\end{gather*}
holds with $a,b \in \llbracket 1,3 \rrbracket$ and $0 \leq a+b+c \leq 3$. 

\begin{remark}[On the construction of $\ddot{Y}$ and $\delta \dot{Y}$]
In this case we are dealing with a second-order ordinary differential equation (ODE), as opposed to a first-order ODE. As we only have access to noisy measurements $\{ \vy(t_i)\}_{i = 1}^{m}$, if we are dealing with the differential formulation, we have to numerically estimate $\{ \ddot{\vy}(t_i)\}_{i = 1}^{m}$, in order to solve
\begin{align*}
\argmin\limits_{\substack{ \Omega \in \mathcal{P}  \\ \Omega \in \rr^{n \times d}}}  \norm{\ddot{Y} - \Omega^T \Psi(Y)}^2_F, %\label{Eq:2ndOrderODE}
\end{align*}
where $\ddot{Y} = \left[ \ddot{\vy}(t_1),\cdots, \ddot{\vy}(t_m)\right] \in\rr^{d\times m}$ is the matrix with the estimates of the second derivatives of $\vy$ with respect to time as columns. On the other hand, if we wish to tackle the problem from an integral perspective, we now use the fact that  $\dot{\vx}(t_{i + 1}) = \dot{\vx}(t_{1}) + \int_{t_1}^{t_{i + 1}} \Xi^T\bm{\psi}(\vx(\tau)) d\tau$, which allows us to phrase the sparse regression problem from an integral perspective as
\begin{align*}
\argmin\limits_{\substack{ \Omega \in \mathcal{P}  \\ \Omega \in \rr^{n \times d}}}  \norm{\delta \dot{Y} - \Omega^T \Gamma(Y)}^2_F, %\label{Eq:2ndOrderODE_integral}
\end{align*}
where $\delta \dot{Y} = \left[ \dot{\vy}(t_2) - \dot{\vy}(t_1),\cdots, \dot{\vy}(t_m) - \dot{\vy}(t_1)\right] \in\rr^{d\times m-1}$ and $\Gamma(Y)$ is computed similarly as in Section~\ref{Section:LearningSparseDynamics}. Note that in this case using the differential formulation requires estimating the second derivative of the noisy data with respect to time to form $\ddot{Y}$, whereas the integral formulation requires estimating the first derivative with respect to time to form $\delta \dot{Y}$.
\end{remark}

\begin{figure*}[]
    \centering
    \hspace{\fill}
    \subfloat{{\includegraphics[width=15.6cm]{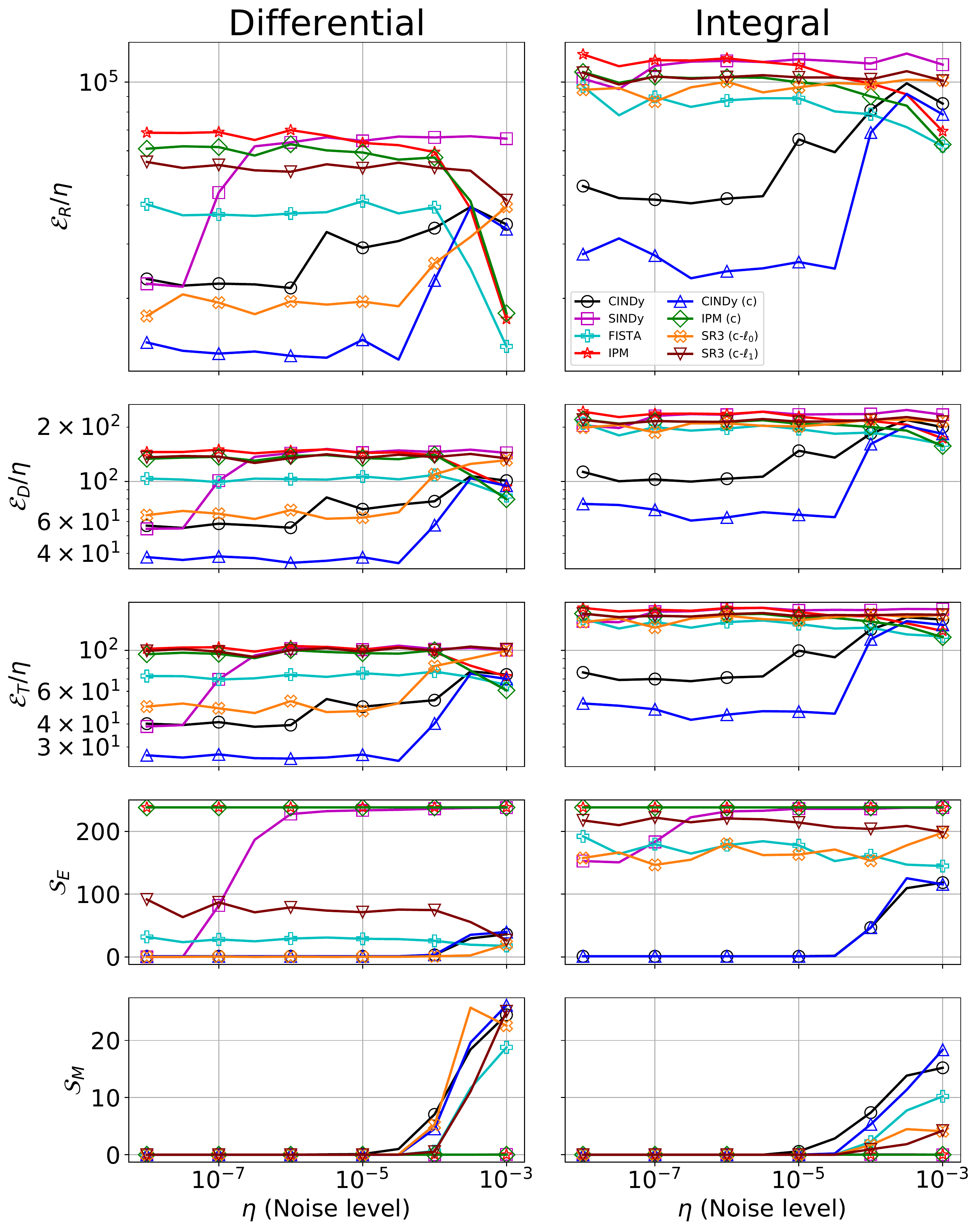} }}%
    \hspace{\fill}
    \caption{\textbf{Sparse recovery of the Fermi-Pasta-Ulam-Tsingou model: } Algorithm comparison for a Fermi-Pasta-Ulam-Tsingou model of dimension $d=5$, with a differential formulation shown on the left column, and with an integral formulation on the right column. The first, second, third, fourth and firth rows of images indicate a comparison of $\mathcal{E}_R/\eta$, $\mathcal{E}_D/\eta$, $\mathcal{E}_T/\eta$, $\mathcal{S}_E$, and $\mathcal{S}_M$, as we vary the noise level $\eta$, respectively.}%
    \label{fig:FPUT5_v5}%
\end{figure*}

\begin{figure*}[]
    \centering
    \hspace{\fill}
    \subfloat{{\includegraphics[width=15.6cm]{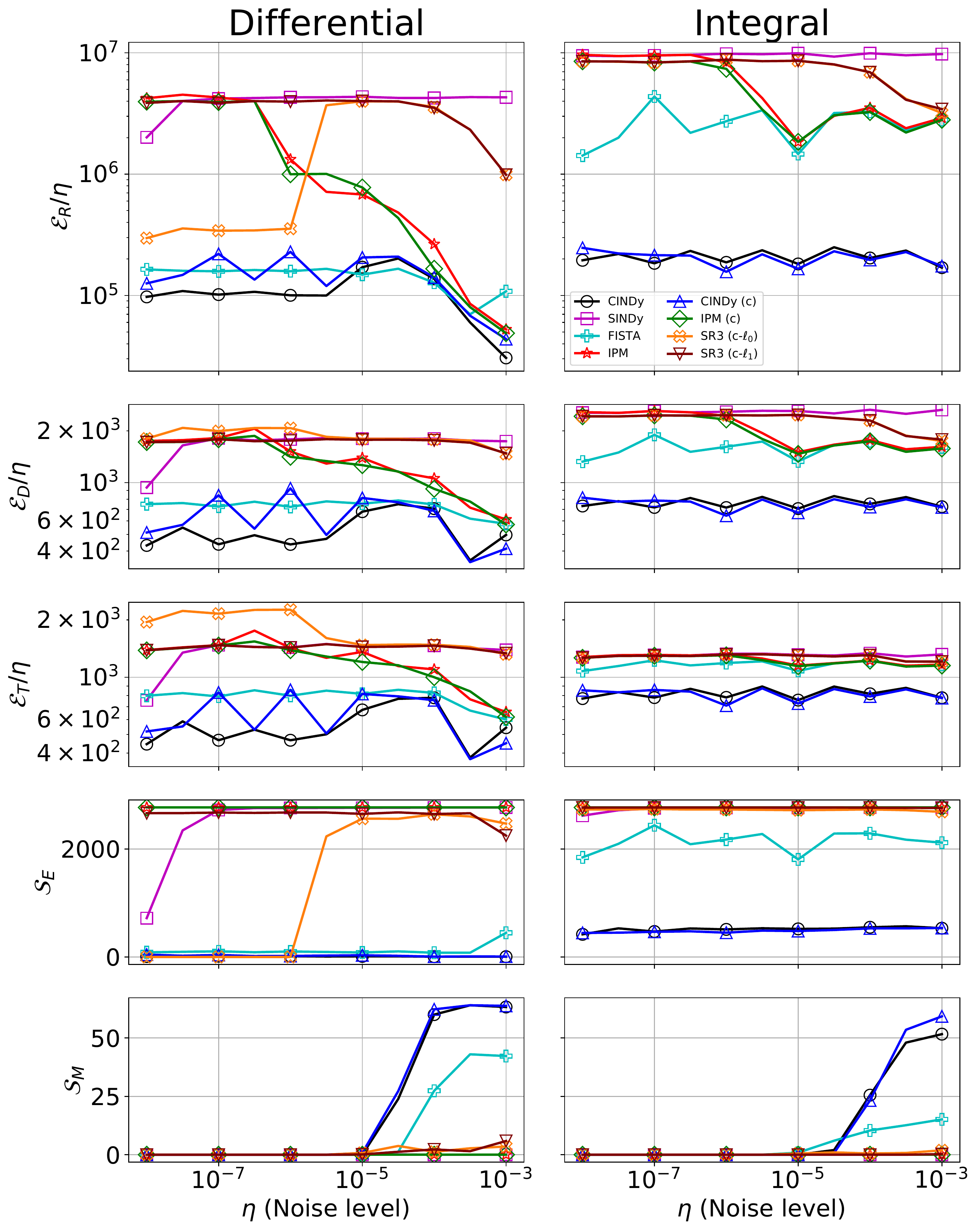} }}%
    \hspace{\fill}
    \caption{\textbf{Sparse recovery of the Fermi-Pasta-Ulam-Tsingou model: } Algorithm comparison for a Fermi-Pasta-Ulam-Tsingou model of dimension $d=10$, with a differential formulation shown on the left column, and with an integral formulation on the right column. The first, second, third, fourth and firth rows of images indicate a comparison of $\mathcal{E}_R/\eta$, $\mathcal{E}_D/\eta$, $\mathcal{E}_T/\eta$, $\mathcal{S}_E$, and $\mathcal{S}_M$, as we vary the noise level $\eta$, respectively.}%
    \label{fig:FPUT10_v5}%
\end{figure*}

The images in Figure~\ref{fig:FPUT5_v5} and \ref{fig:FPUT10_v5} show the recovery results for $d = 5$ and $d = 10$, respectively, when learning the Fermi-Pasta-Ulam-Tsingou with $\beta = 0.7$. The exact dynamic in this case satisfies $\norm{\Xi}_0 = 42$ and $\norm{\Xi}_0 = 92$ for $d = 5$ and $d = 10$, respectively. In the lower dimensional experiment we used $4500$ points, and in the higher dimensional one $9000$. In both cases we spread out the points across $150$ experiments, for a maximum time of $1$ second per experiment. The derivatives used where computed using polynomial interpolation of degree 8. We used polynomial interpolation of degree $8$ to compute the integrals for all the noise levels. The values of the noise level $\alpha$ ranged from $10^{-8}$ to $10^{-2}$. Each test was performed $20$ times. The graphs indicate with lines the average value obtained with $\mathcal{E}_{R}$, $\mathcal{E}_{D}$ and $\mathcal{E}_{T}$, $\mathcal{S}_E$ and $\mathcal{S}_M$ for a given noise level and algorithm. The shaded regions indicate the value obtained after adding and subtracting a standard deviation to the average error for each noise level and algorithm.

For the case with $d = 5$ (Figure~\ref{fig:FPUT5_v5}) we can observe that in the differential formulation the two best performing algorithms are the CINDy (c) and the SR3 (c-$\ell_0$) frameworks, in terms of $\mathcal{E}_{R}$, $\mathcal{E}_{D}$ and $\mathcal{E}_{T}$, which highlight the importance of adding extra constraints, and showcase the potential improvement in performance that can be obtained from using SR3 (c-$\ell_0$) over SINDy. These two algorithms are closely followed by the CINDy framework for noise levels below $10^{-4}$. These three algorithms are also the most successful in correctly recovering the sparsity of the underlying dynamic, as can be seen in the results for $\mathcal{S}_E$. The sparsity realized by all the algorithms except SINDy, IPM and IPM (c) eventually cause them to have some missing basis functions, as we can see that for higher noise levels the value of $\mathcal{S}_M$ increases, which is expected. In terms of the integral formulation there is a large performance boost from using the CINDy and CINDy (c) frameworks as opposed to any of the other frameworks, as measured in all the metrics under consideration. In this case we do not observe SR3 (c-$\ell_0$) performing on par with its CG-based counterparts. Again, we observe that the algorithms that are successful in terms of $\mathcal{S}_E$ have a higher tendency to miss out on some of the basis functions as the level of noise increases, as seen in the graphs that depict $\mathcal{S}_M$. 

For the case with $d = 10$ (Figure~\ref{fig:FPUT10_v5}) we observe again that the best performing frameworks are the CINDy, CINDy (c), and FISTA variants for the differential formulation, while the best performing algorithms for the integral formulation are the CINDy and CINDy (c) frameworks. Note that in this case the FISTA algorithm outperforms all the variants except CINDy and CINDy (c) for the metrics being considered. We also note the change in regime for SINDy and SR3 (c-$\ell_0$) as the noise level increases in the differential formulation. For the former a noise level above $10^{-8}$ seems to indicate a significant increase in the number of extra basis functions, while for the latter the critical noise level seems to be around $10^{-6}$. In terms of performance with respect to $\mathcal{S}_{E}$ CINDy and CINDy (c) are significantly better than the other frameworks (except FISTA for the differential formulation), however for the highest noise levels they miss out on a significant number of basis functions, as measured by $\mathcal{S}_{M}$. For all noise levels except in one, SINDy tends to pick up all the available basis functions, as can be seen in the plot for $\mathcal{S}_{E}$.

In this set of experiments, for $d=5$ and $d = 10$ we can also clearly observe the higher robustness with respect to noise of the SR3 (c-$\ell_0$) algorithm over the SINDy framework, which was one of the main reasons why it was developed. Note that this was not clear in the experiments for the Kuramoto model. Note however that the CG-based frameworks are more robust than SR3 (c-$\ell_0$).

\subsubsection{Sample Efficiency} \label{section:FPUT:sample_efficiency}

\begin{figure*}[]
    \centering
    \vspace{-10pt}
    \hspace{\fill}
    \subfloat{{\includegraphics[width=8cm]{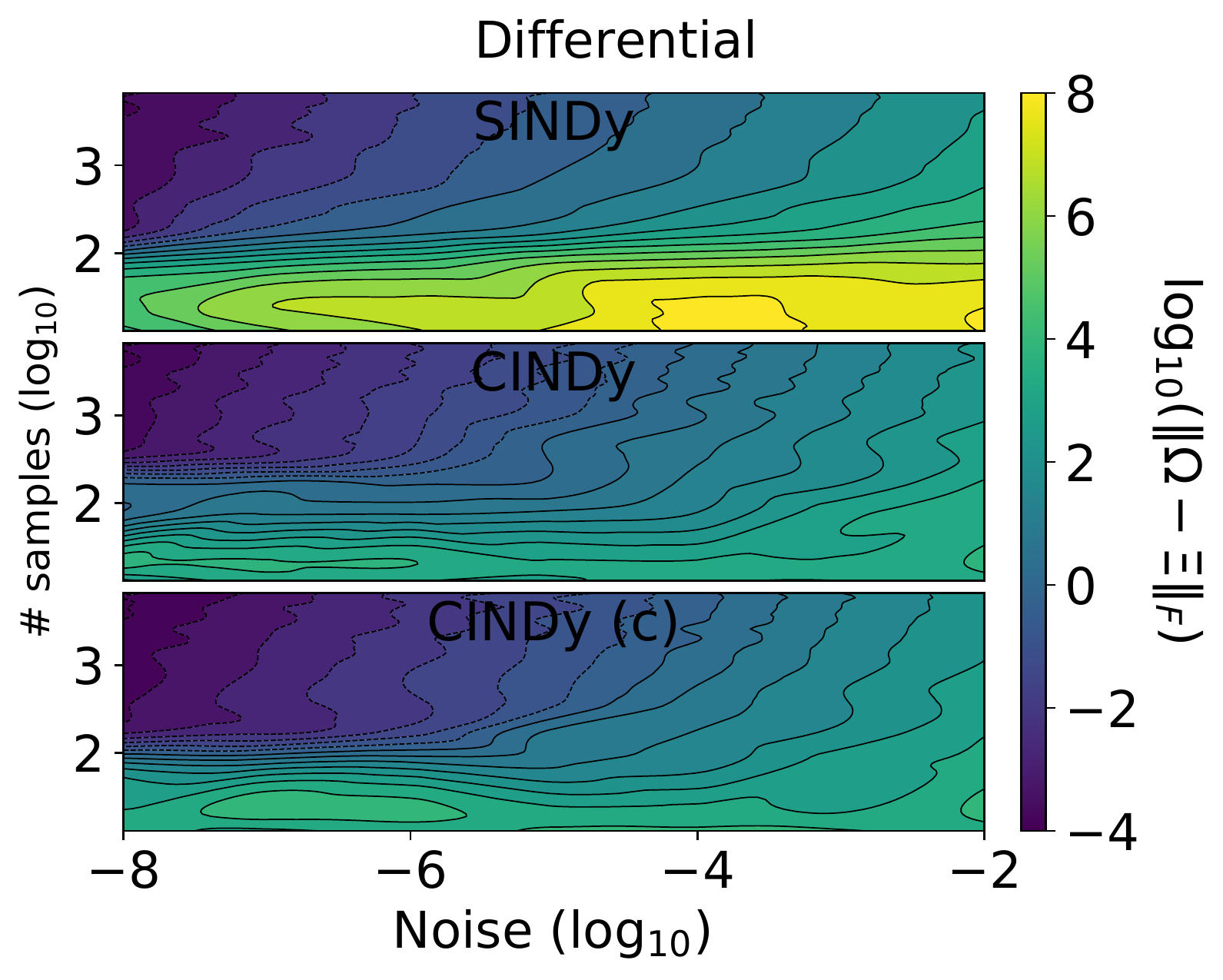} }\label{fig:FPUT:diff:dim5_sampleeff_acc}}%
    %\qquad
    \hspace{\fill}
    \subfloat{{\includegraphics[width=8cm]{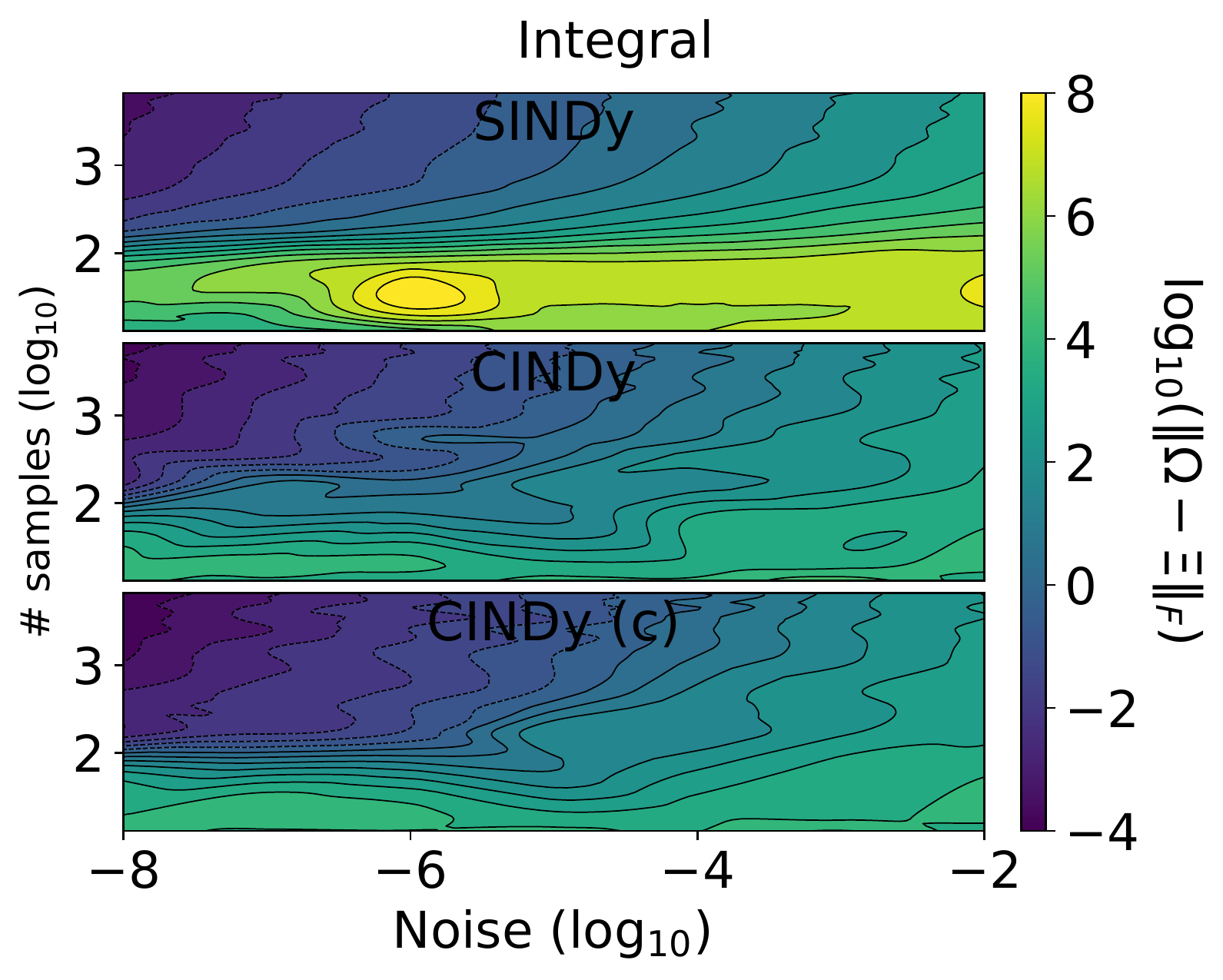} }\label{fig:FPUT:int:dim5_sampleeff_acc}}%
    \hspace{\fill}
    \caption{\textbf{Sample efficiency of the sparse recovery of the
        Fermi-Pasta-Ulam-Tsingou model: } Algorithm comparison for a Fermi-Pasta-Ulam-Tsingou model of
      dimension $d=5$ for both formulations.}%
    \label{fig:FPUT_sampleeff}%
  \end{figure*}

As in Section~\ref{section:kuramoto:sample_efficiency}, we present a heat map in Figure~\ref{fig:FPUT_sampleeff} that compares the accuracy, in terms of $\log(\mathcal{E}_R)$, as we vary the noise levels
(x-axis) and the numbers of training samples (y-axis) for the SINDy and CINDy algorithms. If one
takes a look at the differential formulation, one can see that in the
low training sample regime both CINDy and CINDy (c) perform slightly better
than SINDy at higher noise levels. The difference in performance is less pronounced than in Figure~\ref{fig:kuramoto_sampleeff}, however. We provide an extended
analysis over a broader range of sample sizes in
Appendix~\ref{appx:section:sample_efficiency} for completeness. In the images in the Appendix it is easier to discern the advantage of adding constraints to the system, as we can clearly see that the accuracy of CINDy (c) is higher than that of CINDy for a given noise level and number of samples.

\subsubsection{Simulation of learned trajectories} \label{section:FPUT:trajectories}

The results shown in Figure~\ref{fig:FPUT10_v5} for $\mathcal{E}_T$
(see third row of images) seem to indicate that all four formulations
will perform similarly when predicting trajectories, however, this
stands in contrast to what is shown in the first row of images, where
we can see that the dynamic learned by SINDy is far from the true
dynamic. This discrepancy is due to the fact that all the data is
generated in one regime of the dynamical phenomenon, and the metric
$\mathcal{E}_T$ is computed using noisy testing data from the same
regime. If we were to test the performance in inferring trajectories
with initial conditions that differed from those that had been seen in
the training-testing-validation data, the picture would be quite
different. For example one could test the inference power of the
different dynamics if the initial position of the oscillators is a
sinusoid with unit amplitude.

\begin{figure}[h!]
	\centering
	\includegraphics[width=14cm]{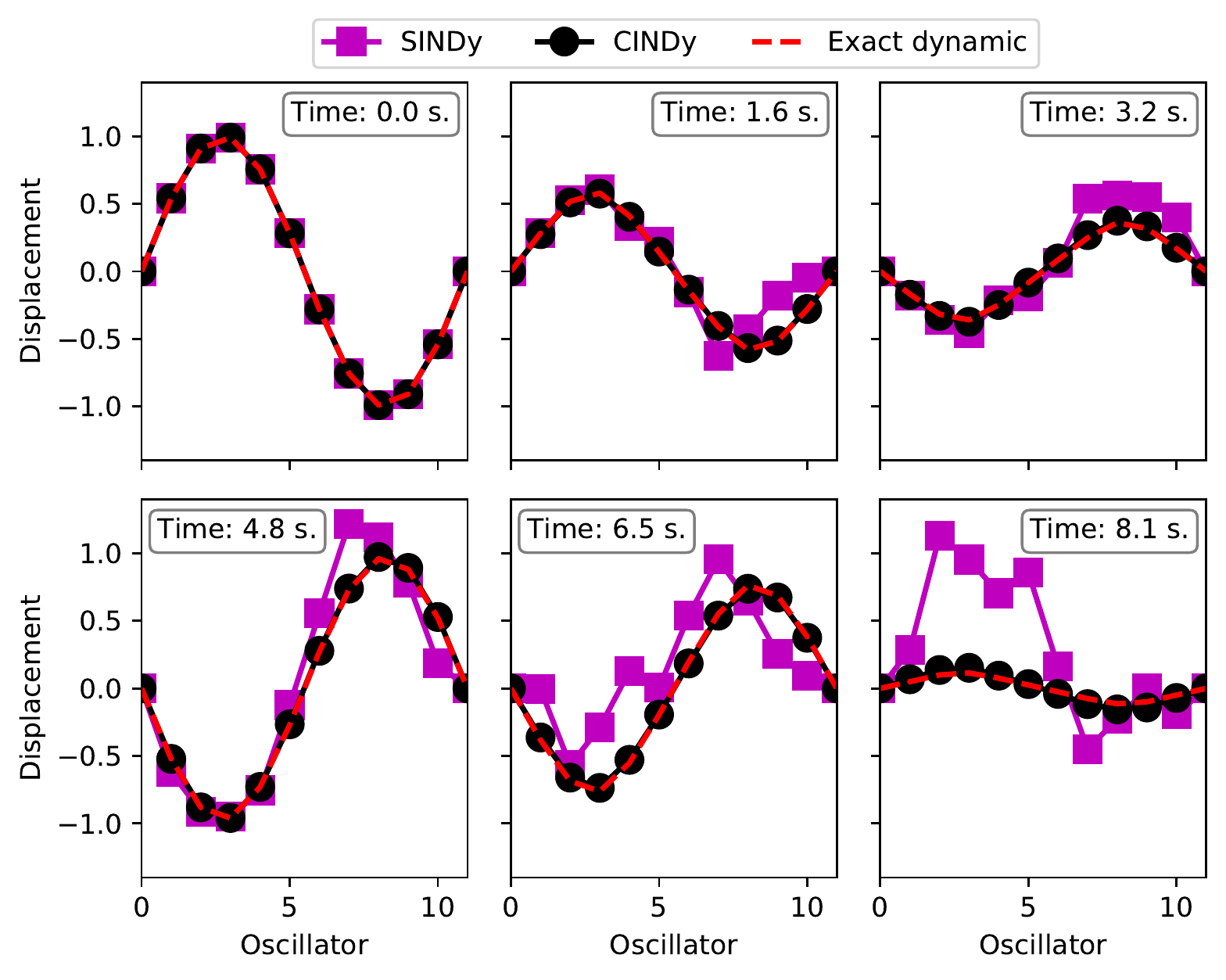} 
	\caption{\textbf{Trajectory comparison:} Fermi-Pasta-Ulam-Tsingou model of dimension $d = 10$.} \label{FPUT:trajectories_example}
\end{figure}

We can see the difference in accuracy between the different learned dynamics by simulating forward in time the dynamic learned by the CINDy algorithm and the SINDy algorithm, and comparing that to the evolution of the true dynamic. The results in Figure~\ref{FPUT:trajectories_example} show the difference in behaviour for different times for the dynamics learnt by the two algorithms in the integral formulation with a noise level of $10^{-4}$ for the example of dimensionality $d = 10$. In keeping with the physical nature of the problem, we present the ten dimensional phenomenon as a series of oscillators suffering a displacement on the vertical y-axis, in a similar fashion as was done in the original paper \cite{fermi1955studies}. Note that we have added to the images the two extremal particles that do not oscillate.

The two orders of magnitude in difference between the SINDy and the CINDy algorithms, in terms of $\mathcal{E}_{R}$, manifests itself clearly when we try to predict trajectories with out-of-sample initial positions that differ from those that have been used for learning. This suggests that the CINDy algorithm has better generalization properties under noise.

\subsection{Michaelis-Menten model}

The Michaelis-Menten model is used to describe enzyme reaction kinetics \cite{michaelis2007kinetik}. We focus on the following derivation \cite{briggs1925note}, in which an enzyme E combines with a substrate S to form an intermediate product ES with a reaction rate $k_{f}$. This reaction is reversible, in the sense that the intermediate product ES can decompose into E and S, with a reaction rate $k_{r}$. This intermediate product ES can also proceed to form a product P, and regenerate the free enzyme E. This can be expressed as
\begin{align*}
  \ch{S + E <>[ $k_f$ ][ $k_r$ ] E.S ->[ $k_{\text{cat}}$ ] E + P}  .
\end{align*}
If we assume that the rate for a given reaction depends proportionately on the concentration of the reactants, and we denote the concentration of E, S, ES and P as $x_{\text{E}}$, $x_{\text{S}}$, $x_{\text{ES}}$ and $x_{\text{P}}$, respectively, we can express the dynamics of the chemical reaction as:
\begin{align*}
    \dot{x}_{\text{E}} &= -k_f x_{\text{E}} x_{\text{S}} + k_r x_{\text{ES}} + k_{\text{cat}} x_{\text{ES}} \\
    \dot{x}_{\text{S}} &= -k_f x_{\text{E}} x_{\text{S}} + k_r x_{\text{ES}} \\
    \dot{x}_{\text{ES}} &= k_f x_{\text{E}} x_{\text{S}} - k_r x_{\text{ES}} - k_{\text{cat}} x_{\text{ES}} \\
    \dot{x}_{P} &=  k_{\text{cat}} x_{\text{ES}}.
\end{align*}
One of the interesting things about the Michaelis-Menten dynamic is that we can use some of the structural constraints described in Section~\ref{Section:Convervation}, that is, the exact dynamic satisfies
\begin{gather*}
    \dot{x}_{\text{S}} + \dot{x}_{\text{ES}} + \dot{x}_{\text{P}} = 0 \\
    \dot{x}_{\text{E}} + \dot{x}_{\text{ES}} = 0.  
\end{gather*}
The exact dynamic $\Xi$ can be expressed using a dictionary of monomials of degree up to two. To be more precise, the dictionary used is
\begin{align*}
   \mathcal{D} = \left\{ x_{\text{E}}^{a_{\text{E}}} x_{\text{S}}^{a_{\text{S}}} x_{\text{ES}}^{a_{\text{ES}}} x_{\text{P}}^{a_{\text{P}}} \mid  a_{\text{E}}, a_{\text{S}}, a_{\text{ES}}, a_{\text{P}} \in \llbracket 0, 2 \rrbracket,  0 \leq a_{\text{E}} + a_{\text{S}} + a_{\text{ES}} + a_{\text{P}} \leq 2\right\}.
\end{align*}
With this dictionary in mind, and denoting the coefficient vector associated with the chemical E as $\xi_{\text{E}}$, and likewise for the other chemicals, we can impose the following additional series of constraints for all $i$:

\begin{gather*}
 \xi_{\text{E}} + \xi_{\text{S}} + \xi_{\text{P}} = 0 \\
 \xi_{\text{E}} + \xi_{\text{ES}} = 0.
\end{gather*}

\begin{figure*}[]
    \centering
    \hspace{\fill}
    \subfloat{{\includegraphics[width=15.6cm]{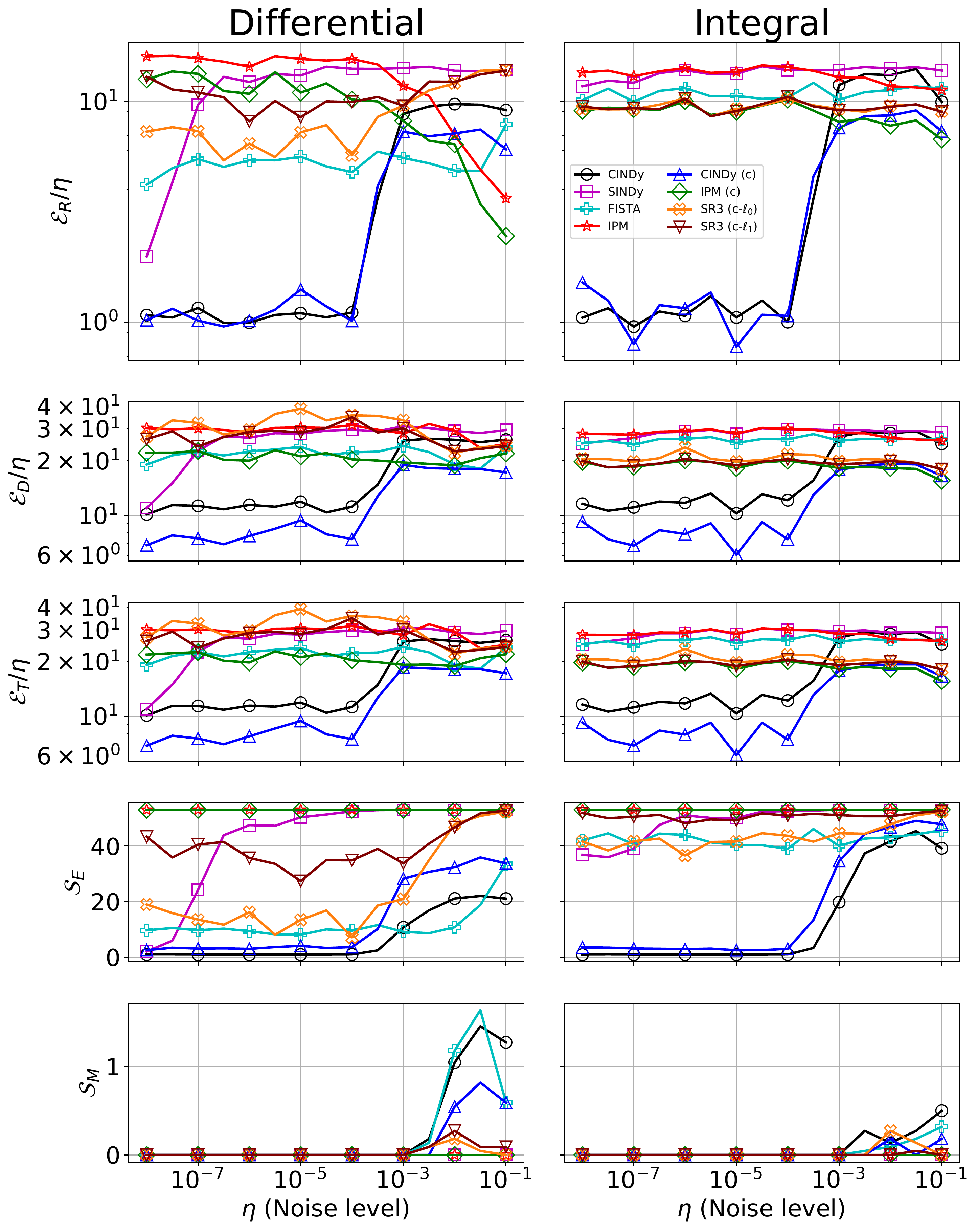} }\label{fig:MMeasy}}%
    \hspace{\fill}
    \caption{\textbf{Sparse recovery of the Michaelis-Menten model: } Algorithm comparison for a Michaelis-Menten model of dimension $d=4$, with a differential formulation shown on the left column, and with an integral formulation on the right column. The first, second, third, fourth and firth rows of images indicate a comparison of $\mathcal{E}_R/\eta$, $\mathcal{E}_D/\eta$, $\mathcal{E}_T/\eta$, $\mathcal{S}_E$, and $\mathcal{S}_M$, as we vary the noise level $\eta$, respectively.}%
    \label{fig:MM4}%
\end{figure*}

The images in Figure~\ref{fig:MM4} show the recovery results for $k_f = 0.01$ $k_r = 1$, $k_{\text{cat}}= 1$ and $d = 4$. The derivatives used were computed using a polynomial interpolation of degree 8. We used polynomial interpolation of degree $8$ to compute the integrals for all the noise levels. A total of $6000$ points were used to infer the dynamic, spread over $150$ experiments for a maximum time of $0.01$ seconds, for both cases. The initial state of the system for the $j$-th experiment was selected randomly as $\vx^j(t_0) \sim \mathcal{U}[0, 1]^4$. 

%\begin{gather*}
%\left\lvert \left( \xi_{\text{E}} + \xi_{\text{S}} + \xi_{\text{P}}\right)^T  \bm{\psi}(\vy(t_i))   \right\rvert \leq \varepsilon \\
%\left\lvert \left( \xi_{\text{E}} + \xi_{\text{ES}} \right)^T  \bm{\psi}(\vy(t_i))   \right\rvert \leq \varepsilon'.
%\end{gather*}

%\todo[inline]{AC: Right now the constraints I've imposed in the code are $\xi_{\text{E}} + \xi_{\text{S}} + \xi_{\text{P}} = 0$ and $\xi_{\text{E}} + \xi_{\text{ES}} = 0$.}

%Where the $\varepsilon>0$ can be selected based on the performance on testing data. The images in Figure~\ref{fig:MM4} show the recovery results for $k_f = 0.01$ $k_r = 1$, $k_{\text{cat}}= 1$ and $d = 5$. The number of samples generated, the number of experiments performed and all the other experiment details are the same as in Section~\ref{section:FPUT}. The derivatives used where computed using a polynomial interpolation of degree 8. We used polynomial interpolation of degree $8$ to compute the integrals for all the noise levels.

In this experiment we observe that for both the integral and differential formulation there is a significant performance advantage from using the CINDy and CINDy (c) frameworks for noise levels below $10^{-3}$, in terms of all the metrics under consideration. As in the previous examples, there is also a huge difference in the number of extra basis functions picked up, as measured by $\mathcal{S}_E$, with CINDy and CINDy (c) producing in general the sparsest solutions, followed by FISTA and SR3 (c-$\ell_0$). The IPM, IPM (c) and SINDy frameworks tend to pick up all the available basis functions for noise levels above $10^{-6}$. In terms of correct basis functions that have not been picked up, FISTA, CINDy and CINDy (c) only miss out on average on less than one of the basis functions for the highest noise levels in the differential formulation. Similar comments can be made regarding the integral formulation, where we observe that the aforementioned three algorithms only miss out on some basis functions for noise levels above $10^{-2}$.

\section*{Conclusion}
\label{sec:conclusion}
We have presented a CG-based optimization algorithm, namely the Blended Conditional Gradients algorithm, that can be used to solve a convex sparse recovery problem formulation, resulting in the CINDy framework. In comparison with other existing frameworks, CINDy shows more accurate recovery under noise, while also outperforming other algorithms in terms of sparsity. Moreover, as the underlying optimization algorithm relies on a linear optimization algorithm, we can easily encode linear inequality constraints into the learning formulation to achieve dynamics that are consistent with the physical phenomenon. Existing algorithms, on the other hand, were only able to deal with linear equality constraints.

\section*{Acknowledgments}
\label{sec:acknowledgments}
This research was partially funded by Deutsche Forschungsgemeinschaft (DFG) through the DFG Cluster of Excellence MATH+ and Project A05 in CRC TRR 154.

\bibliography{references}

\bibliographystyle{icml2021}

\appendix

\section{Preliminaries} 
Given a differentiable function $f(\vx): \rr^d \rightarrow \rr$ we say the function $f(\vx)$ is:
\begin{definition}[$L$-smooth] A function is $L$-smooth if for any $\vx,\vy \in \rr^d$ we have:
\begin{align*}
    f(\vx) \leq f(\vy) + \innp{\nabla f(\vy), \vy - \vx} + \frac{L}{2}\norm{\vx - \vy}^2.
\end{align*}
This is equivalent to the gradient of $f(\vx)$ being $L$-Lipschitz. If the function is twice-differentiable this is equivalent to:
\begin{align}
    0 < L = \max\limits_{\vx,\vy \in \rr^d} \frac{(\vx - \vy)^T \nabla^2 f(\vx)(\vx - \vy)}{\norm{\vx - \vy}^2}
\end{align}
\end{definition}
\begin{definition}[Convex] A function is convex  if for any $\vx,\vy \in \rr^d$ we have:
\begin{align*}
    f(\vx) \geq f(\vy) + \innp{\nabla f(\vy), \vy - \vx} 
\end{align*}
If the function is twice-differentiable this is equivalent to:
\begin{align}
    0 = \min\limits_{\vx,\vy \in \rr^d} \frac{(\vx - \vy)^T \nabla^2 f(\vx)(\vx - \vy)}{\norm{\vx - \vy}^2}
\end{align}
\end{definition}
\begin{definition}[$\mu$-strongly convex] A function is $\mu$-strongly convex if for any $\vx,\vy \in \rr^d$ we have:
\begin{align*}
    f(\vx) \geq f(\vy) + \innp{\nabla f(\vy), \vy - \vx} + \frac{\mu}{2}\norm{\vx - \vy}^2.
\end{align*}
If the function is twice-differentiable this is equivalent to:
\begin{align}
    0 < \mu = \min\limits_{\vx,\vy \in \rr^d} \frac{(\vx - \vy)^T \nabla^2 f(\vx)(\vx - \vy)}{\norm{\vx - \vy}^2}
\end{align}
\end{definition}
Given a compact convex set $\cx \subset \rr^d$ we define the constrained optimization problem:
\begin{align}
    \min\limits_{\vx \in \mathcal{X}} f(\vx). \label{eq:generic_problem}
\end{align}

\section{Accelerated Projected Gradient Descent} \label{NAGD}
In this section we will focus in the case where $f(\vx)$ is $L$-smooth and convex (or potentially $\mu$-strongly convex). One of the key characteristics of convex problems is that any local minima to the problem in Equation~\eqref{eq:generic_problem} is a global minima. Moreover, there exist efficient algorithms for computing the minima of these problems. In order to tackle the problem shown in Equation~\eqref{eq:generic_problem}, one can, for example, use either a projection-based or a projection-free methods, depending on how computationally difficult it is to compute projections onto $\cx$. We denote the \emph{Euclidean projection} of $\vx$ onto $\cx$ as $\Pi_{\cx}\left( \vx\right): \rr^n \rightarrow \cx$, which is defined as:
\begin{align*}
    \Pi_{\cx}\left( \vx\right) \defeq \argmin\limits_{\vy \in \cx} \frac{1}{2} \norm{\vx - \vy}^2.
\end{align*}
In general, computing these projections is non-trivial. However, for a series of structured feasible regions there are closed-form expressions for these projections, which can be computed efficiently (see Table~\ref{Table:complexity_projections}):

\begin{table*}[h!]
 \begin{center} 
 \begin{tabular}{lll} 
  \toprule
  \textbf{Feasible region $\mathcal{X}$}&\textbf{Mathematical expression}&\textbf{Projection}\\
  \midrule
   Unit probability simplex&$\{\vx\in\mathbb{R}^d\mid \mathbf{1}_d^\top\vx=1,\vx\geq\mathbf{0}\}$ &$\mathcal{O}(d)$\\
  $\ell_p$-ball, $p\in\{1,2,+\infty\}$&$\{\vx\in\mathbb{R}^d\mid\|\vx\|_p\leq1\}$&$\mathcal{O}(d)$\\
  Nuclear norm-ball&$\{\mathbf{X}\in\mathbb{R}^{m\times n}\mid\|\mathbf{X}\|_
  {\operatorname{nuc}}\leq1\}$&$\mathcal{O}(mn\min\{m,n\})$\\
  Matroid polytope&$\{\vx\in\mathbb{R}^d\mid\forall S\in\mathcal{P}(E),\mathbf{1}_S^\top\vx\leq r(S),\vx\geq\mathbf{0}\}$&$\mathcal{O}(\operatorname{poly}(d))$\\
  \bottomrule
 \end{tabular}
 \end{center}
 \caption{Complexities of projections onto several feasible regions.}
\label{Table:complexity_projections}
\end{table*}

Fortunately, projections onto the probability simplex can be computed very efficiently, which makes \emph{accelerated projected descent algorithms} an attractive alternative when solving constrained convex problems over the probability simplex (as is the case in Line~\ref{alg:quadratic_problemFCFW} of Algorithm~\ref{algo:FCFW}). These algorithms are termed \emph{accelerated} because they are able to improve upon the convergence guarantees offered by the standard projected gradient descent algorithm, both in the convex case (see Algorithm~\ref{algo:NAGD}) and in the strongly convex case (see Algorithm~\ref{algo:NAGD_strcvx}). If we measure optimality by the number of iterations $k$ needed for the algorithms to achieve an $\epsilon$-optimal accuracy (which means that $f(\vx_k) - \min_{\vx \in \cx} f(\vx) \leq \epsilon$), then in the smooth convex case the accelerated projected gradient descent algorithm is able to reach an $\epsilon$-optimal solution in $\mathcal{O}(1/\sqrt{\epsilon})$ iterations, as opposed to the $\mathcal{O}(1/\epsilon)$ iterations needed with standard projected gradient descent. In the smooth and strongly convex case the accelerated projected gradient descent achieves an $\epsilon$-optimal solution in $\mathcal{O}(\sqrt{\mu/L}\log 1/\epsilon)$ iterations, as opposed to the $\mathcal{O}(\mu/L\log 1/\epsilon)$ iterations needed for the standard projected gradient descent algorithm.

\begin{algorithm}[]
\SetKwInOut{Input}{Input}\SetKwInOut{Output}{Output}
\Input{Objective function $f(\vx)$, feasible region $\cx$, initial point $\vx_0 \in \cx$.}
\Output{Point $\vx_{K+1} \in \cx$.}
\hrulealg
$\vy_0 \leftarrow \vx_0$\;
$L \leftarrow \max_{\vx, \vy\in \cx} (\vx - \vy)^T\nabla^2 f(\vx) (\vx - \vy)/\norm{\vx - \vy}^2$\;
$\gamma_0 \leftarrow 0 $\;
\For{$k = 1$ to $K$}{
$\vx_{k+1} \leftarrow \Pi_{\cx} \left( \vy_k -\frac{1}{L} \nabla f(\vy_k) \right)$ \;
$\gamma_{k+1} \leftarrow \frac{1 + \sqrt{1 + 4 \gamma_k^2}}{2}$\;
$\vy_{k+1} \leftarrow \vx_{k+1} + \frac{\gamma_{k} - 1}{\gamma_{k+1}} \left( \vx_{k+1} - \vx_{k}\right)$\;
}
\caption{Accelerated gradient descent for smooth convex problems.} \label{algo:NAGD}
\end{algorithm}

\begin{algorithm}[]
\SetKwInOut{Input}{Input}\SetKwInOut{Output}{Output}
\Input{Objective function $f(\vx)$, feasible region $\cx$, initial point $\vx_0 \in \cx$.}
\Output{Point $\vx_{K+1} \in \cx$.}
\hrulealg
$\vy_0 \leftarrow \vx_0$\;
$\mu \leftarrow \min_{\vx, \vy\in \cx} (\vx - \vy)^T\nabla^2 f(\vx) (\vx - \vy)/\norm{\vx - \vy}^2$\;
$L \leftarrow \max_{\vx, \vy\in \cx} (\vx - \vy)^T\nabla^2 f(\vx) (\vx - \vy)/\norm{\vx - \vy}^2$\;
\For{$k = 1$ to $K$}{
$\vx_{k+1} \leftarrow \Pi_{\cx} \left( \vy_k -\frac{1}{L} \nabla f(\vy_k) \right)$ \;
$\vy_{k+1} \leftarrow \vx_{k+1} + \frac{1 - \sqrt{\mu/L}}{1 + \sqrt{\mu/L}} \left( \vx_{k+1} - \vx_{k}\right)$\;
}
\caption{Accelerated gradient descent for smooth strongly-convex problems.} \label{algo:NAGD_strcvx}
\end{algorithm}

\section{On computing integrals and derivatives from noisy data} \label{appx:derivs_and_ints}

One of the key requirements for the success of any sparse recovery algorithm is the accurate estimation of integrals and derivatives. We are typically given a dictionary of basis functions $\mathcal{D} = \left\{\psi_i \mid i \in \llbracket 1, n \rrbracket \right\}$ that can be used to represent our dynamic, encoded by $\dot{\vx}(t) = \Xi^T \bm{\psi}(\vx(t))$, where $\Xi \in \rr^{n \times d}$ and $\bm{\psi}(\vx(t)) = \left[ \psi_1(\vx(t)), \cdots, \psi_n(\vx(t)) \right]^T \in \rr^{n}$ along with some data points. Ideally, we would like to observe a series of noise-free data points from the physical system $\left\{ \vx(t_i)\right\}_{i=1}^m$. We would also like to observe $\left\{ \dot{\vx}(t_i)\right\}_{i=1}^m$, if we follow the differential approach, so that we can construct the matrix $\dot{X}$ in the left-hand side of Equation~\eqref{eq:appx:l1_minimization}. Or we would like to observe $\{\int_{t_1}^{t_{j+1}} \psi_i(\vx(\tau)) d\tau\}_{j=1}^{m-1}$ for all $i\in\llbracket 1, n \rrbracket$, if we follow the integral approach, so that we can construct the $\Gamma (X)$ matrix on the right-hand side of Equation~\eqref{eq:appx:l1_minimization}. 

\begin{mdframed}[linewidth=0.5mm]
\begin{minipage}{0.5\textwidth}
\begin{nscenter}
\textbf{Noise-free LASSO Differential approach}
\end{nscenter}
\begin{gather*}
\argmin\limits_{\substack{ \norm{\Omega}_{1,1} \leq \alpha  \\ \Omega \in \rr^{n \times d}}}  \norm{\dot{X} - \Omega^T \Psi(X)}^2_F
\end{gather*}
\end{minipage}
\vline
\begin{minipage}{0.5\textwidth}
\begin{nscenter}
\textbf{Noise-free LASSO Integral approach}
\end{nscenter}
\begin{gather}
\argmin\limits_{\substack{ \norm{\Omega}_{1,1} \leq \alpha  \\ \Omega \in \rr^{n \times d}}}  \norm{\delta X - \Omega^T \Gamma(X)}^2_F \label{eq:appx:l1_minimization}
\end{gather}
\end{minipage}
\end{mdframed}

However, in general we do not observe either of these quantities, and have to estimate $\dot{X}$ and $\Gamma (X)$ from $\left\{ \vx(t_i)\right\}_{i=1}^m$. The estimation of these matrices is made even more difficult if we are only able to observe noise-corrupted data points $\left\{ \vy(t_i)\right\}_{i=1}^m$ such that $\vy(t_i) = \vx(t_i) + \nu(t_i)$, where $\left\{ \nu(t_i)\right\}_{i=1}^m$ is a set of i.i.d. noise vectors, which is typically the case. In this case, we can only work with a noisy matrix $Y$, with which we either have to estimate $\dot{Y} \in \rr^{d \times m}$ or $\Gamma(Y) \in \rr^{n \times m - 1}$. The error in the formulation now comes from both the noisy data, and from the estimation of the derivatives/integrals from noisy data. The question then becomes, how do $\norm{\dot{Y} - \dot{X}}_F$ and $\norm{\Gamma(Y) - \Gamma(X) }_F$ evolve as we increase the noise? Given data contaminated with noise, is it easier to get an accurate estimate of $\dot{Y}$ than of $\Gamma(Y)$? If so, we might favor one approach over the other. 

\begin{figure*}[h!]
    \centering
    \vspace{-10pt}
    \hspace{\fill}
    \subfloat{{\includegraphics[width=8cm]{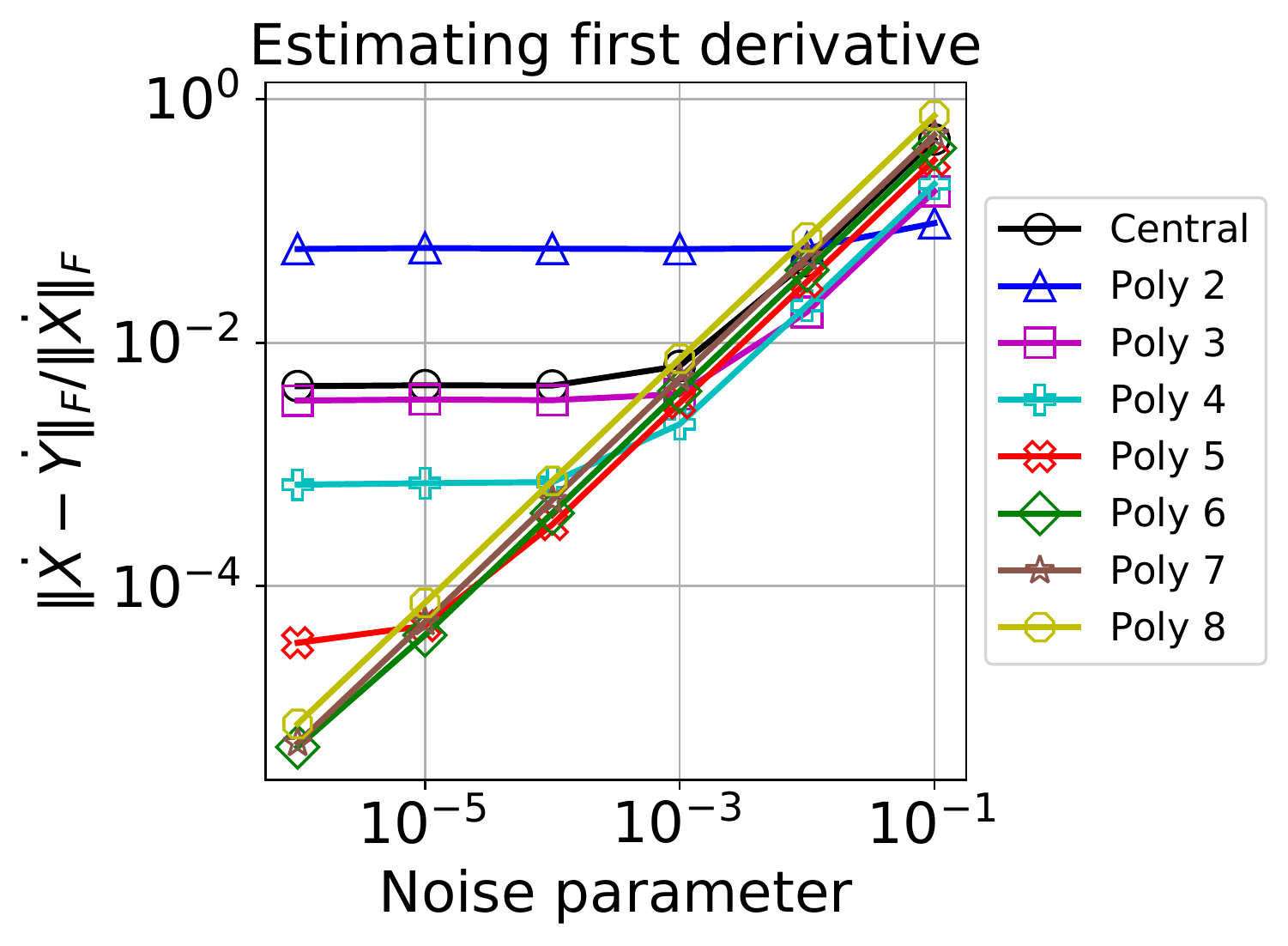} }\label{fig:first_derivative}}%
    %\qquad
    \hspace{\fill}
    \subfloat{{\includegraphics[width=8cm]{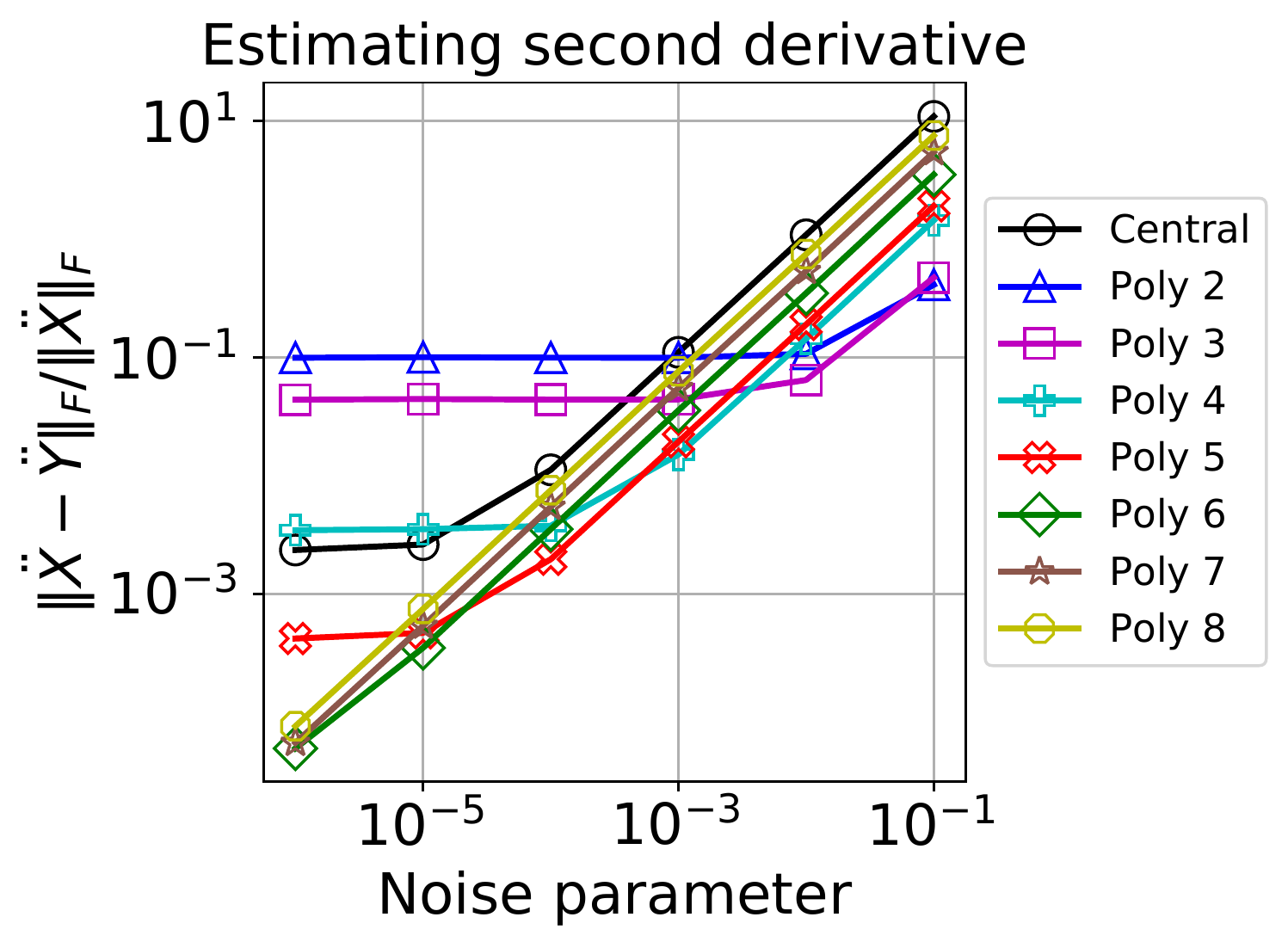} }\label{fig:second_derivative}}%
    \hspace{\fill}
    \caption{\textbf{Fermi-Pasta-Ulam-Tsingou:} Comparison of estimates of first and second derivatives with respect to time with exact first and second order derivatives.  }%
    \label{fig:comparison_derivative}%
\end{figure*}

For example if we let the Fermi-Pasta-Ulam-Tsingou system with $d = 5$ evolve from a random initial state, and we sample the system $T/c$ times at regularly spaced intervals, and we repeat this experiment $c = 60$ times, we generate $T = 2400$ data points $\vx^j(t_i)$ with $i\in \llbracket 1, T/c\rrbracket$ and $j\in \llbracket 1, c\rrbracket$, with $t_{T/c} = 3$ seconds. We can proceed to corrupt these data points with Gaussian noise as in Section~\ref{Section:Numerical_experiments}, generating data points with $\vy^j(t_i) = \vx^j(t_i) + \alpha \mathcal{N}\left(0, \Sigma \right)$. In this case we need to estimate $\ddot{X}$ if we want to formulate the problem from a differential perspective, or $\dot{X}$ if we want to formulate the problem from an integral perspective. Using the same dictionary of basis functions as the ones used in Section~\ref{section:FPUT}, we use the first-order central difference rule as well as differentiation of local polynomial interpolations to estimate the first derivative of $X$ with respect to time, using $Y$. The results are shown on the image of the left in Figure~\ref{fig:comparison_derivative}. We also use the second-order central difference formula, as well as differentiation of local polynomial interpolations to estimate the second derivative of $X$ with respect to time, using $Y$. The results can be seen on the image to the right in Figure~\ref{fig:comparison_derivative}.

\begin{figure*}[h!]
    \centering
    \subfloat{{\includegraphics[width=8cm]{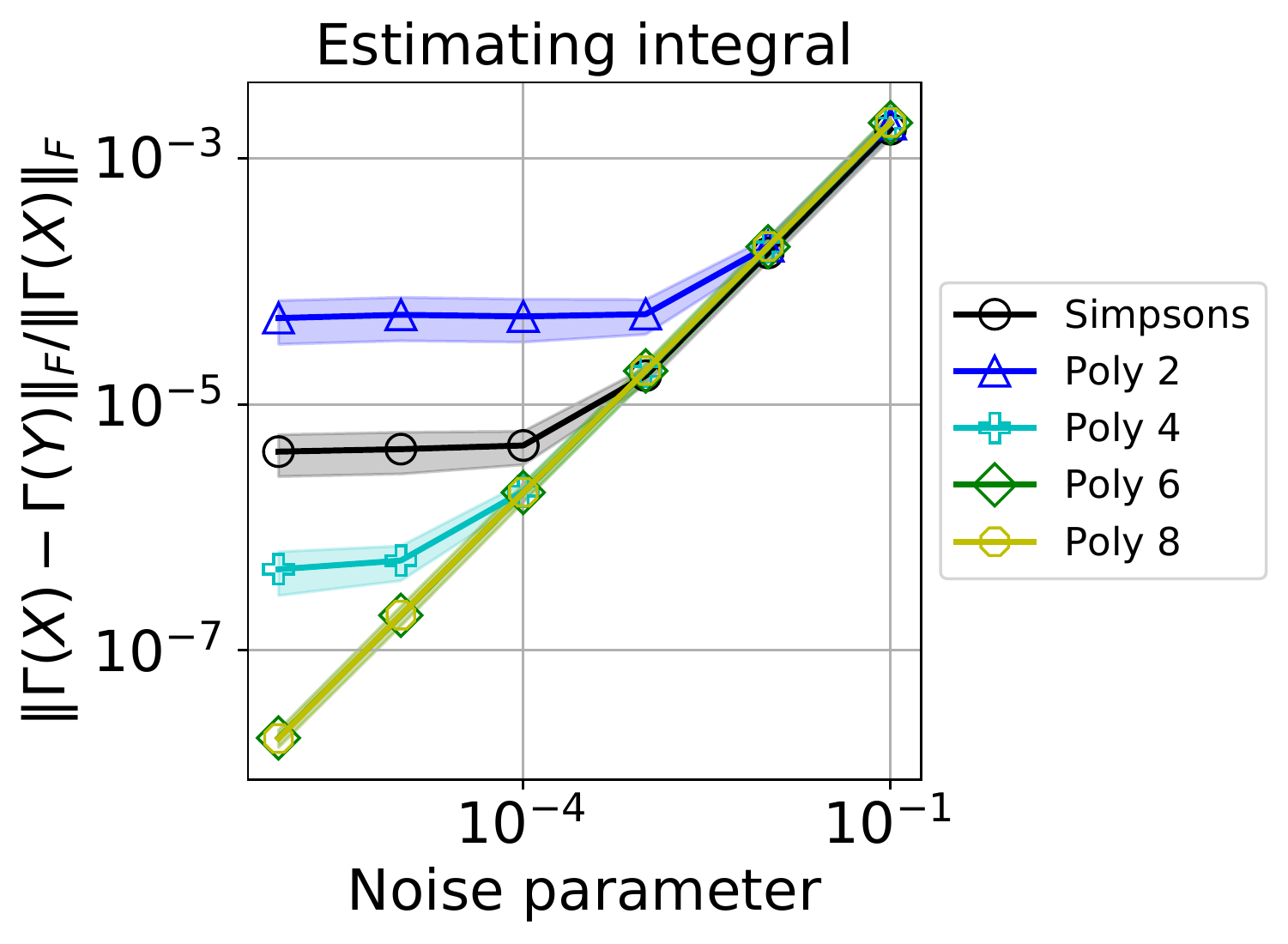} }}%
    \caption{\textbf{Fermi-Pasta-Ulam-Tsingou:} Comparison of estimates of first and second derivatives with respect to time with exact first and second order derivatives.  }%
    \label{fig:comparison_integral}%
\end{figure*}

We  also use Simpsons quadrature rules as well as integration of local polynomial interpolations to estimate the integral $\Gamma(\ddot{X})$ using $\dot{Y}$, which is our estimate of $\dot{X}$, the errors when computing these integrals with the aforementioned methods are shown in Figure~\ref{fig:comparison_integral}.

Regarding the errors shown between approximating $\dot{X}$ and $\ddot{X}$, for a given method we expect the errors estimating $\ddot{X}$ to be higher than the ones in estimating $\dot{X}$, which is what we observe in the experiments. There is no direct way to make a fair comparison between an approximation to the matrices $\dot{X}$ and $\ddot{X}$ and an approximation to the matrix $\Gamma (\dot{X})$, given their different nature (and even size), however if we use the metric $\norm{\dot{X} - \dot{Y}}_F / \norm{\dot{X}}_F$, $\norm{\ddot{X} - \ddot{Y}}_F / \norm{\ddot{X}}_F$ and $\norm{\Gamma (\dot{X}) - \Gamma(\dot{Y})}_F / \norm{\Gamma(\dot{X})}_F$ to compare the two approximations, the data seems to suggest that it is easier to estimate the matrix $\Gamma (\dot{X})$ than it is to estimate $\ddot{X}$, at least in the current experiment with the Fermi-Pasta-Ulam-Tsingou model. %Note for example that with Simpsons quadrature rule to estimate $\Gamma (X)$ we can get a level of accuracy that far surpasses the one achieved with the differentiation of an $8$-th order local polynomial to estimate $\ddot{X}$, even though, as we said, the direct comparison of the $\norm{\ddot{X} - \ddot{Y}}_F / \norm{\ddot{X}}_F$ and $\norm{\Gamma (\dot{X}) - \Gamma(\dot{Y})}_F / \norm{\Gamma(X)}_F$ metrics is problematic.

\newpage

\section{Additional figures} \label{appx:section:additional_figures}

\begin{figure*}[h!]
    \centering
    \vspace{-10pt}
    \hspace{\fill}
    \subfloat[]{{\includegraphics[width=3.95cm]{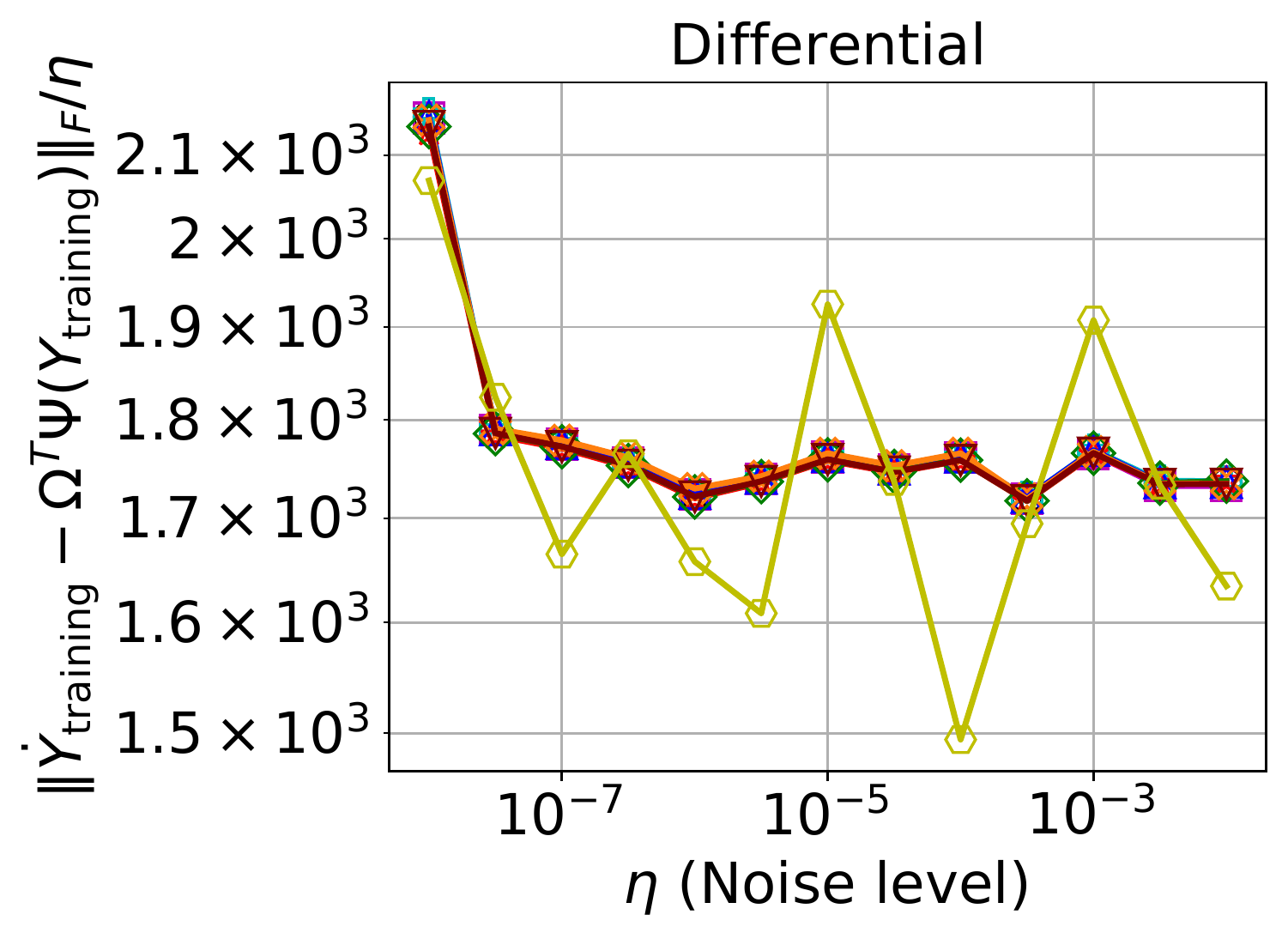} }\label{fig:appx:Kura5:diff:training}}%
    %\qquad
    \hspace{\fill}
    \subfloat[]{{\includegraphics[width=3.85cm]{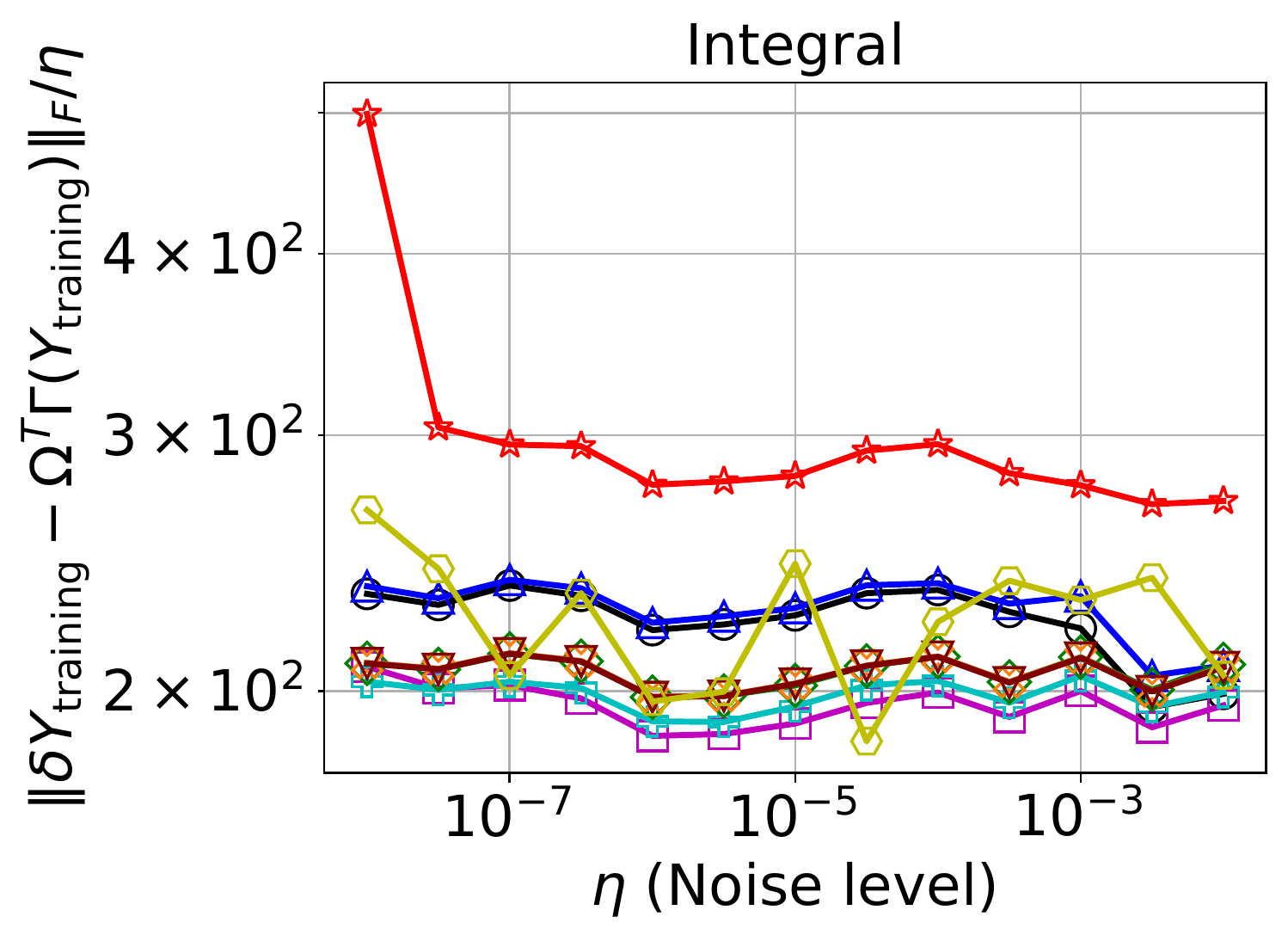} }\label{fig:appx:Kura5:int:training}}%
    \hspace{\fill}
    \subfloat[Iteration]{{\includegraphics[width=3.9cm]{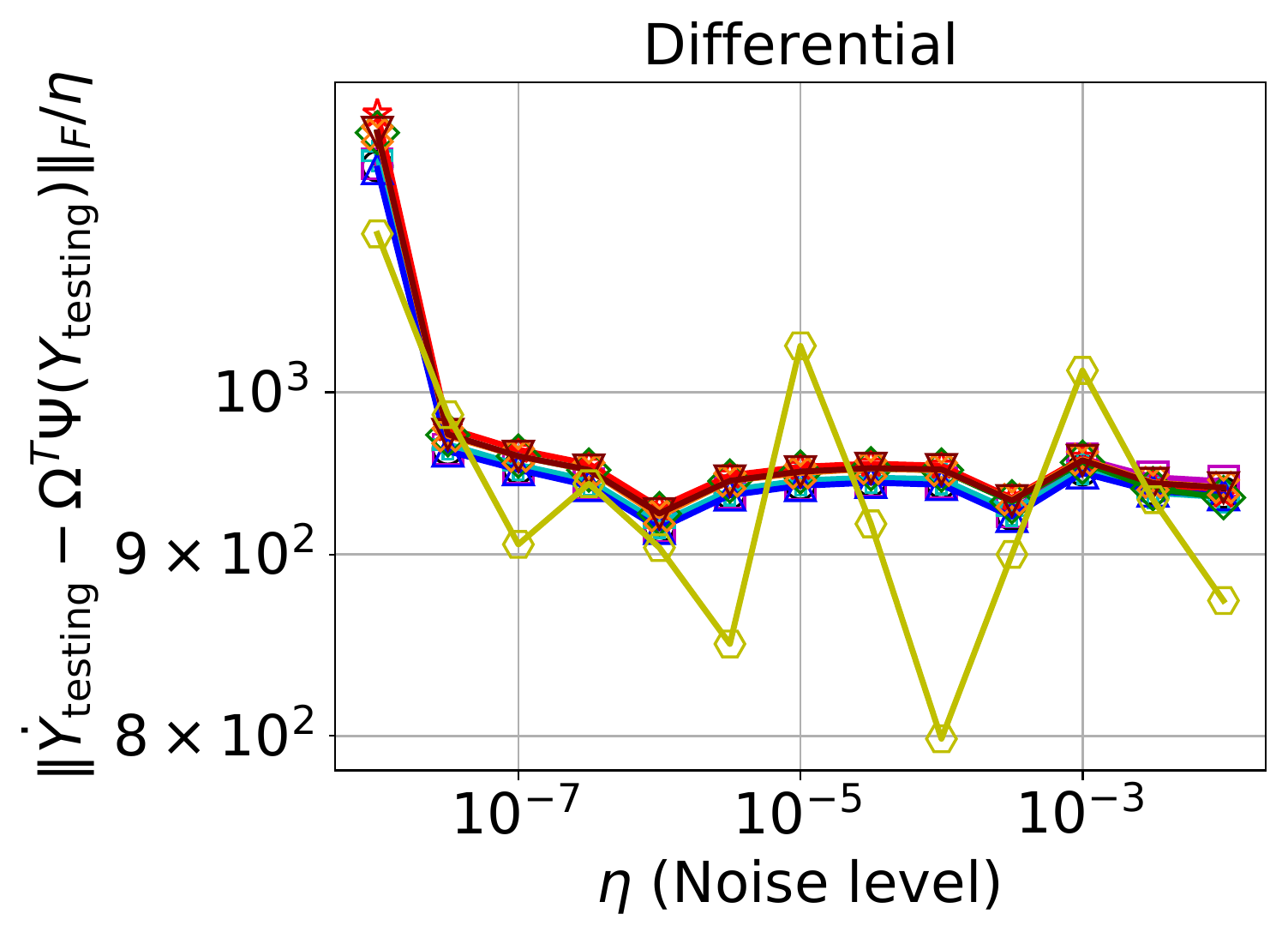} }\label{fig:appx:Kura5:diff:val}}%
    %\qquad
    \hspace{\fill}
    \subfloat[]{{\includegraphics[width=3.85cm]{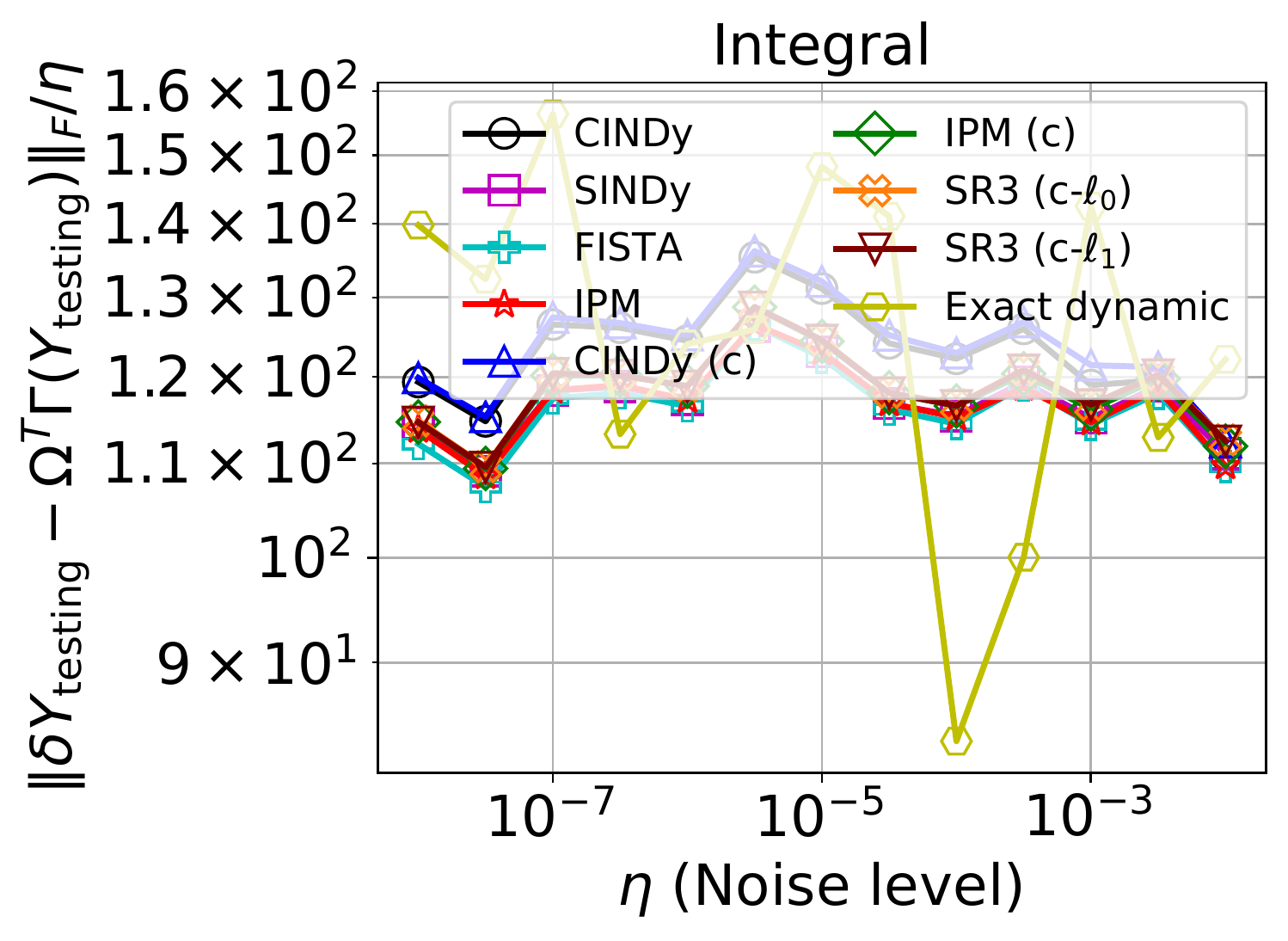} }\label{fig:appx:Kura5:int:val}}%
    \hspace*{\fill}
    \caption{Kuramoto: Evaluation of $\norm{\ddot{Y}_{\mathrm{training}} -  \Omega^T \Psi(Y_{\mathrm{training}})}_F$ \protect\subref{fig:appx:Kura5:diff:training} for the differential formulation, and $\norm{ \delta \dot{Y}_{\mathrm{training}} -  \Omega^T \Gamma(Y_{\mathrm{training}})}_F$ \protect\subref{fig:appx:Kura5:int:training} for the integral formulation for the experiments with $d = 5$, and evaluation of $\norm{\ddot{Y}_{\mathrm{validation}} -  \Omega^T \Psi(Y_{\mathrm{validation}})}_F$ \protect\subref{fig:appx:Kura5:diff:val} for the differential formulation, and $\norm{ \delta \dot{Y}_{\mathrm{validation}} -  \Omega^T \Gamma(Y_{\mathrm{validation}})}_F$ \protect\subref{fig:appx:Kura5:int:val} for the integral formulation for the experiments with $d = 5$.}%
    \label{fig:appx:kura5}%
\end{figure*}

\begin{figure*}[h!]
    \centering
    \vspace{-10pt}
    \hspace{\fill}
    \subfloat[]{{\includegraphics[width=3.95cm]{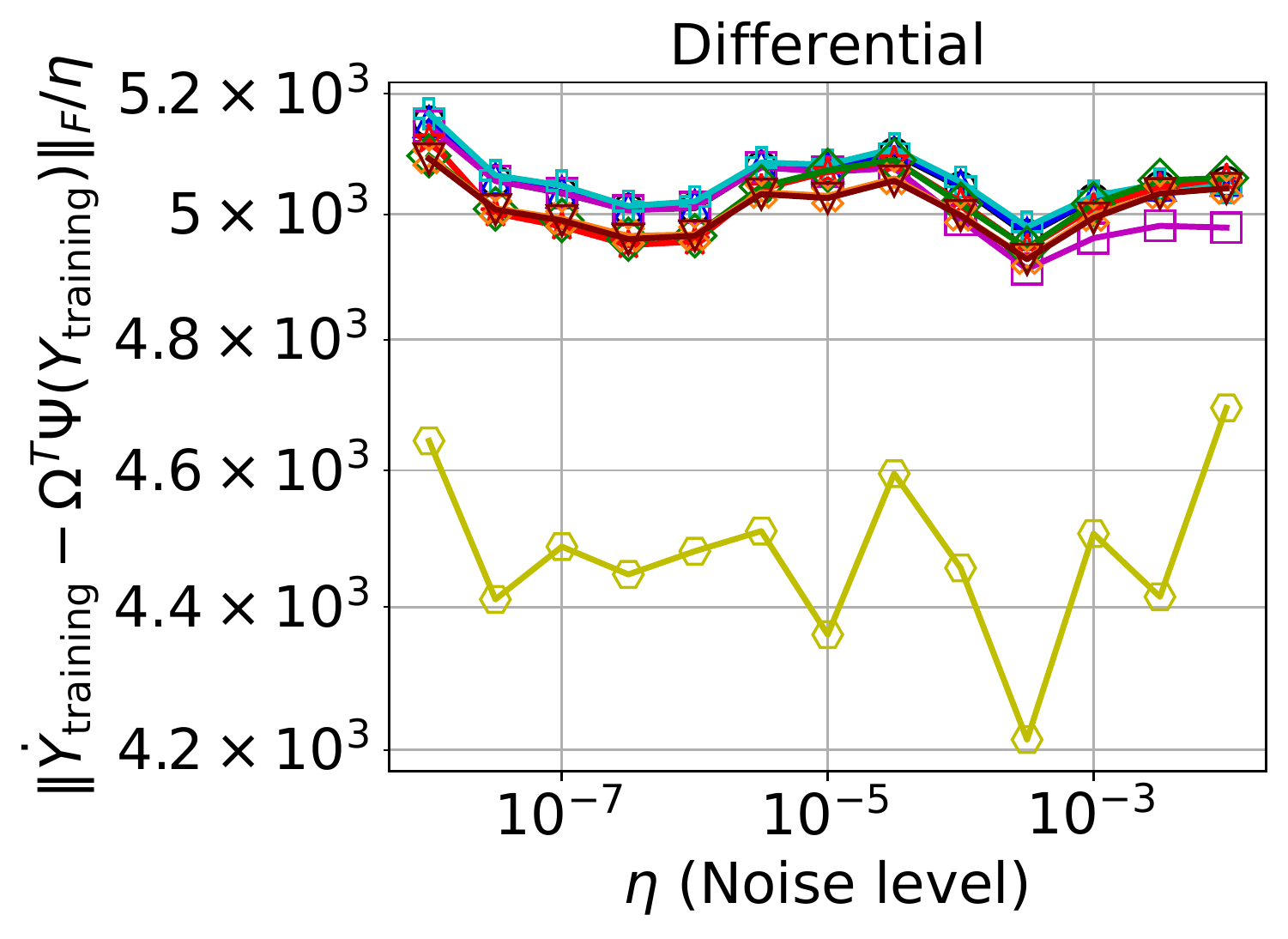} }\label{fig:appx:Kura10:diff:training}}%
    %\qquad
    \hspace{\fill}
    \subfloat[]{{\includegraphics[width=3.85cm]{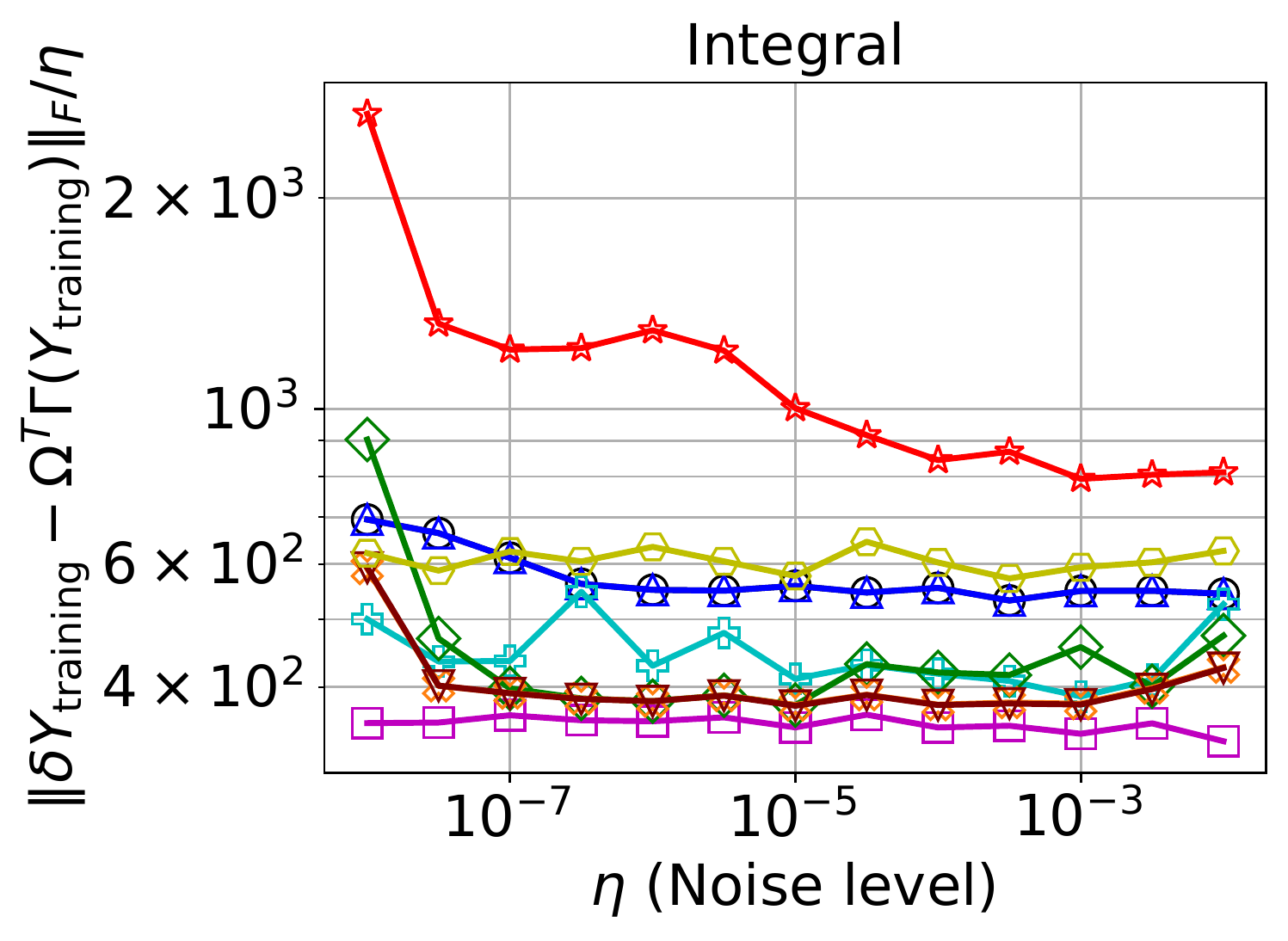} }\label{fig:appx:Kura10:int:training}}%
    \hspace{\fill}
    \subfloat[]{{\includegraphics[width=3.9cm]{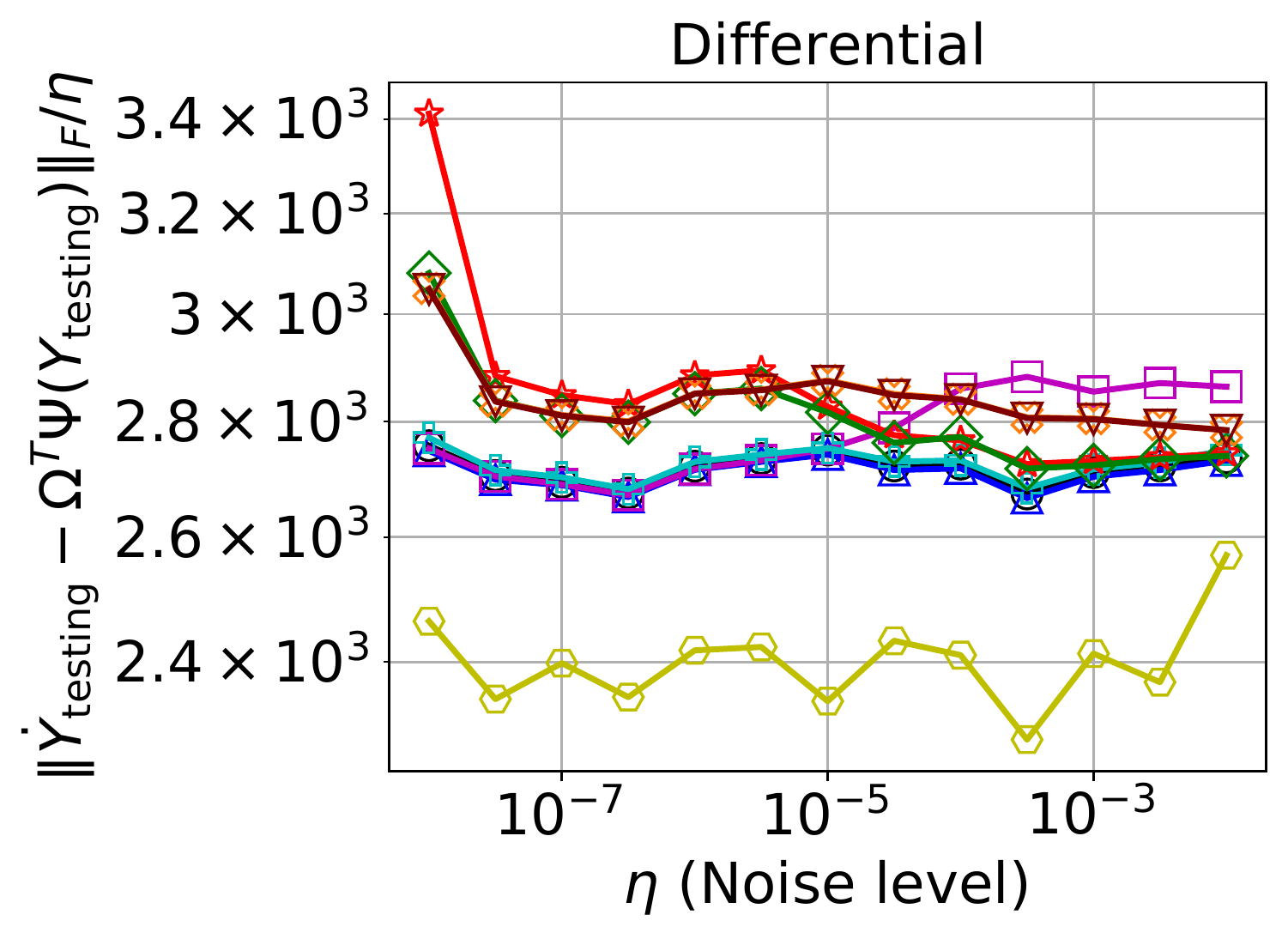} }\label{fig:appx:Kura10:diff:val}}%
    %\qquad
    \hspace{\fill}
    \subfloat[]{{\includegraphics[width=3.85cm]{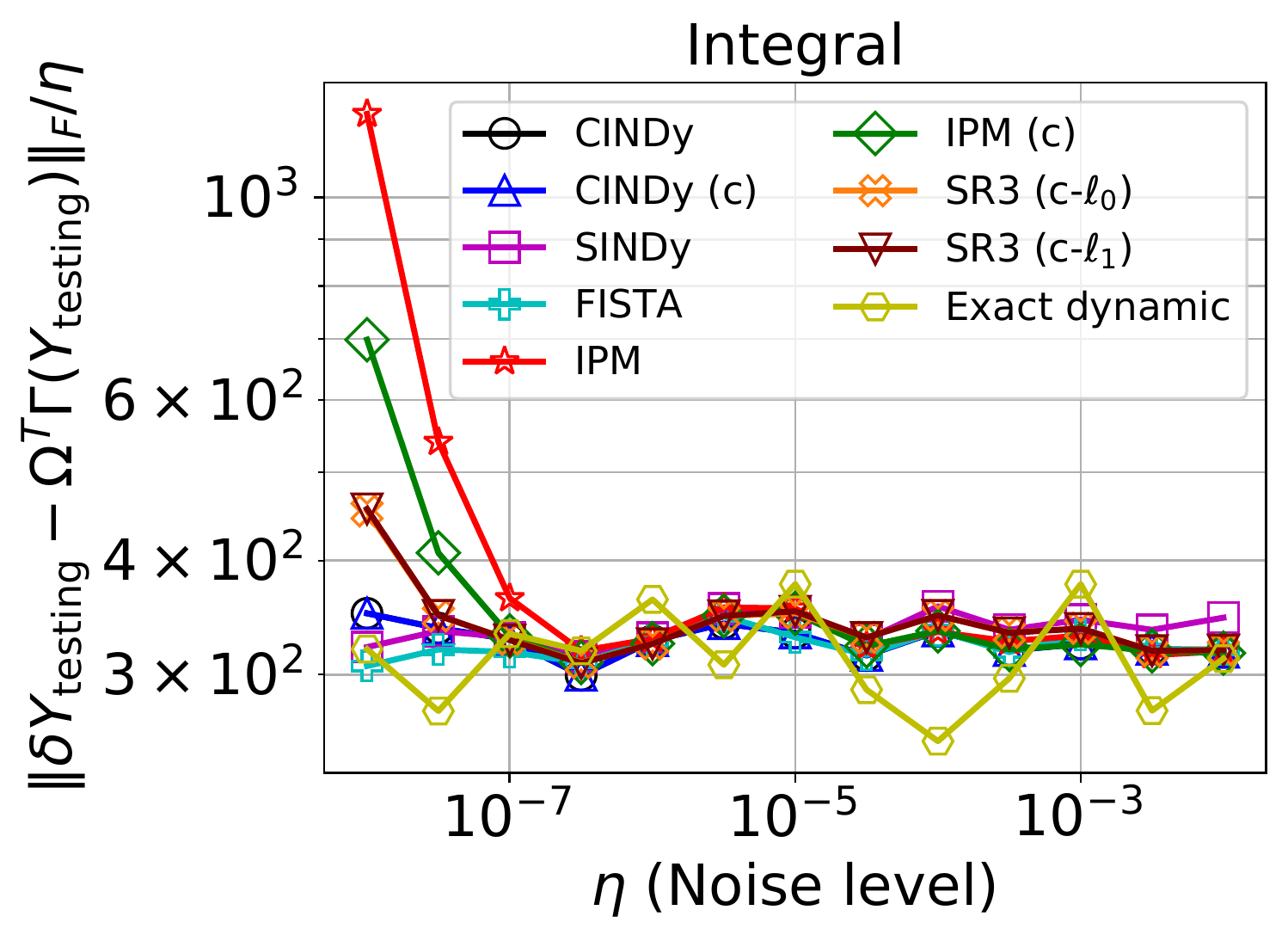} }\label{fig:appx:Kura10:int:val}}%
    \hspace*{\fill}
    \caption{Kuramoto: Evaluation of $\norm{\ddot{Y}_{\mathrm{training}} -  \Omega^T \Psi(Y_{\mathrm{training}})}_F$ \protect\subref{fig:appx:Kura10:diff:training} for the differential formulation, and $\norm{ \delta \dot{Y}_{\mathrm{training}} -  \Omega^T \Gamma(Y_{\mathrm{training}})}_F$ \protect\subref{fig:appx:Kura10:int:training} for the integral formulation for the experiments with $d = 10$, and evaluation of $\norm{\ddot{Y}_{\mathrm{validation}} -  \Omega^T \Psi(Y_{\mathrm{validation}})}_F$ \protect\subref{fig:appx:Kura10:diff:val} for the differential formulation, and $\norm{ \delta \dot{Y}_{\mathrm{validation}} -  \Omega^T \Gamma(Y_{\mathrm{validation}})}_F$ \protect\subref{fig:appx:Kura10:int:val} for the integral formulation for the experiments with $d = 10$.}%
    \label{fig:appx:kura10}%
\end{figure*}

The images shown in Figure~\ref{fig:appx:kura5} and \ref{fig:appx:kura10} show the objective function evaluation for the constrained version of the Kuramoto LASSO problem with $d = 5$ and $d = 10$, respectively. The objective function is evaluated for the different methods, for both the training data and the testing data and the differential and integral formulation.

\begin{figure*}[h!]
    \centering
    \vspace{-10pt}
    \hspace{\fill}
    \subfloat[]{{\includegraphics[width=3.95cm]{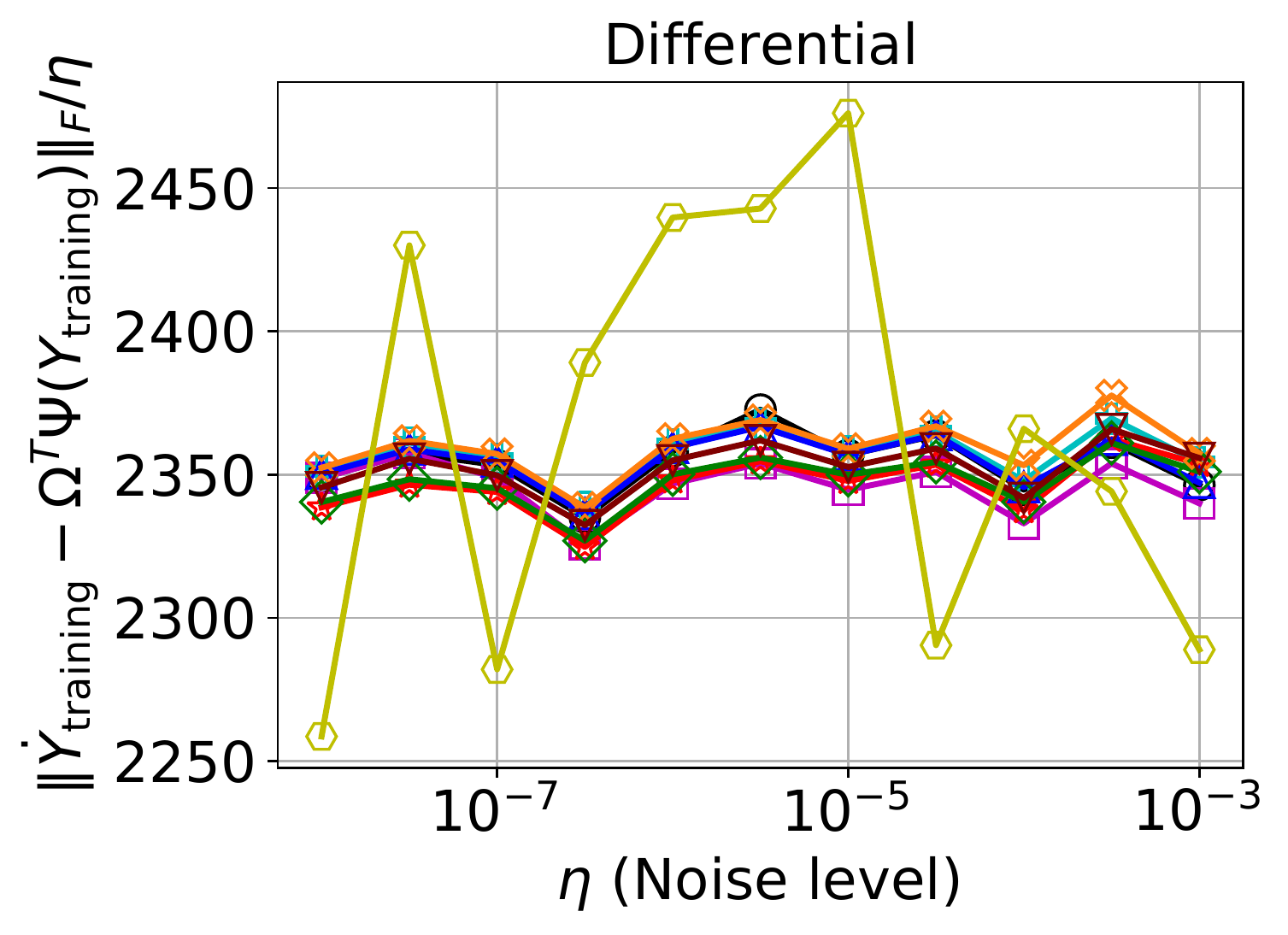} }\label{fig:appx:FPUT5:diff:training}}%
    %\qquad
    \hspace{\fill}
    \subfloat[]{{\includegraphics[width=3.85cm]{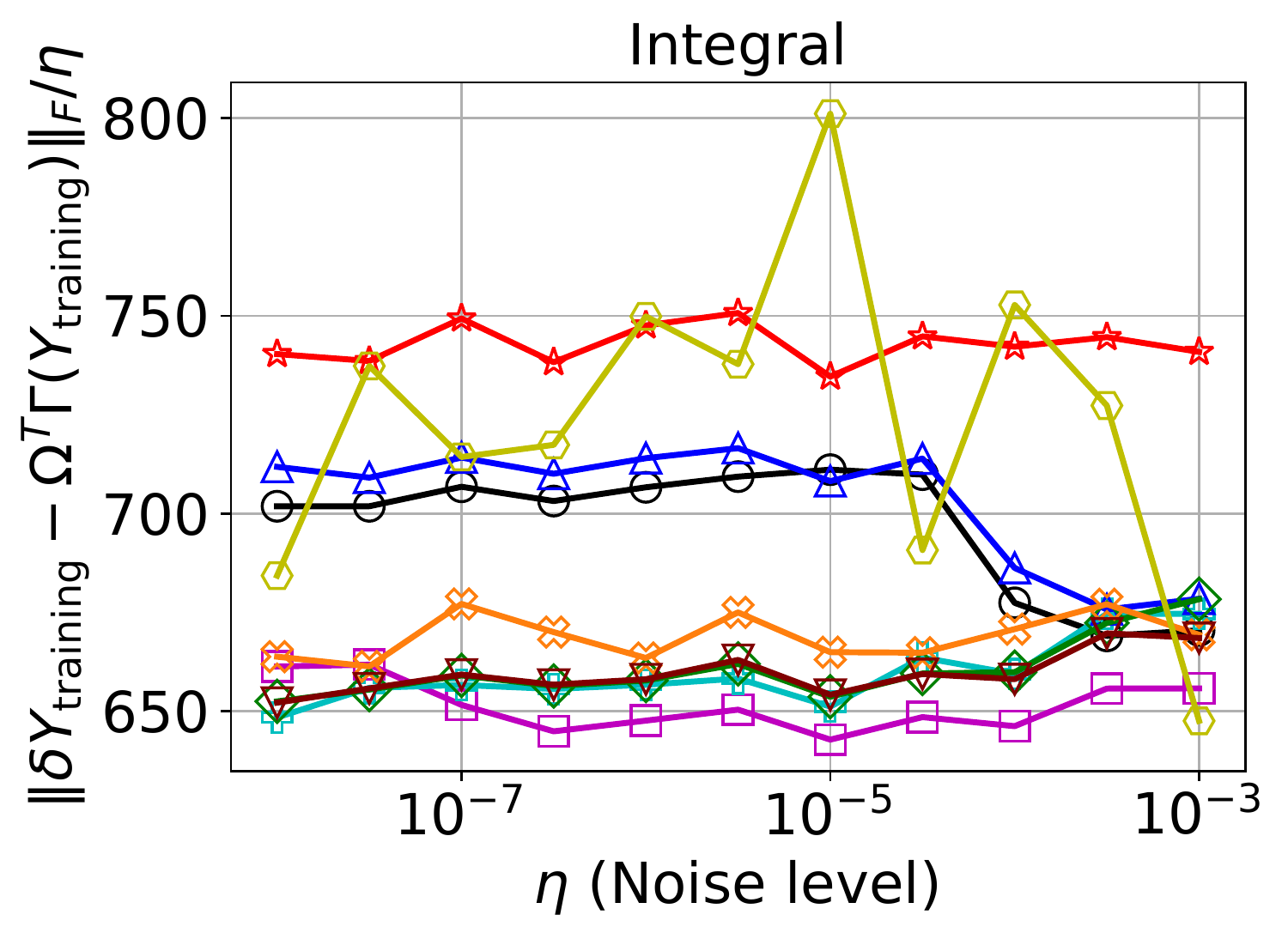} }\label{fig:appx:FPUT5:int:training}}%
    \hspace{\fill}
    \subfloat[]{{\includegraphics[width=3.9cm]{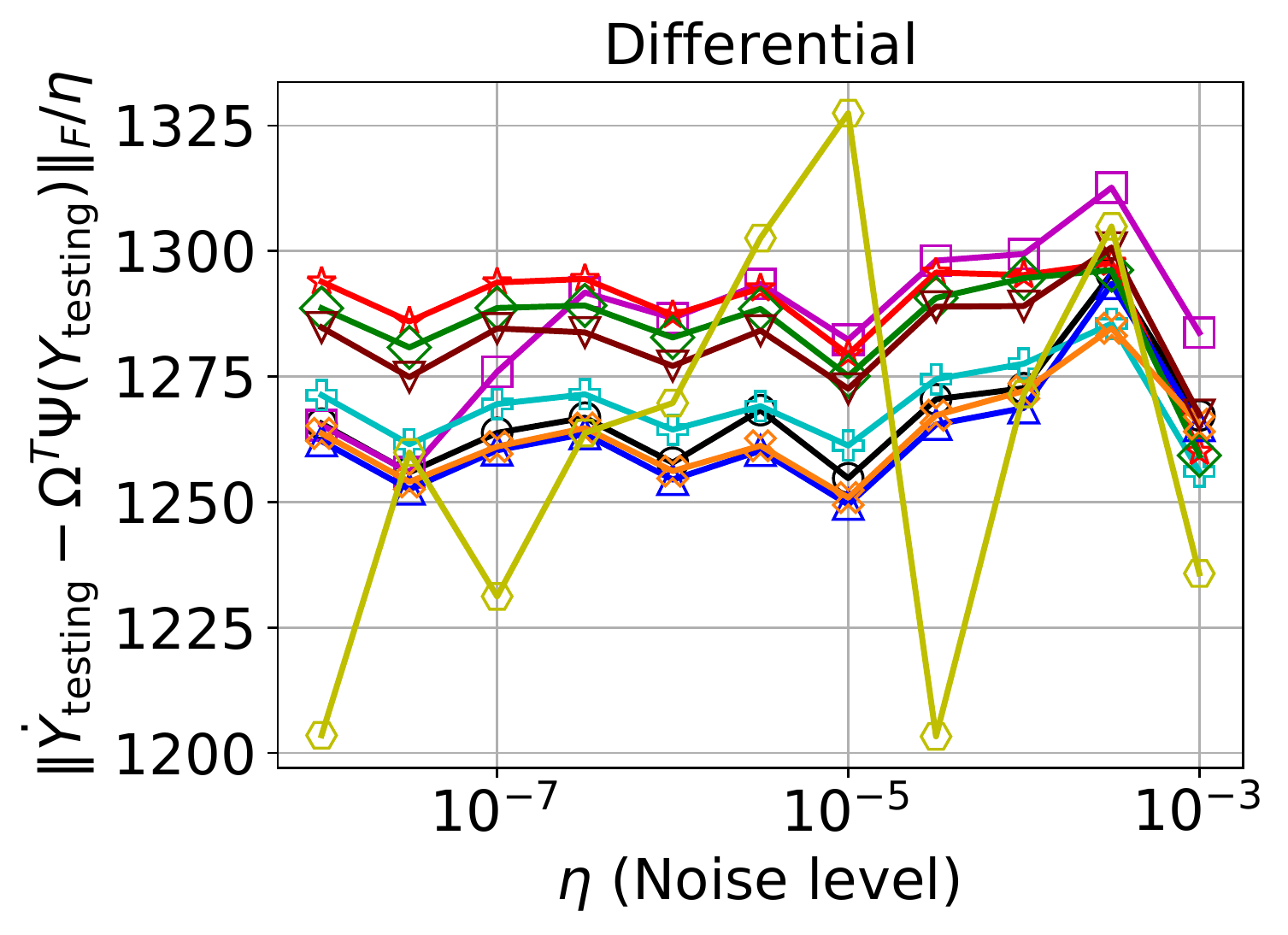} }\label{fig:appx:FPUT5:diff:val}}%
    %\qquad
    \hspace{\fill}
    \subfloat[]{{\includegraphics[width=3.85cm]{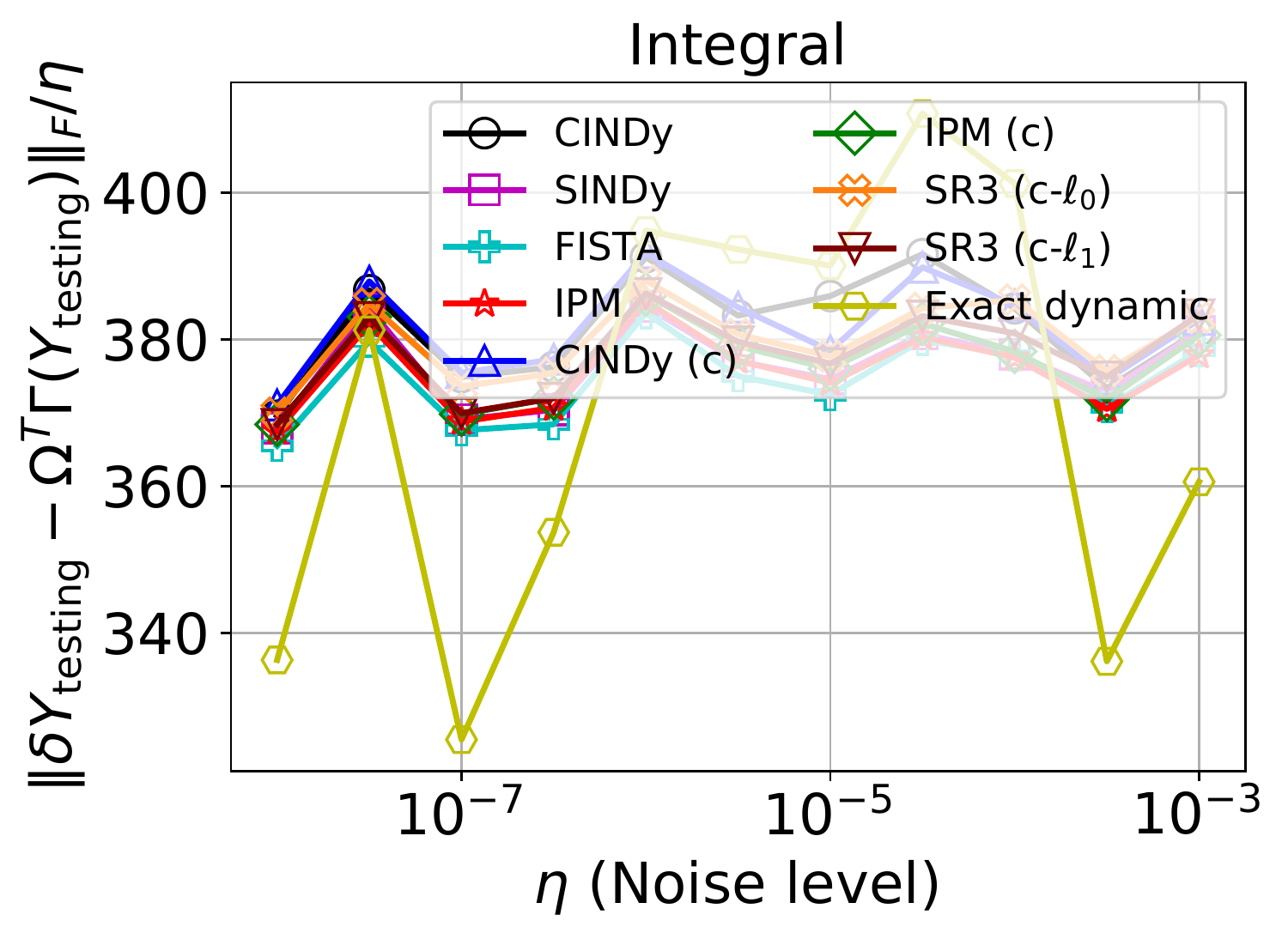} }\label{fig:appx:FPUT5:int:val}}%
    \hspace*{\fill}
    \caption{Fermi-Pasta-Ulam-Tsingou: Evaluation of $\norm{\ddot{Y}_{\mathrm{training}} -  \Omega^T \Psi(Y_{\mathrm{training}})}_F$ \protect\subref{fig:appx:FPUT5:diff:training} for the differential formulation and $\norm{ \delta \dot{Y}_{\mathrm{training}} -  \Omega^T \Gamma(Y_{\mathrm{training}})}_F$ \protect\subref{fig:appx:FPUT5:int:training} for the integral formulation with $d = 5$, and evaluation of $\norm{\ddot{Y}_{\mathrm{validation}} -  \Omega^T \Psi(Y_{\mathrm{validation}})}_F$ \protect\subref{fig:appx:FPUT5:diff:val} for the differential formulation, and $\norm{ \delta \dot{Y}_{\mathrm{validation}} -  \Omega^T \Gamma(Y_{\mathrm{validation}})}_F$ \protect\subref{fig:appx:FPUT5:int:val} for the integral formulation with $d= 5$.}%
    \label{fig:appx:FPUT5}%
\end{figure*}

\begin{figure*}[h!]
    \centering
    \vspace{-10pt}
    \hspace{\fill}
    \subfloat[]{{\includegraphics[width=3.95cm]{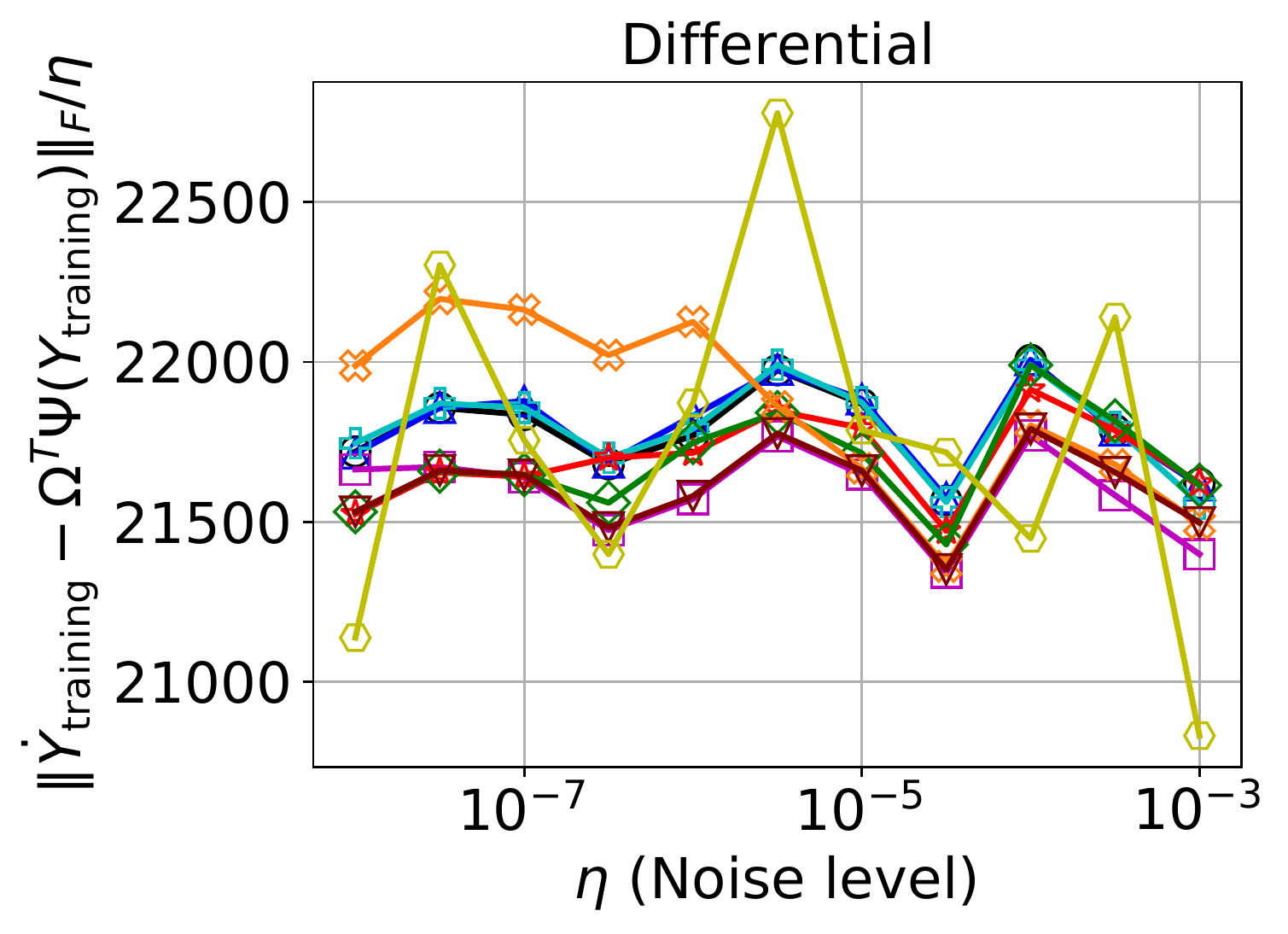} }\label{fig:appx:FPUT10:diff:training}}%
    %\qquad
    \hspace{\fill}
    \subfloat[]{{\includegraphics[width=3.85cm]{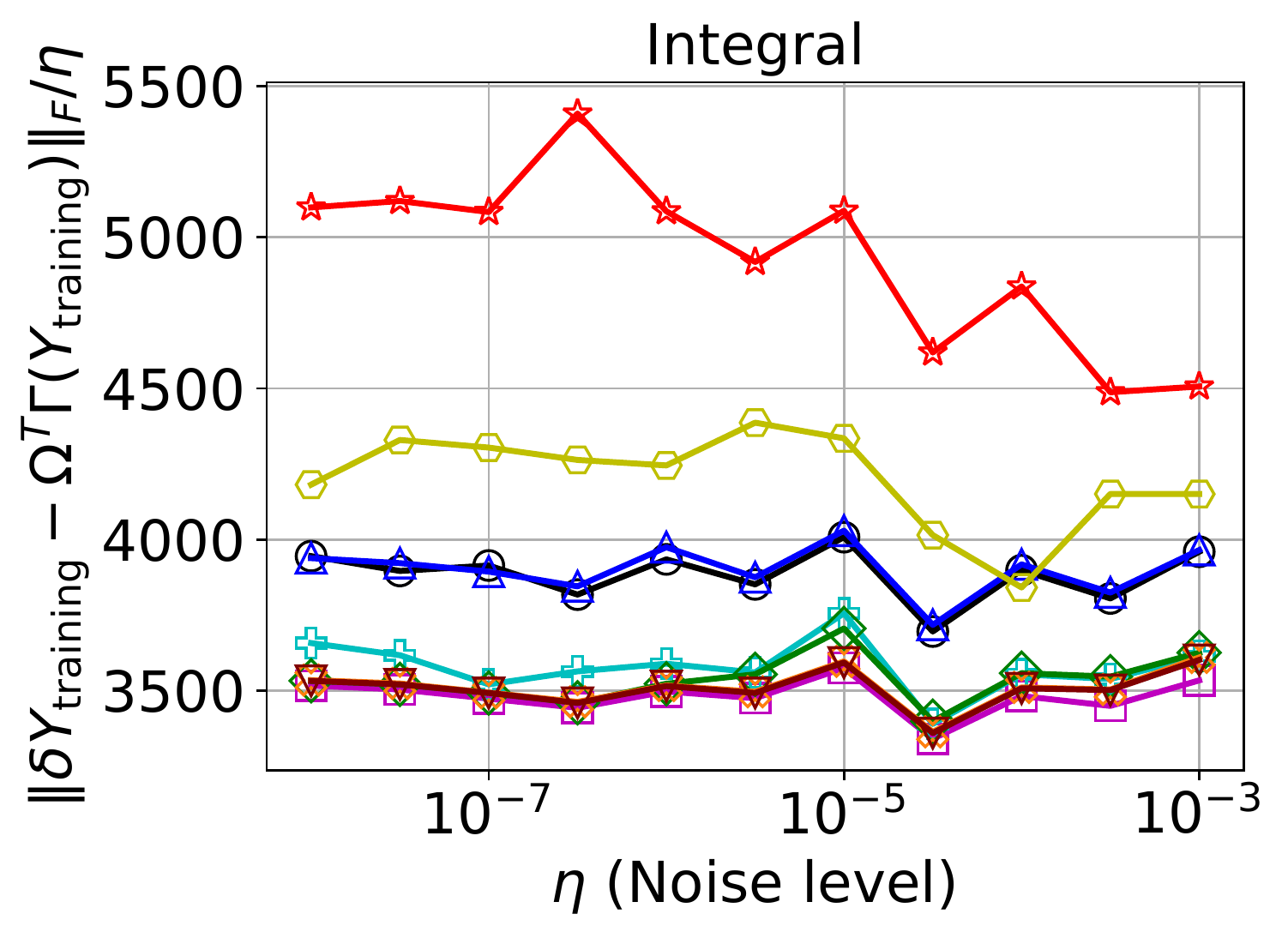} }\label{fig:appx:FPUT10:int:training}}%
    \hspace{\fill}
    \subfloat[]{{\includegraphics[width=3.9cm]{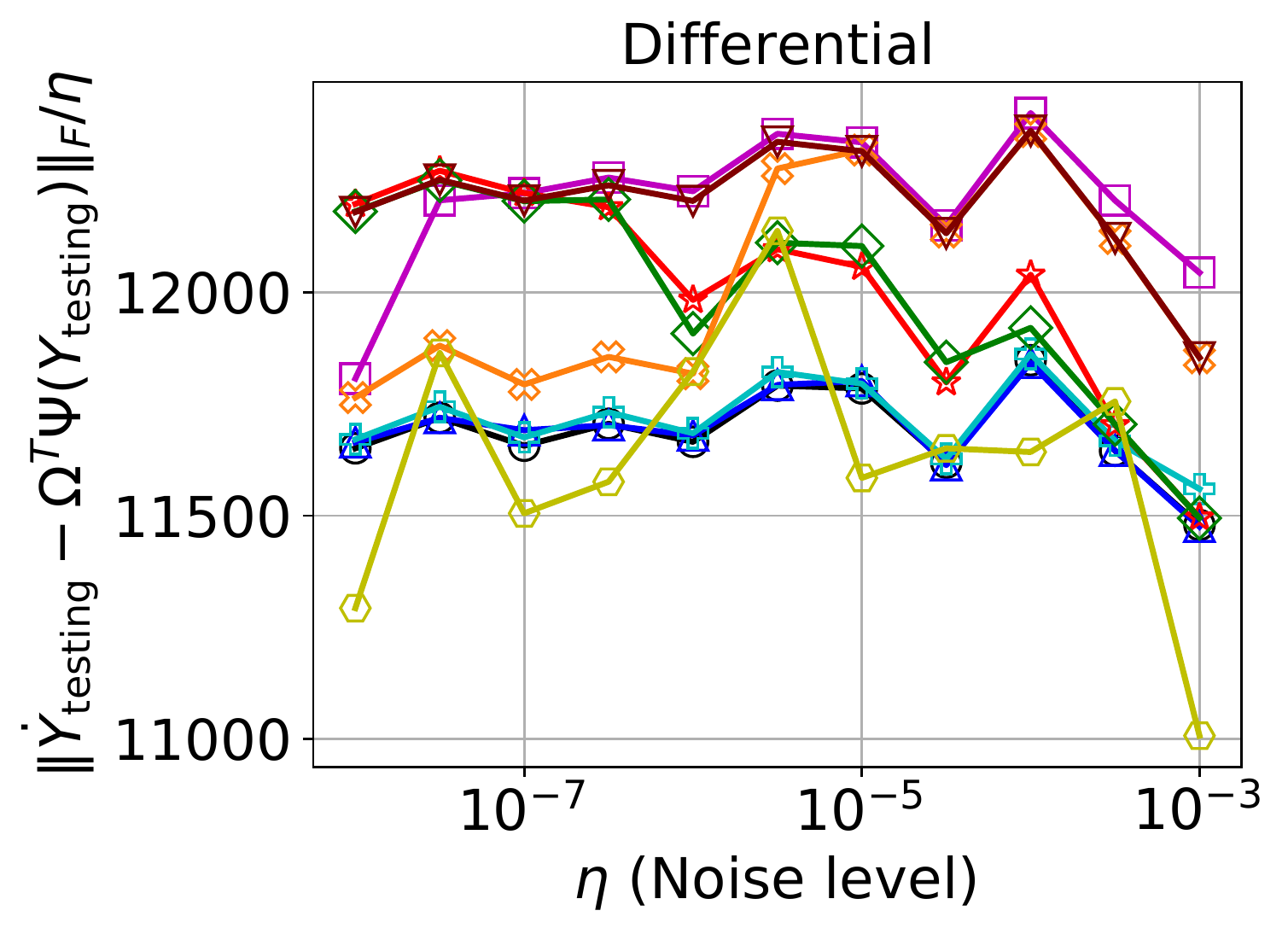} }\label{fig:appx:FPUT10:diff:val}}%
    %\qquad
    \hspace{\fill}
    \subfloat[]{{\includegraphics[width=3.85cm]{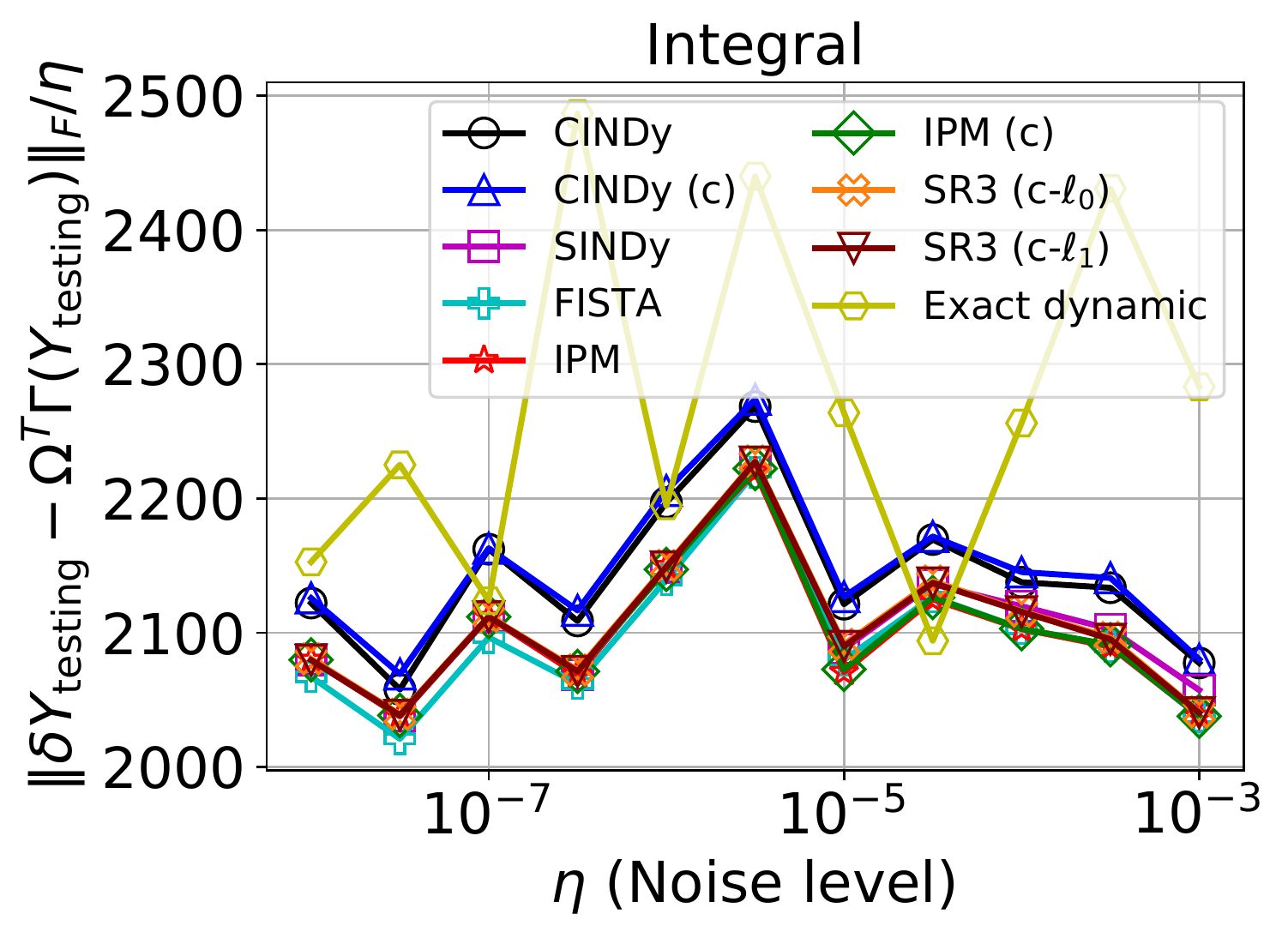} }\label{fig:appx:FPUT10:int:val}}%
    \hspace*{\fill}
    \caption{Fermi-Pasta-Ulam-Tsingou: Evaluation of $\norm{\ddot{Y}_{\mathrm{training}} -  \Omega^T \Psi(Y_{\mathrm{training}})}_F$ \protect\subref{fig:appx:FPUT10:diff:training} for the differential formulation and $\norm{ \delta \dot{Y}_{\mathrm{training}} -  \Omega^T \Gamma(Y_{\mathrm{training}})}_F$ \protect\subref{fig:appx:FPUT10:int:training} for the integral formulation with $d = 10$, and evaluation of $\norm{\ddot{Y}_{\mathrm{validation}} -  \Omega^T \Psi(Y_{\mathrm{validation}})}_F$ \protect\subref{fig:appx:FPUT10:diff:val} for the differential formulation, and $\norm{ \delta \dot{Y}_{\mathrm{validation}} -  \Omega^T \Gamma(Y_{\mathrm{validation}})}_F$ \protect\subref{fig:appx:FPUT10:int:val} for the integral formulation with $d= 10$.}%
    \label{fig:appx:FPUT10}%
\end{figure*}

The images shown in Figure~\ref{fig:appx:FPUT5} and \ref{fig:appx:FPUT10} show the objective function evaluation for the constrained version of the Fermi-Pasta-Ulam-Tsingou LASSO problem with $d = 5$ and $d = 10$, respectively. The objective function is evaluated for the different methods, for both the training data and the testing data and the differential and integral formulation.

\begin{figure*}[h!]
    \centering
    \vspace{-10pt}
    \hspace{\fill}
    \subfloat[]{{\includegraphics[width=3.95cm]{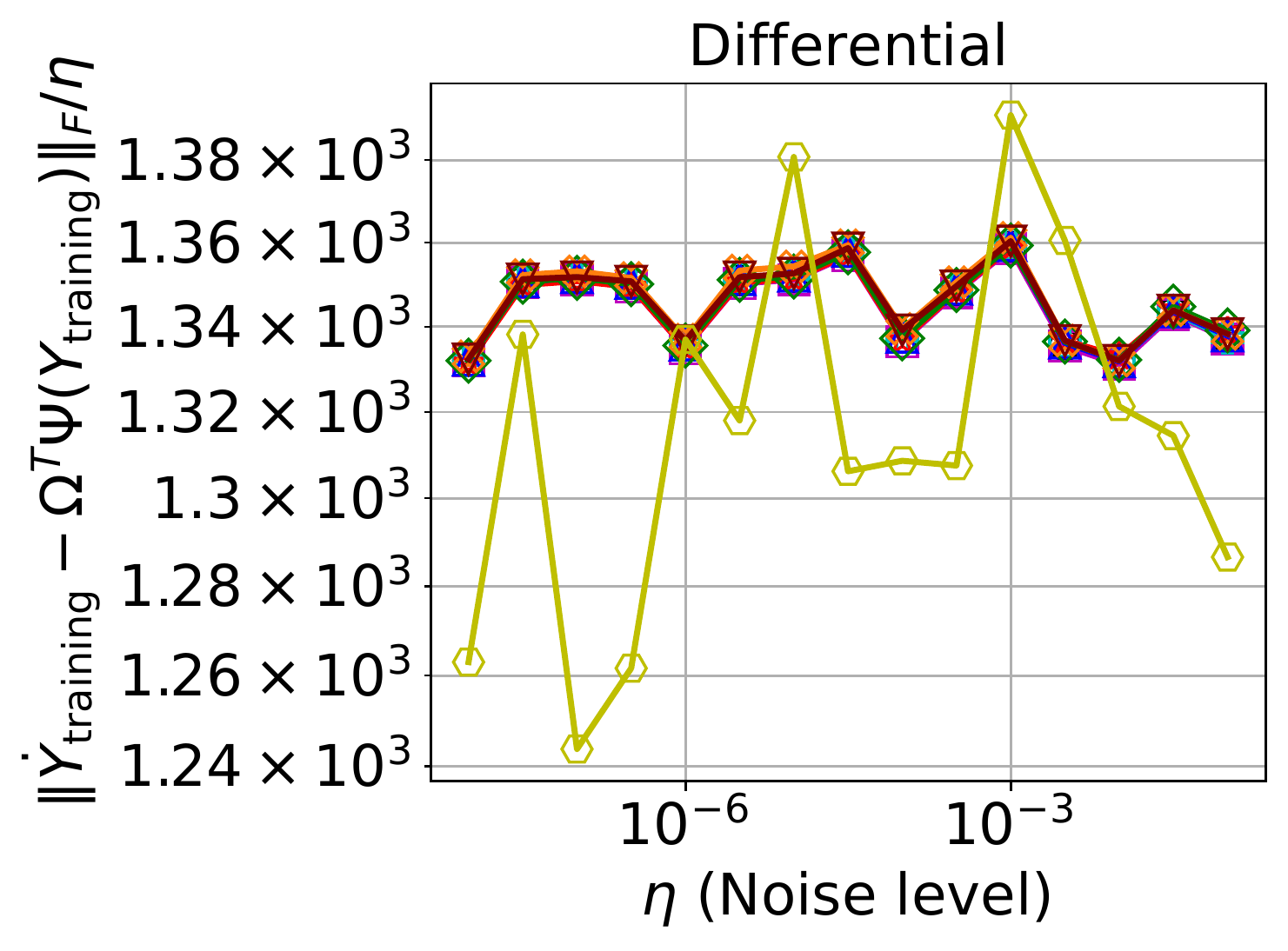} }\label{fig:appx:MM4:diff:training}}%
    %\qquad
    \hspace{\fill}
    \subfloat[]{{\includegraphics[width=3.85cm]{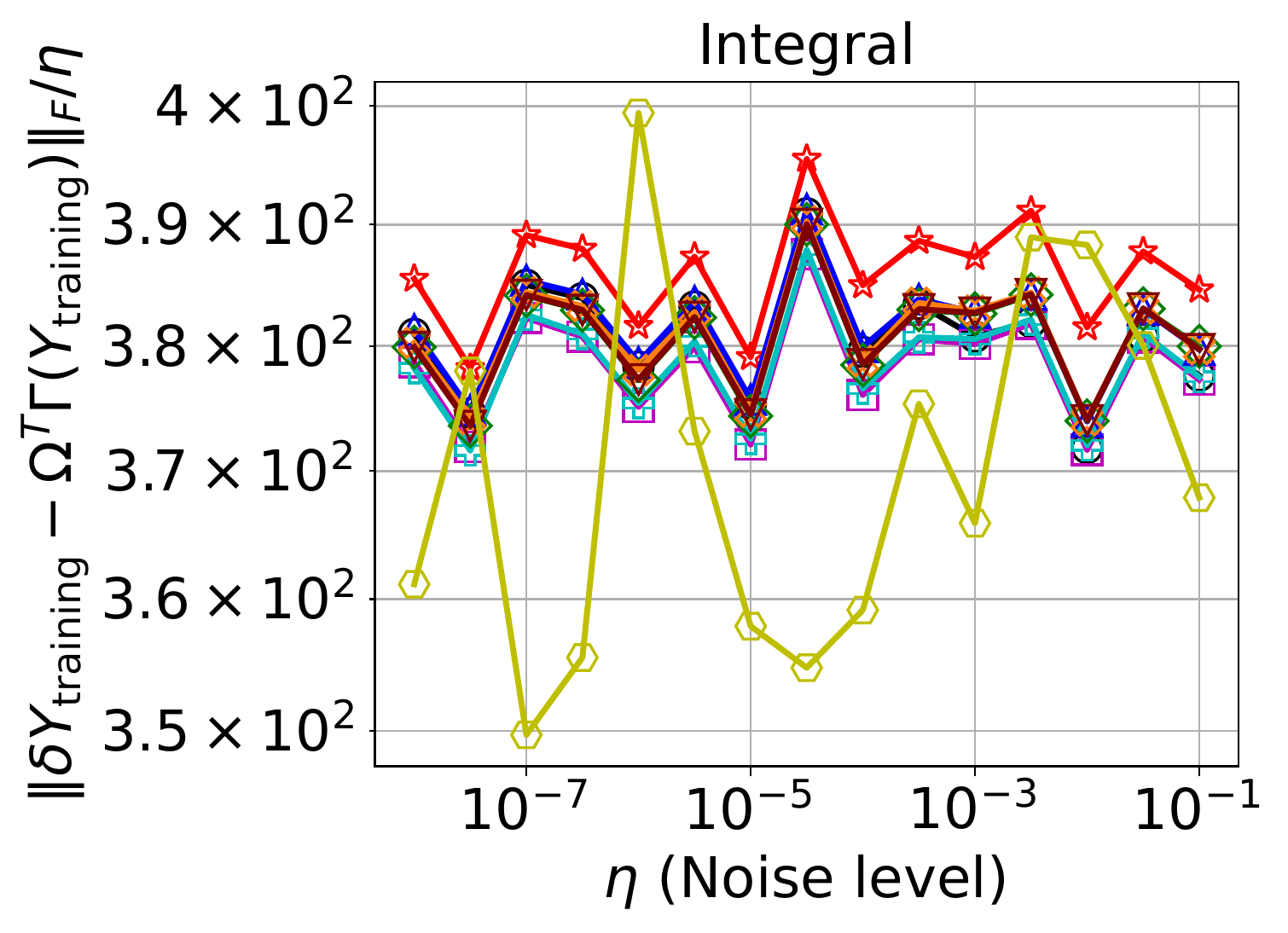} }\label{fig:appx:MM4:int:training}}%
    \hspace{\fill}
    \subfloat[]{{\includegraphics[width=3.9cm]{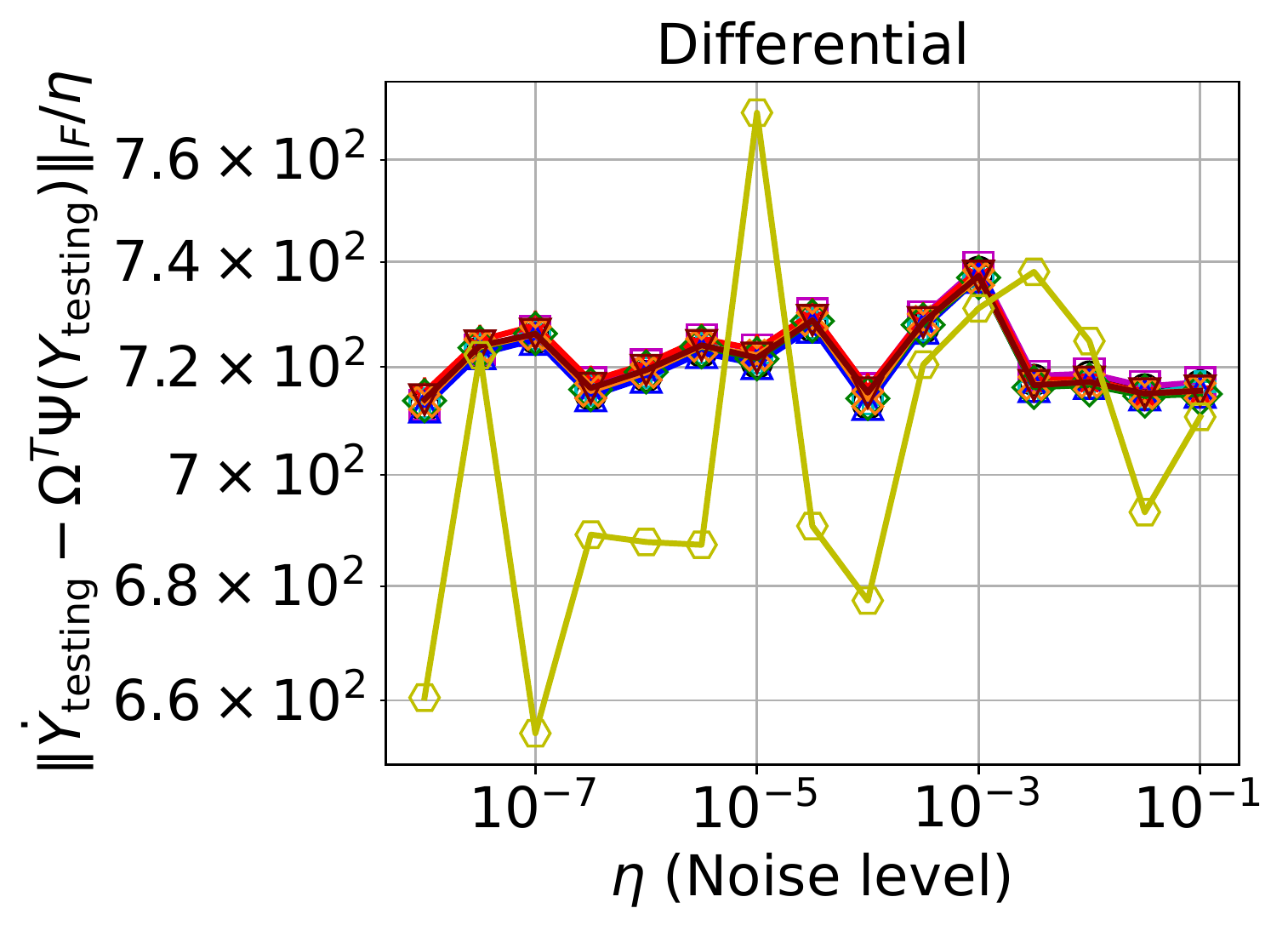} }\label{fig:appx:MM4:diff:val}}%
    %\qquad
    \hspace{\fill}
    \subfloat[]{{\includegraphics[width=3.85cm]{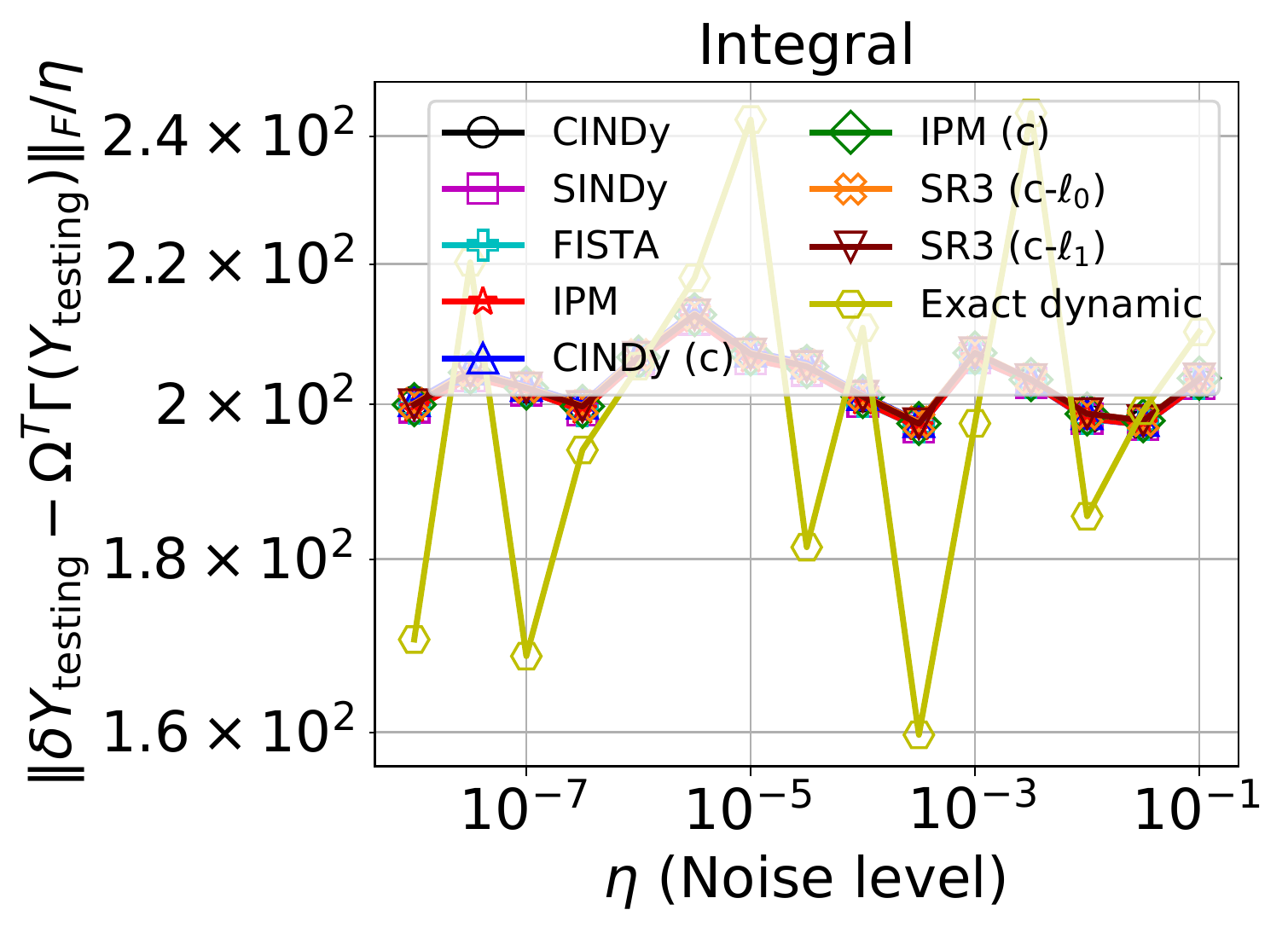} }\label{fig:appx:MM4:int:val}}%
    \hspace*{\fill}
    \caption{Michaelis-Menten: Evaluation of $\norm{\ddot{Y}_{\mathrm{training}} -  \Omega^T \Psi(Y_{\mathrm{training}})}_F$ \protect\subref{fig:appx:MM4:diff:training} for the differential formulation, and $\norm{ \delta \dot{Y}_{\mathrm{training}} -  \Omega^T \Gamma(Y_{\mathrm{training}})}_F$ \protect\subref{fig:appx:MM4:int:training} for the integral formulation with $d = 4$, and evaluation of $\norm{\ddot{Y}_{\mathrm{validation}} -  \Omega^T \Psi(Y_{\mathrm{validation}})}_F$ \protect\subref{fig:appx:MM4:diff:val} for the differential formulation, and $\norm{ \delta \dot{Y}_{\mathrm{validation}} -  \Omega^T \Gamma(Y_{\mathrm{validation}})}_F$ \protect\subref{fig:appx:MM4:int:val} for the integral formulation with $d= 4$.}%
    \label{fig:appx:MM}%
\end{figure*}

The images shown in Figure~\ref{fig:appx:MM} show the objective function evaluation for the constrained version of the Michaelis-Menten LASSO problem with $d = 4$. The objective function is evaluated for the different methods, for both the training data and the testing data and the differential and integral formulation.

\newpage

\section{Sample efficiency} \label{appx:section:sample_efficiency}

\begin{figure*}[]
    \centering
    \vspace{-10pt}
    \hspace{\fill}
    \subfloat{{\includegraphics[width=8cm]{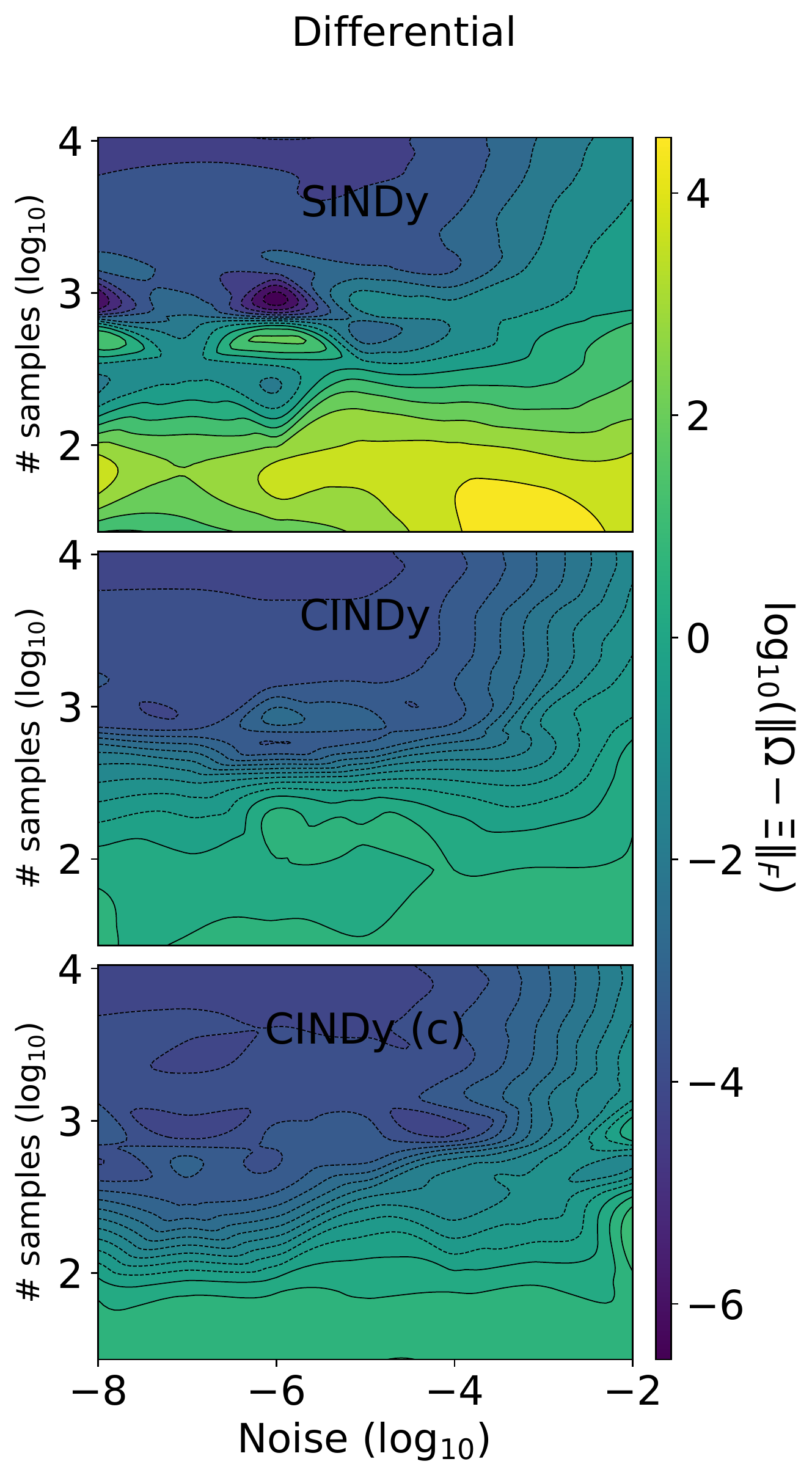} }\label{fig:appx:kuramoto:diff:dim5_sampleeff_acc}}%
    %\qquad
    \hspace{\fill}
    \subfloat{{\includegraphics[width=8cm]{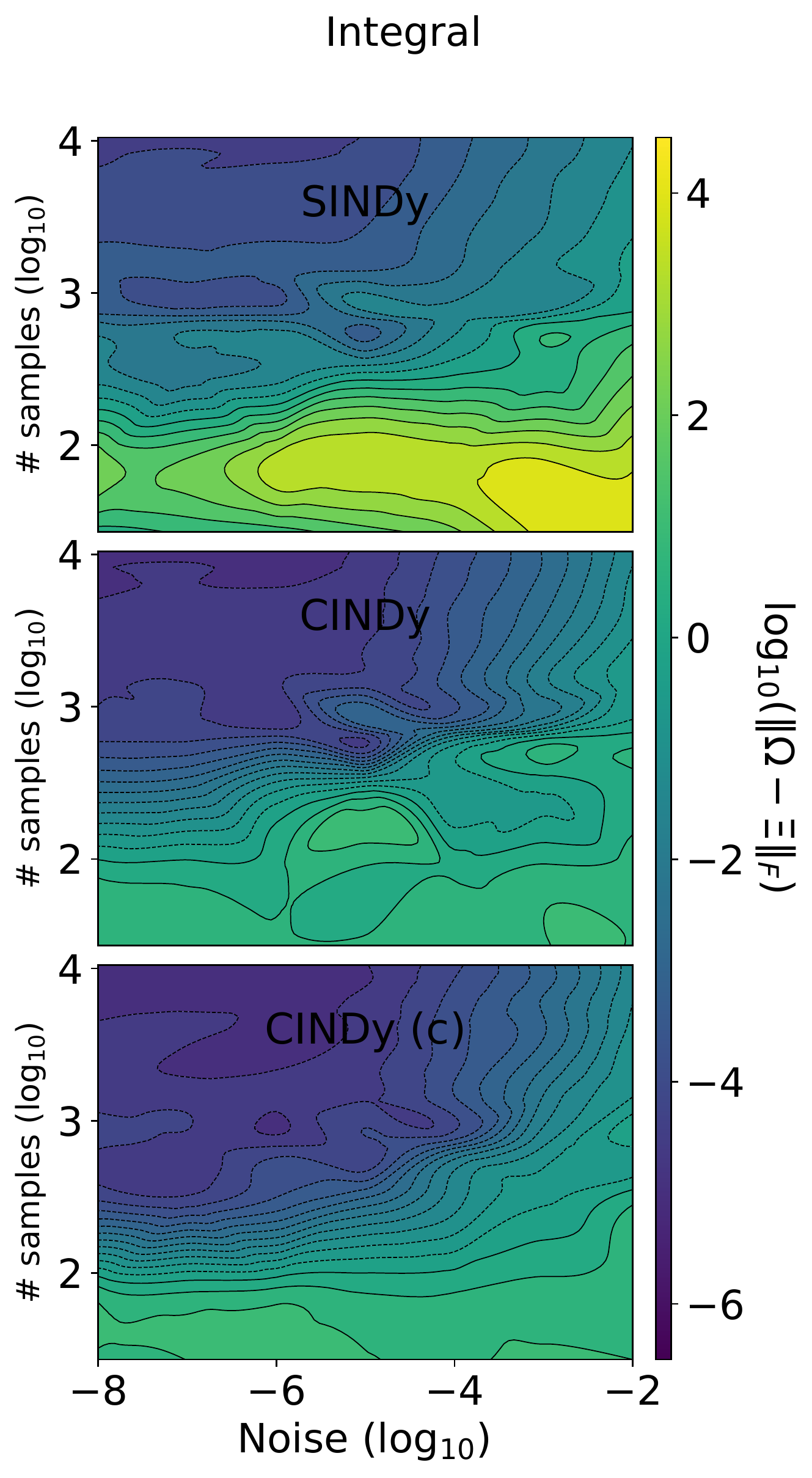} }\label{fig:appx:kuramoto:int:dim5_sampleeff_acc}}%
    \hspace{\fill}
    \caption{\textbf{Sample efficiency of the sparse recovery of the
        Kuramoto model: } Algorithm comparison in terms of $\mathcal{E}_R$ for a Kuramoto model of
      dimension $d=5$ for the differential formulation (left) and the integral formulation (right).}%
    \label{fig:appx:kuramoto_sampleeff}%
  \end{figure*}

\begin{figure*}[]
    \centering
    \vspace{-10pt}
    \hspace{\fill}
    \subfloat{{\includegraphics[width=8cm]{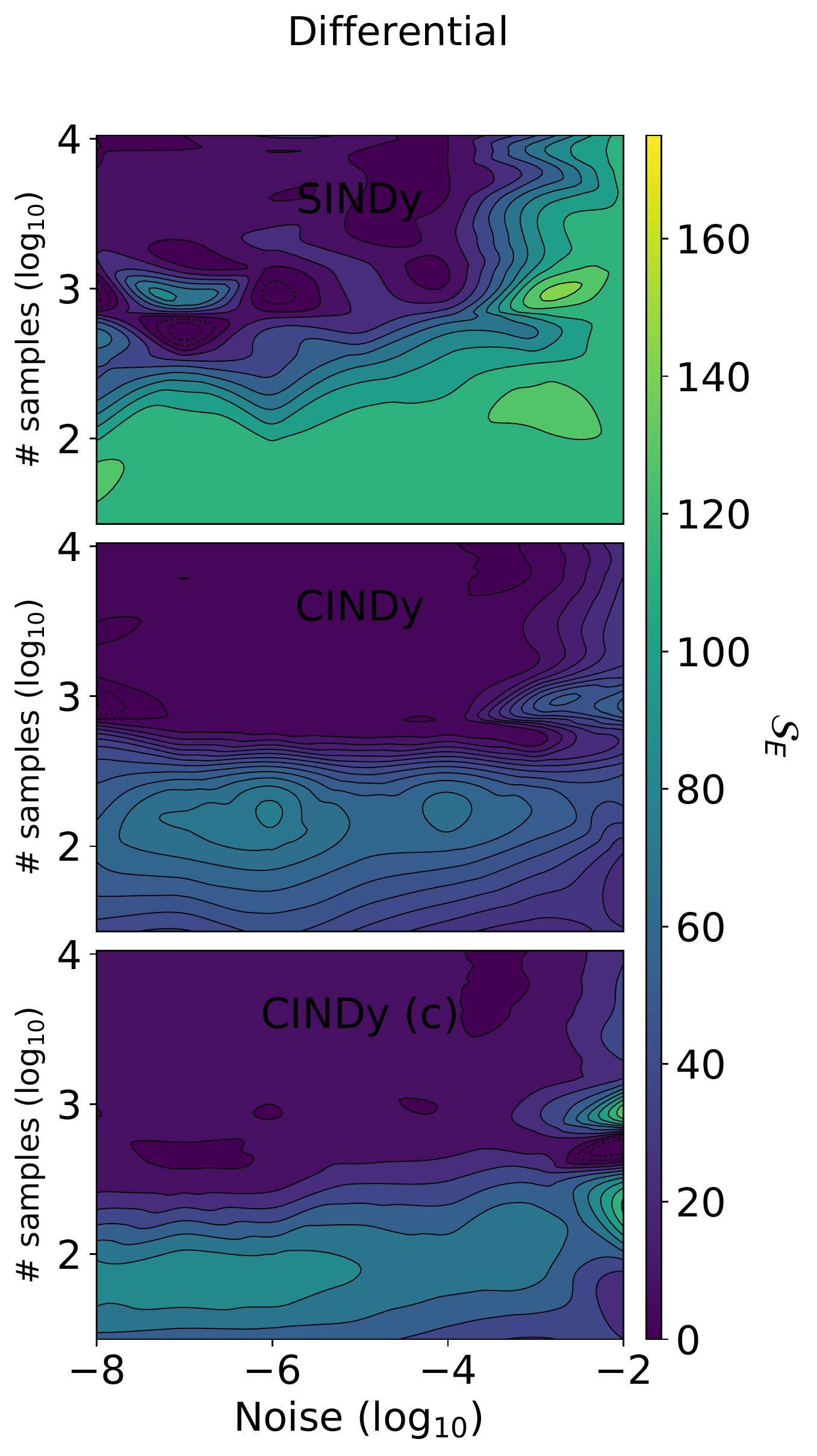} }\label{fig:appx:kuramoto:diff:dim5_sampleeff_extra}}%
    %\qquad
    \hspace{\fill}
    \subfloat{{\includegraphics[width=8cm]{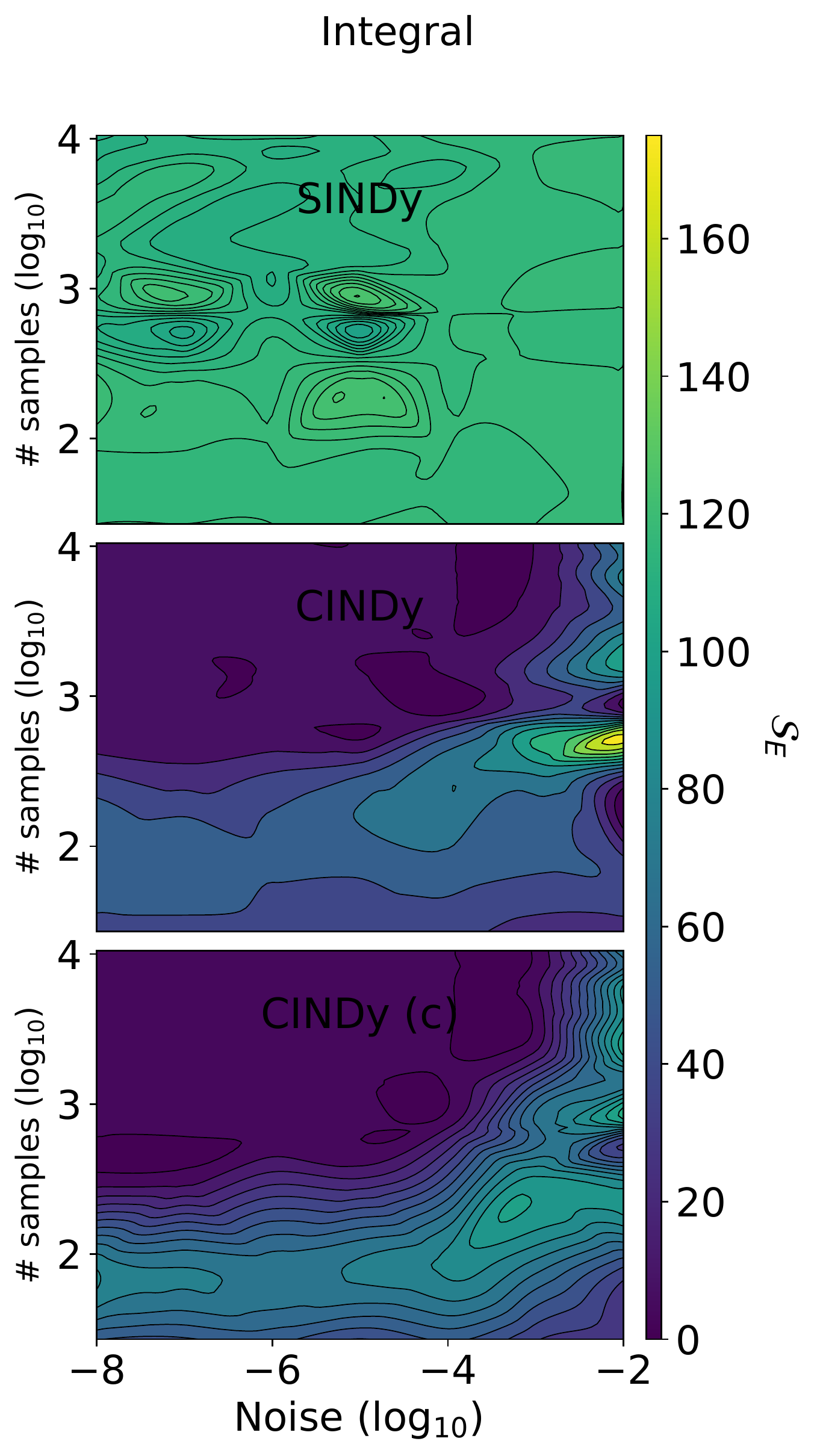} }\label{fig:appx:kuramoto:int:dim5_sampleeff_extra}}%
    \hspace{\fill}
    \caption{\textbf{Sample efficiency of the sparse recovery of the
        Kuramoto model: } Algorithm comparison in terms of $\mathcal{S}_E$ for a Kuramoto model of
      dimension $d=5$ for the differential formulation (left) and the integral formulation (right).}%
    \label{fig:appx:kuramoto_sampleeff_extra}%
  \end{figure*}

\begin{figure*}[]
    \centering
    \vspace{-10pt}
    \hspace{\fill}
    \subfloat{{\includegraphics[width=8cm]{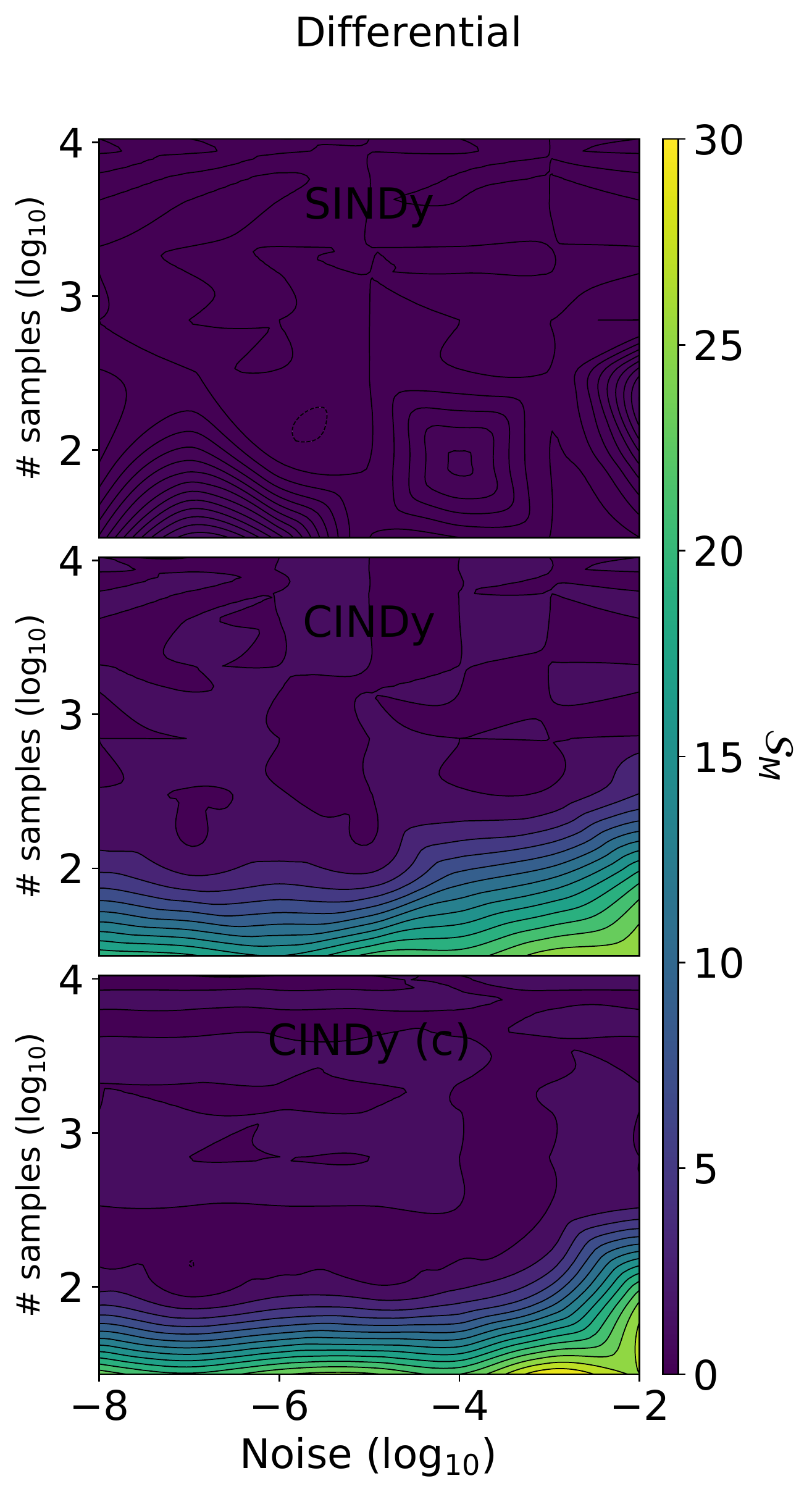} }\label{fig:appx:kuramoto:diff:dim5_sampleeff_missing}}%
    %\qquad
    \hspace{\fill}
    \subfloat{{\includegraphics[width=8cm]{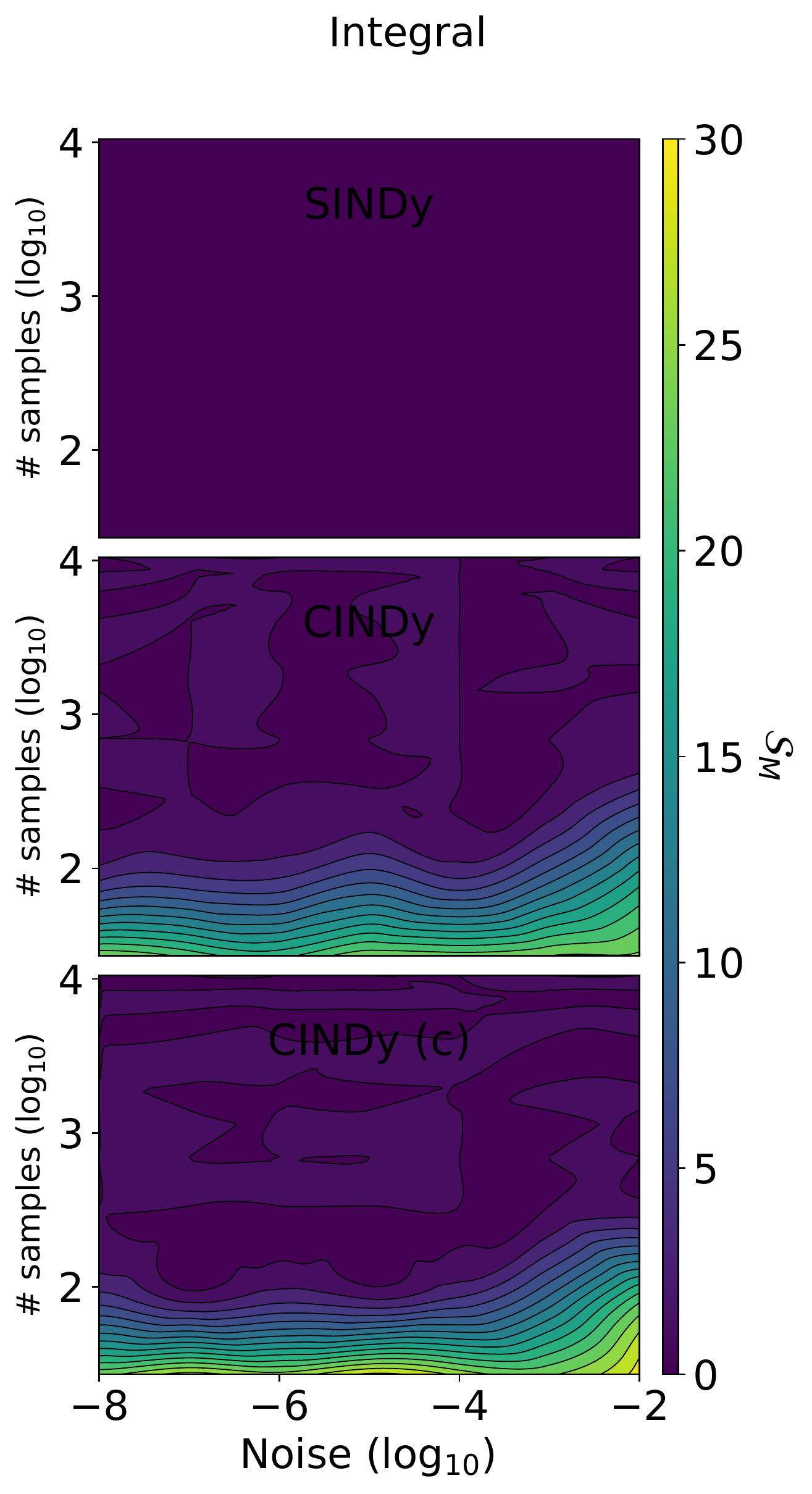} }\label{fig:appx:kuramoto:int:dim5_sampleeff_missing}}%
    \hspace{\fill}
    \caption{\textbf{Sample efficiency of the sparse recovery of the
        Kuramoto model: } Algorithm comparison in terms of $\mathcal{S}_M$ for a Kuramoto model of
      dimension $d=5$ for the differential formulation (left) and the integral formulation (right).}%
    \label{fig:appx:kuramoto_sampleeff_missing}%
  \end{figure*}

The images shown in Figure~\ref{fig:appx:kuramoto_sampleeff}, \ref{fig:appx:kuramoto_sampleeff_extra}, and \ref{fig:appx:kuramoto_sampleeff_missing} show the evolution of $\mathcal{E}_R$, $\mathcal{S}_E$ and $\mathcal{S}_M$ as we vary the number of training data points when learning the Kuramoto dynamic ($d = 5$) with the dictionary described in Section~\ref{section:kura}. The images show the resulting metrics when generating $50$ data points per experiment and using local polynomial interpolation of degree $8$ to compute the derivatives and the integrals. Note that the artifacts present in the images are caused by the use of cubic interpolation to generate the images which results in an oscillatory behaviour in the heat maps.

\begin{figure*}[]
    \centering
    \vspace{-10pt}
    \hspace{\fill}
    \subfloat{{\includegraphics[width=8cm]{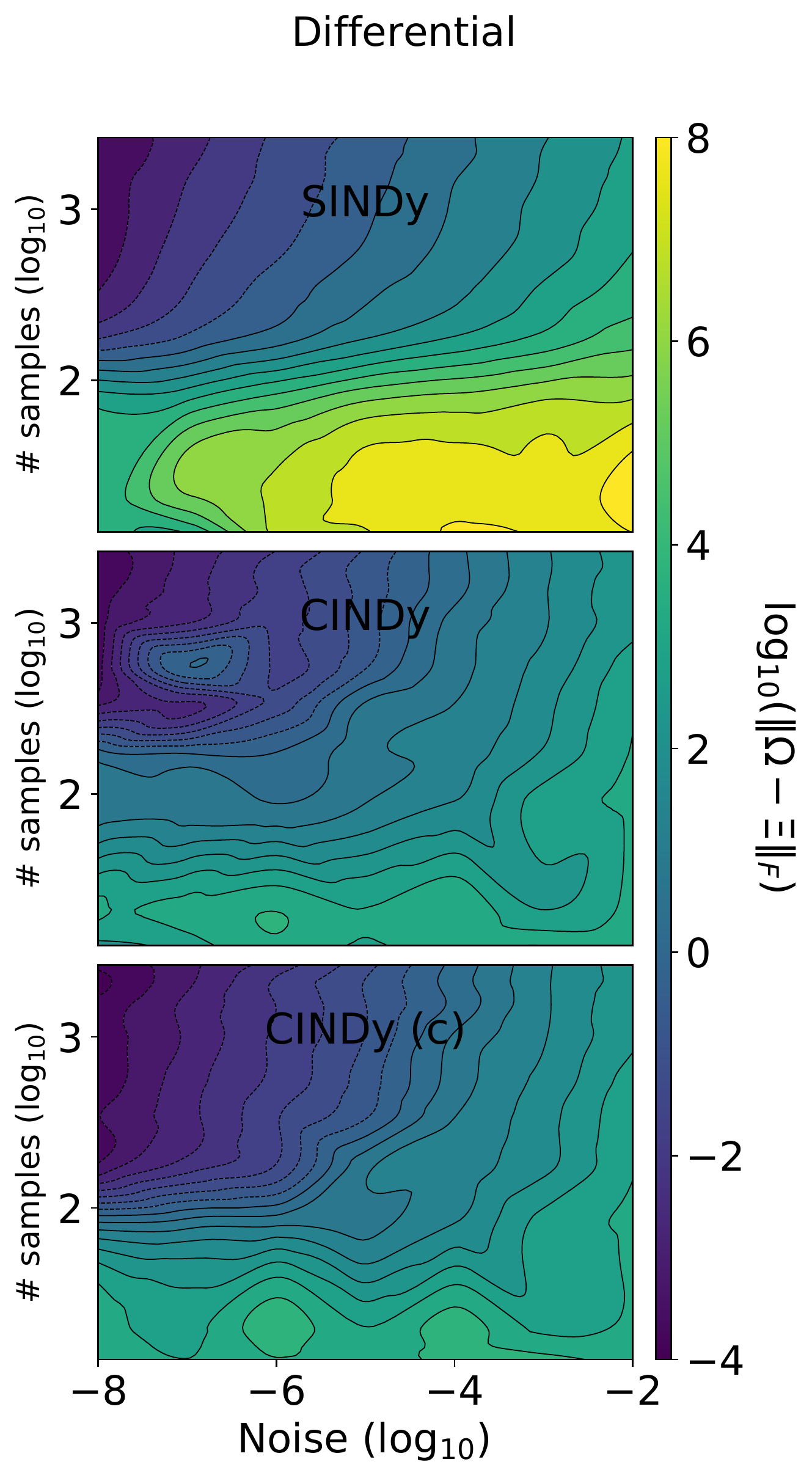} }\label{fig:appx:FPUT:diff:dim5_sampleeff_acc}}%
    %\qquad
    \hspace{\fill}
    \subfloat{{\includegraphics[width=8cm]{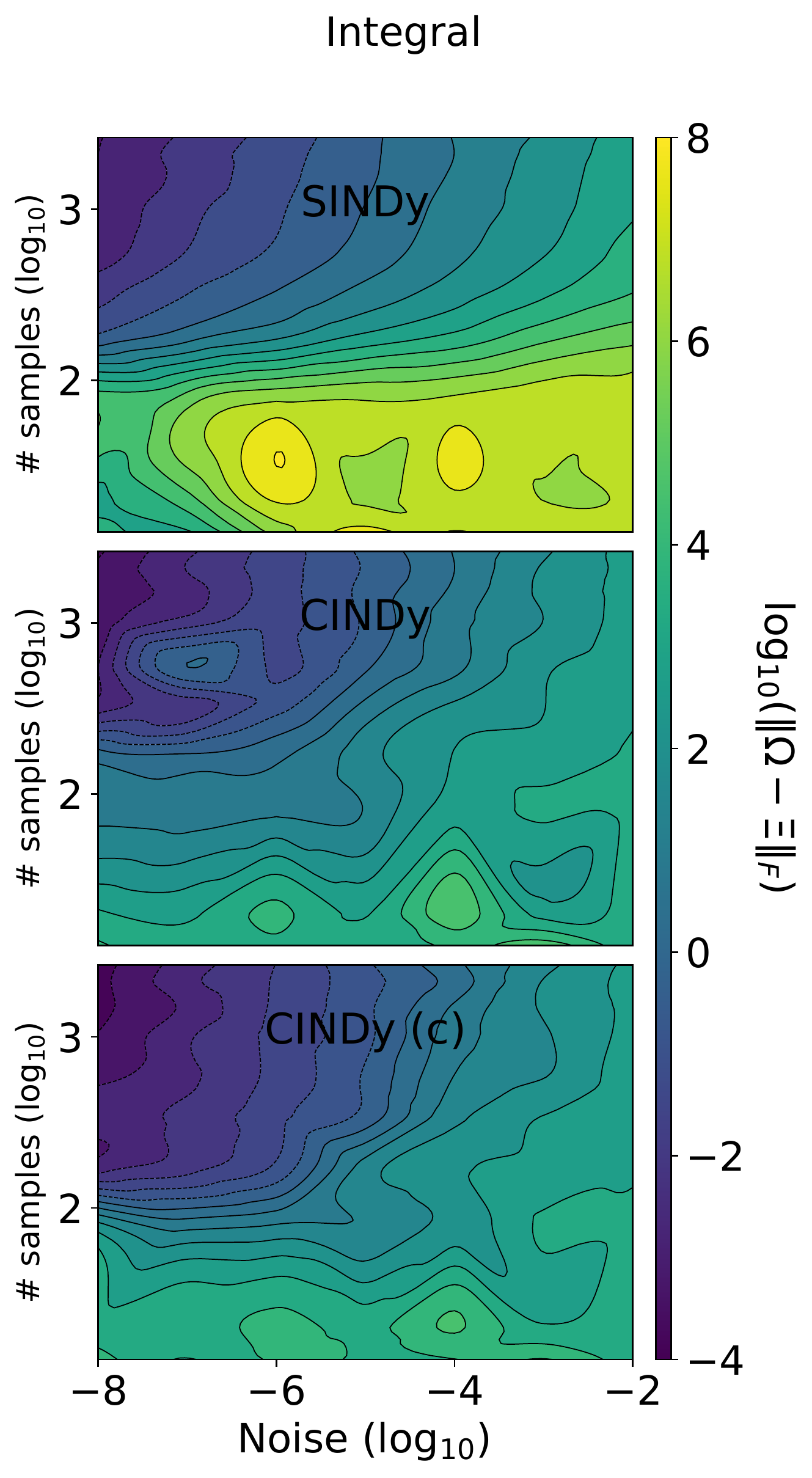} }\label{fig:appx:FPUT:int:dim5_sampleeff_acc}}%
    \hspace{\fill}
    \caption{\textbf{Sample efficiency of the sparse recovery of the
        Fermi-Pasta-Ulam-Tsingou model: } Algorithm comparison in terms of $\mathcal{E}_R$ for a Fermi-Pasta-Ulam-Tsingou model of
      dimension $d=5$ for the differential formulation (left) and the integral formulation (right).}%
    \label{fig:appx:FPUT_sampleeff}%
  \end{figure*}

\begin{figure*}[]
    \centering
    \vspace{-10pt}
    \hspace{\fill}
    \subfloat{{\includegraphics[width=8cm]{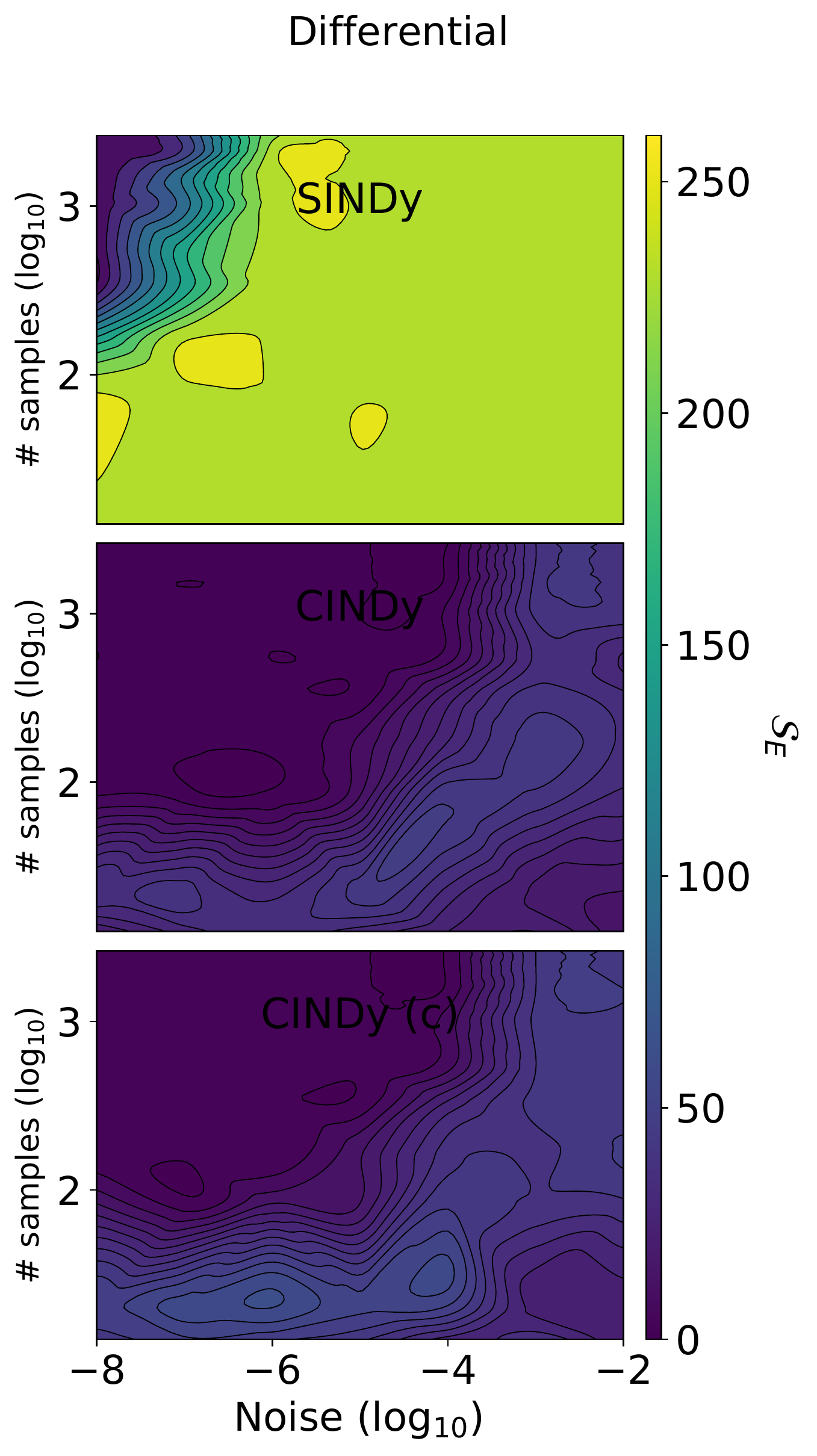} }\label{fig:appx:FPUT:diff:dim5_sampleeff_extra}}%
    %\qquad
    \hspace{\fill}
    \subfloat{{\includegraphics[width=8cm]{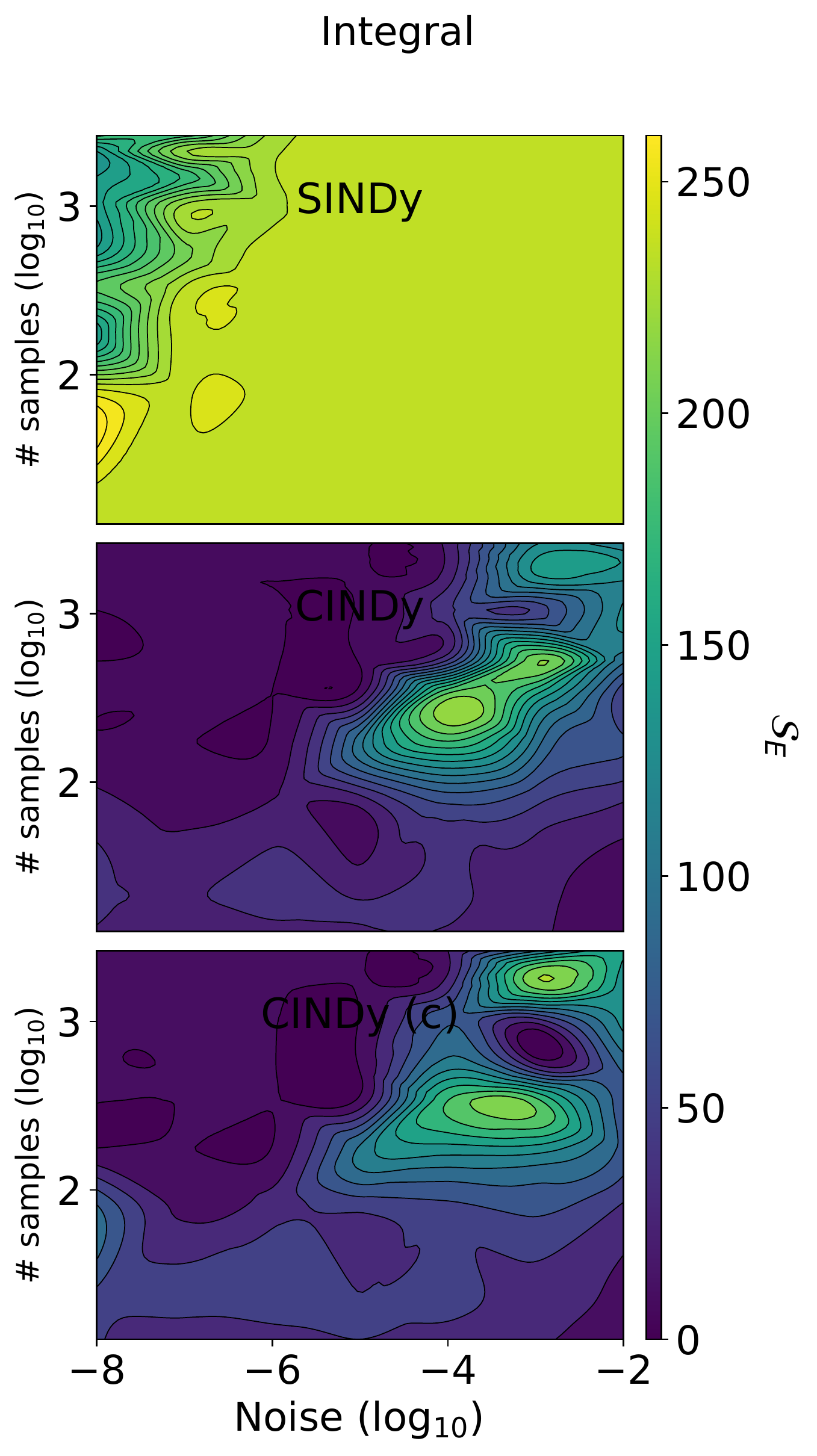} }\label{fig:appx:FPUT:int:dim5_sampleeff_extra}}%
    \hspace{\fill}
    \caption{\textbf{Sample efficiency of the sparse recovery of the
        Fermi-Pasta-Ulam-Tsingou model: } Algorithm comparison in terms of $\mathcal{S}_E$ for a Fermi-Pasta-Ulam-Tsingou model of
      dimension $d=5$ for the differential formulation (left) and the integral formulation (right).}%
    \label{fig:appx:FPUT_sampleeff_extra}%
  \end{figure*}

\begin{figure*}[]
    \centering
    \vspace{-10pt}
    \hspace{\fill}
    \subfloat{{\includegraphics[width=8cm]{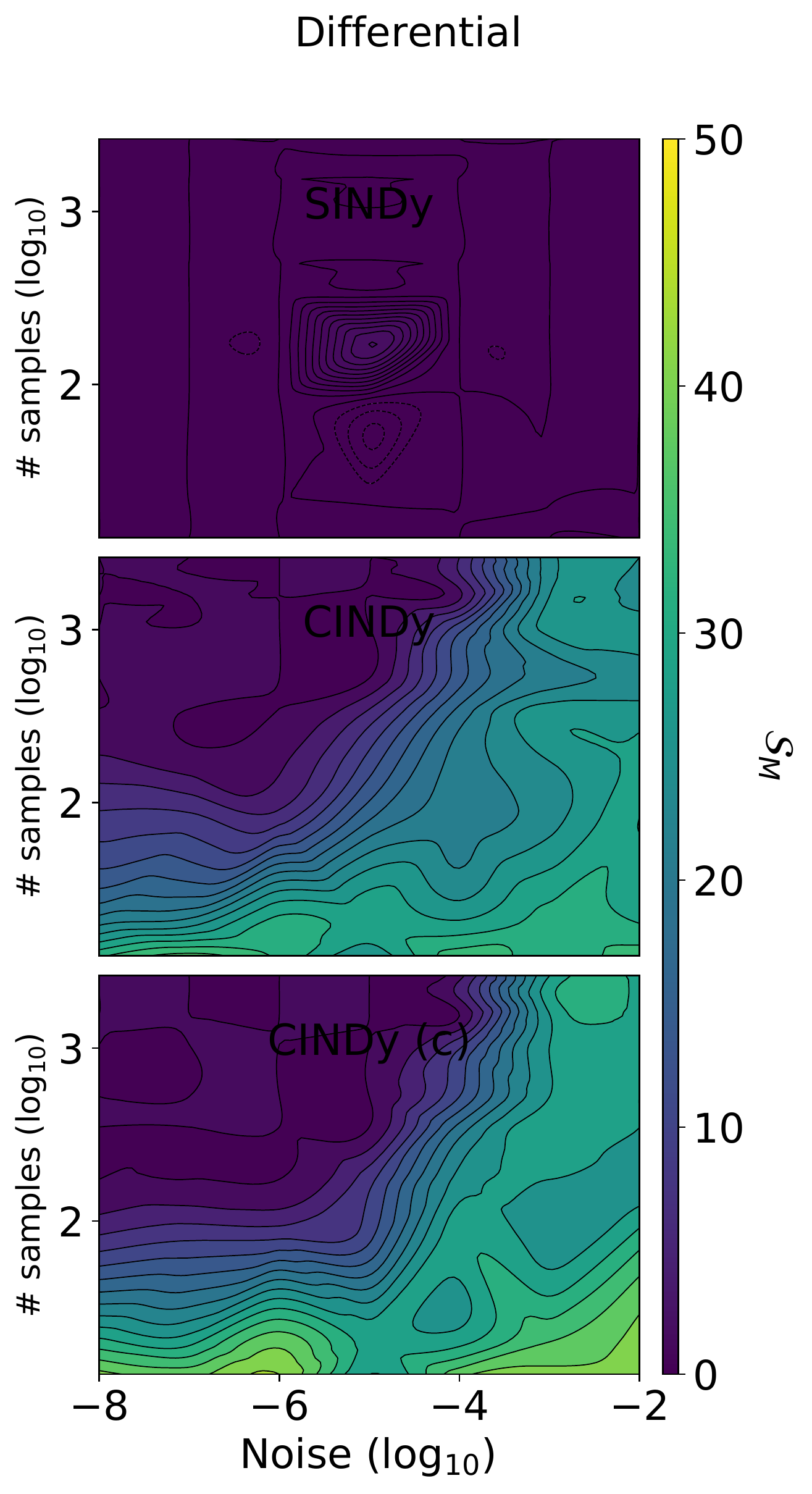} }\label{fig:appx:FPUT:diff:dim5_sampleeff_missing}}%
    %\qquad
    \hspace{\fill}
    \subfloat{{\includegraphics[width=8cm]{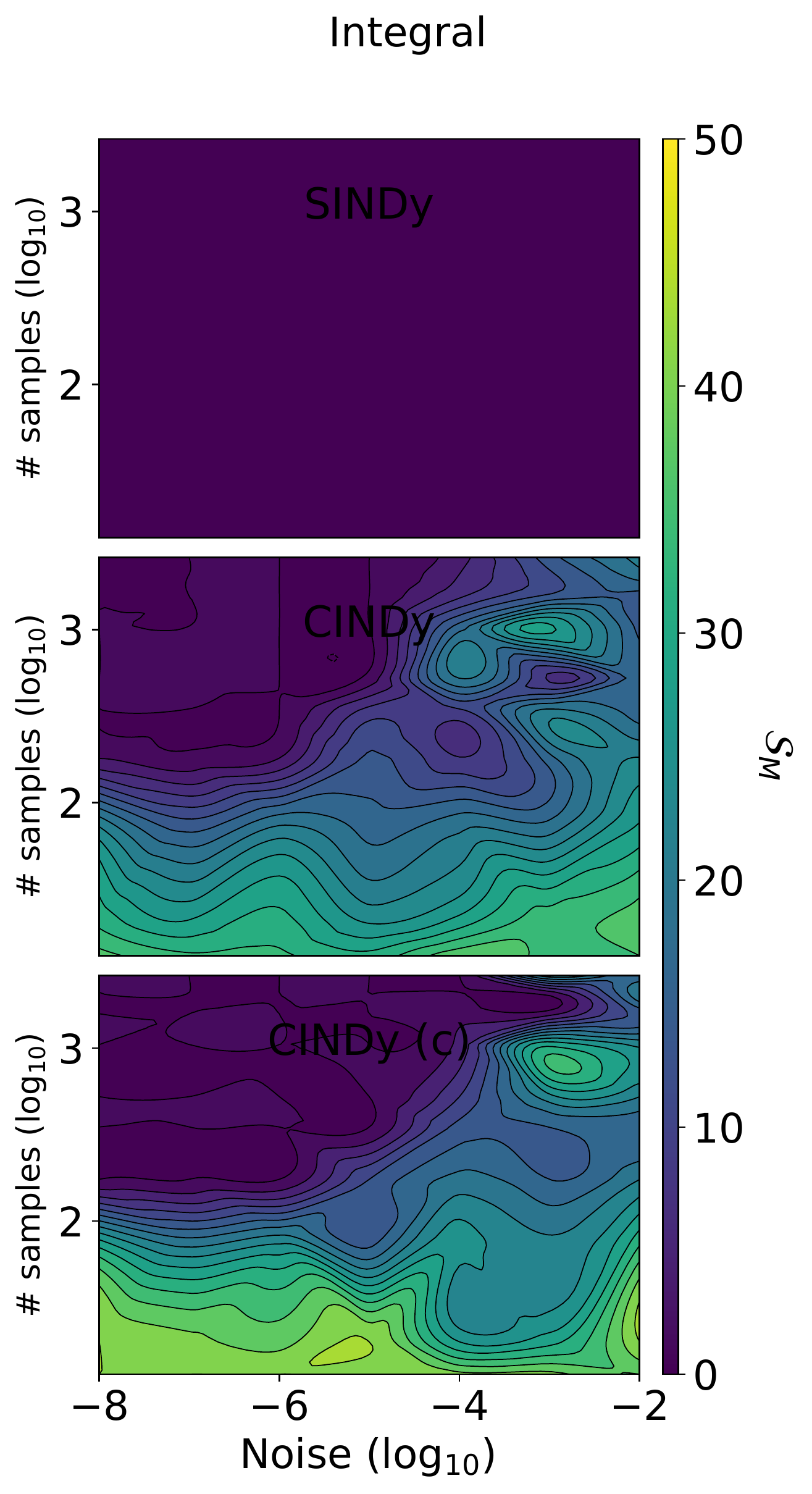} }\label{fig:appx:FPUT:int:dim5_sampleeff_missing}}%
    \hspace{\fill}
    \caption{\textbf{Sample efficiency of the sparse recovery of the
        Fermi-Pasta-Ulam-Tsingou model: } Algorithm comparison in terms of $\mathcal{S}_M$ for a Fermi-Pasta-Ulam-Tsingou model of
      dimension $d=5$ for the differential formulation (left) and the integral formulation (right).}%
    \label{fig:appx:FPUT_sampleeff_missing}%
  \end{figure*}

The images shown in Figure~\ref{fig:appx:FPUT_sampleeff}, \ref{fig:appx:FPUT_sampleeff_extra}, and \ref{fig:appx:FPUT_sampleeff_missing} show the evolution of $\mathcal{E}_R$, $\mathcal{S}_E$ and $\mathcal{S}_M$ as we vary the number of training data points when learning the Fermi-Pasta-Ulam-Tsingou dynamic ($d = 5$) with the dictionary described in Section~\ref{section:FPUT}. The images show the resulting metrics when generating $30$ data points per experiment and using local polynomial interpolation of degree $8$ to compute the derivatives and the integrals.

\end{document}